\documentclass[10pt,draft]{report}   
\usepackage{setspace}
\usepackage{buthesis}
\usepackage[bf,uppercase,center]{titlesec}

\titleformat{\chapter}[hang]{\Huge\bfseries\filcenter \singlespacing }{\thechapter.}{0 pt}{ \Huge\bfseries       \uppercase }
\titleformat{\section}[hang]{\large\bfseries }{\thesection.}{2 pt}{ \large\bfseries  }
\titlespacing{\chapter}{0pt}{0 in}{*1.5}

\usepackage{amsmath,amssymb,graphicx,amsfonts,amsthm}
\usepackage{verbatim}
\usepackage{bbm}
\usepackage{epsfig}
\usepackage{wrapfig}
\usepackage{latexsym,amsfonts,amscd}
\usepackage[usenames,dvipsnames]{xcolor}

\newcommand{\mcV}{}

\newcommand{\resie}[1]{\stackrel[#1]{}{\mbox{\emph{Res}} }}
\usepackage{changebar}
\usepackage{enumerate}
\usepackage{tikz}
\usetikzlibrary{decorations}
\usetikzlibrary{shadows}
\usepackage{latexsym}
\usepackage{mathrsfs}
\usepackage{amsfonts}
\usepackage{amsmath}
\usepackage{amssymb}
\usepackage{bm}
\usepackage{yfonts}
\usepackage{stackrel}

\newcommand{\mc}{\mathcal}
\newcommand{\vep}{\varepsilon}
\newcommand{\ep}{\epsilon}

\newcommand{\lt}{\left}
\newcommand{\rt}{\right}
\DeclareMathOperator{\Res}{Res}
\newcommand{\resi}[1]{\stackrel[#1]{}{\Res }}

\newcommand{\pd}{\partial}

\newcommand{\vho}{\varrho}

\newcommand{\mf}{\mathfrak}
\newcommand{\re}{\mf{Re}\,}
\newcommand{\im}{\mf{Im}\,}

\newcommand{\ol}{\overline}

\newcommand{\qtr}{\frac{1}{4}}

\newcommand{\mt}{\mathbb}

\newcommand{\wt}{\widetilde}

\newcommand{\LA}{\langle}
\newcommand{\res}{\mbox{Res}}
\newcommand{\RA}{\rangle} 

\newcommand{\bal}{\begin{align}}
\newcommand{\eam}{\end{align}}
\newcommand{\Or}{\mc{O}}

\numberwithin{equation}{section}
\newcommand{\bpm}{\begin{pmatrix}}
\newcommand{\epm}{\end{pmatrix}}
\newcommand{\bbm}{\lt| \begin{matrix}}
\newcommand{\ebm}{\end{matrix} \rt|}
\newcommand{\beq}{\begin{equation}}
\newcommand{\eeq}{\end{equation}}
\newcommand{\hf}{\frac{1}{2}}
\newcommand{\thf}{\tfrac{1}{2}}
\newcommand{\pif}{\frac{\pi}{2}}
\newcommand{\eal}{\end{align}}

\newcommand{\kt}{\frac{k}{2}}
\newcommand{\tkt}{\tfrac{k}{2}}

\newtheorem{theorem}{Theorem}[section]
\newtheorem{lemma}[theorem]{Lemma}
\newtheorem{proposition}[theorem]{Proposition}
\newtheorem{corollary}[theorem]{Corollary}
\newtheorem{remark}[theorem]{Remark}
\newcommand{\lf}{\lfloor}
\newcommand{\rf}{\rfloor}

\usepackage{lipsum}

\usepackage[footnotesize,bf]{caption}  

\usepackage[titles]{tocloft}

\dissertation 

\title{  \huge Triple Shifted Sums of \\ Automorphic L-functions }
\author{Thomas A. Hulse}
\degree{Doctor of Philosophy}
\department{the Department of Mathematics}
\previousdegrees{
		B.A., Colby College; Waterville, ME, 2007\\ 
		Sc.M., Brown University; Providence, RI, 2009}	
\thesismonth{May} \thesisyear{2013}
\begin{document}

\doublespacing
\begin{preliminaries}
\maketitle

\copyrightpage

\begin{signature}
  \director{Jeffrey Hoffstein, Ph.D., Advisor}
  \reader{Min Lee, Ph.D., Reader}
  \reader{Adrian Diaconu, Ph.D., Reader}
\end{signature}

\begin{vita}
Thomas Hulse was born in Hanover, New Hampshire on February 18th, 1985 to Muriel Anne Hulse (n\'{e}e $\mbox{St.}$ $\mbox{Sauveur}$) and James Edward Hulse. He received his high school diploma from Lebanon High School, New Hampshire in 2003, his B.A. in Mathematics and Physics from Colby College in 2007, and his Sc.M. in Mathematics from Brown University in 2009. During his time at Brown, Thomas taught several courses in calculus and one class in honors linear algebra and will be a recipient of a departmental teaching award upon graduation. He was the coordinator for the Math Resource Center, a walk-in help center for students, from 2009-2012 and worked as a Teaching Consultant for the Sheridan Center for Teaching and Learning from 2012-2013 as part of completing Certificate IV.

He also orchestrated weekly graduate student and postdoctoral lunches at Kabob and Curry.

\end{vita}

\begin{acknowledgments}

\noindent To The Automorphic $N$, thank you for making $N \geq 2$, because the $N=1$ case was uninteresting. 

\noindent To my advisor Jeff, thank you for supporting me to the end. 

\noindent To my committee members, Adrian and Min, thank you for willingly being required to read this.

\noindent To my parents and grandparents, thank you for helping me get to graduate school.

\noindent  To my entire family, thank you for taking for granted that I would finish graduate school.

\noindent  To the friends I made while at Brown, thank you for making me not want to leave graduate school. 

\noindent  To all my friends, thank you for making me who I am.

\noindent  To my cat Yitze, thank you for standing by me (or sitting on top of me) as I wrote this thesis (and much of this sentence), even when you were a distraction and I wanted you to go away.

\noindent  To my future wife Joanna, thanks for putting up with me. I love you.

\end{acknowledgments}


\begin{spacing}{1}
  \tableofcontents
  \clearpage{\pagestyle{empty}\cleardoublepage}

  \footnotesize
  \fontsize{11.5pt}{12.5pt}\selectfont
  \clearpage{\pagestyle{empty}\cleardoublepage}

  \clearpage{\pagestyle{empty}\cleardoublepage}
  \normalsize
\end{spacing}

\end{preliminaries}

\pagestyle{myheadings}

\chapter{Introduction}

\section{Shifted Sums and Subconvexity}
In analytic number theory, a problem that has generated a great deal of activity is the Generalized Lindel\"{o}f Hypothesis (GLH). In his originating statement \cite{Lindelof} concerning the Riemann zeta function, Lindel\"{o}f conjectured that on the critical line $s=\hf+it$ for $t \in \mt{R}$, in the order estimate of $\zeta(s)$ 
\beq
\zeta(\thf+it) \ll_\epsilon (1+|t|)^{\alpha+\epsilon},
\eeq
we can let $\alpha=0$. This statement is a corollary of the Riemann Hypothesis, and though it appears to be a strictly weaker statement, it remains stubbornly unproven. Lindel\"{o}f was able to use the functional equation of $\zeta(s)$, Stirling's approximation and the Phragm\'{e}n-Lindel\"{o}f convexity principle to show that $\alpha \leq \frac{1}{4}$. Over the past century, gradual progress has been made chipping away at $\alpha$, with the current record due to Huxley \cite{Huxley} standing at $\alpha \leq 32/205$.

The wider, also largely unproven statement of the GLH  concerns similarly characterized growth on the critical line of global $L$-functions. Consider $L(s,f)$,  an $L$-function of degree $d\geq 1$, where $f$ refers to some arithmetic object such as a number field or an automorphic form. In particular, when $\re s >1$ we take $L(s,f)$ to be expressible by an absolutely convergent Dirichlet series and
an Euler product. There is also a meromorphic continuation to all $s \in \mt{C}$ where $L(s,f)$ satisfies a functional equation of the form
\beq
 Q^\frac{s}{2} \lt[ \prod_{j=1}^d \pi^{-\frac{s}{2}}\Gamma\lt(\frac{s+\mu_j}{2}\rt)\rt]L(s,f)=\Lambda(s,f)=\vep(f) \Lambda(\ol{f},1-s).
\eeq
Here $\vep(f)$ is the root number with $|\vep(f)|=1$, the $\mu_j$ are the local parameters at infinity with $\re \mu_j >-1$ and $Q \in \mt{Z}_{\geq 1}$ is the conductor of $L(s,f)$. The object $\wt{f}$ is the dual of $f$, with different (perhaps conjugate) Dirichlet series coefficients, and all other parameters are unchanged excepting that $\vep(\wt{f})=\ol{\vep(f)}$. The GLH, which is itself a corollary of the Generalized Riemann Hypothesis \cite{Conrey}, conjectures that on the critical line in the order estimate,
\beq
L(\thf+it,f) \ll_{\epsilon,d} \lt(Q\prod_{j=1}^d (1+|t+u_j|) \rt)^{\alpha+\epsilon},
\eeq
we can take $\alpha =0$.  As with the zeta function, we can use the Phragm\'{e}n-Lindel\"{o}f convexity principle to establish that $\alpha\leq 1/4$. As $f$ varies over some family, any order estimate that improves upon this bound in any of the aspects of the conductor of $f$ is referred to as a subconvexity bound. 

This thesis is interested in the case where $L(s,f)$ comes from classical automorphic forms. When $f$ is an even weight $k> 0$ holomorphic cusp form on $\Gamma_0(N)\backslash \mt{H}$, given by the Fourier expansion
\beq
f(z) =\sum_{m=1}^\infty a(m) e^{2\pi i mz}=\sum_{m=1}^\infty A(m)m^{\frac{k-1}{2}}e^{2\pi i mz} 
\eeq
then $L(s,f)$ is a degree two $L$-function with the Dirichlet series expansion
\beq
L(s,f)=\sum_{m=1}^\infty \frac{a(m)}{m^{s+\frac{k-1}{2}}}=\sum_{m=1}^\infty \frac{A(m)}{m^{s}},
\eeq
that has conductor $N$ and local parameters determined by $k$. 

Improved knowledge of bounds of central values, that is bounds on $L(\hf,f)$, of these $L$-functions and other related higher-degree $L$-functions, is also interesting independent of GLH. Waldspurger's Theorem \cite{Waldspurger}, made explicit in the case of classical automorphic forms by Kohnen and Zagier \cite{KZ}, ties Fourier coefficients of a half-integral weight holomorphic Hecke cusp eigenform, $\mathfrak{f}$, to the central values of the $L$-functions of real primitive character twists of $f$, the integral-weight Shimura lift of $\mathfrak{f}$. Notably, the convexity bound on the conductor aspect of these central values gives the trivial bound on the growth of the corresponding Fourier coefficients, and improved subconvexity bounds of this type are on equal footing with progress toward the Ramanujan-Petersson Conjecture on related half-integral weight forms.


More recently, subconvexity bounds of higher degree $L$-functions of automorphic forms have been closely tied to progress toward the Quantum Unique Ergodicity (QUE) conjecture on arithmetic surfaces \cite{QUE}. Indeed, citing Watson's formula \cite{Watson} relating integrals of triple products of automorphic forms to central values of degree-eight automorphic $L$-functions, Sarnak demonstrated \cite{Sarnak} that a subconvexity estimate on these $L$-functions would in turn yield this conjecture. He then proceeded to use subconvexity bounds on Rankin-Selberg $L$-functions in the weight aspect to prove QUE in the case of ``dihedral forms''. Ultimately, the remainder of QUE on arithmetic surfaces was proven in a joint work by Holowinsky and Soundararajan \cite{HS}, each proving conditional results that excluded any exceptions to the main conjecture. Soundararajan \cite{Sound} proved a ``weak'' subconvexity result for a broad family of $L$-functions, in which he managed to save a logarithm term from the general convexity bound. Holowinsky's contribution \cite{Holowinsky2} took a wholly different route and resulted from using sieve methods to bound shifted convolution sums of multiplicative functions.   

The term \emph{shifted convolution sum} generally refers to an object resembling 
\beq
\sum_{m=1}^\infty \lambda_1(m+h)\lambda_2(m)W(m,h),
\eeq
where $W(m,h)$ is a weight function that, ideally, has manageable coupling between $m$ and $h$ and decays reasonably over specified ranges. These sums are sometimes referred to as ``sums of the additive divisor problem type'' \cite{Jutila}, as the case when the $\lambda_i$s are additive divisor functions has been long studied in classical analytic number theory. They were notably used by Atkinson \cite{Atkinson} to compute asymptotics of the fourth moment of the Riemann zeta function. Other sums of this form, where the $\lambda_i$s arose from automorphic forms, were first constructed by Selberg \cite{Selberg2} where they resembled
\beq
\sum_{n=1}^\infty \frac{a(n)\ol{b(n+h)}}{(2n+h)^s},
\eeq
when $\re s >1$ and $a(n)$ and $b(n)$ the Fourier coefficients of holomorphic cusp forms. Selberg studied these by means of replacing the real-analytic Eisenstein series in the Rankin-Selberg convolution with a Poincar\'{e} series and then untiling as in the convolution. 

Since then, shifted automorphic Dirichlet series have been frequently used as a means to derive subconvexity bounds of $L$-functions. To get a sense of why this is so, consider the approximate functional equation of the critical value of an $L$-function
\beq
L(\thf,f) = \sum_n \frac{A(n)}{\sqrt{n}}V_f\lt( \frac{n}{X}\rt), 
\eeq
where $V_f(y)$ is a smooth function that decays as $y$ gets large. If we are to consider $L(\thf,f)^2$ in this context, we get
\begin{align}
L(\thf,f)^2 &= \sum_{n,m =1}^\infty  \frac{A(n)A(m)}{\sqrt{nm}}V_f\lt( \frac{n}{X}\rt)  V_f\lt( \frac{m}{X}\rt) \\
 & = \sum_{n=1}^\infty \frac{A(n)^2}{n} V_f\lt( \frac{n}{X}\rt)^2 +2 \sum_{n,h=1}^\infty \frac{A(n+h)A(n)}{\sqrt{n+h}\sqrt{n}}V_f\lt( \frac{n+h}{X}\rt)V_f\lt( \frac{n}{X}\rt) \notag
\end{align}
and so we have a ``diagonal'' term'' that appears to be treatable via the Rankin-Selberg convolution and ``off-diagonal'' terms, which are shifted convolution sums. 

From the above heuristic we get a sense that information about shifted sums, particularly their estimated size, can be parlayed into information about the $L$-function. Indeed, various techniques for bounding shifted convolutions and then producing improved subconvexity results for Rankin-Selberg $L$-functions are nicely catalogued in a survey paper by Lau, Liu and Ye \cite{LLY}. Beyond that, such estimates have also been used to derive information about the growth of the Fourier coefficients themselves, as Good did \cite{Good}, and to investigate their non-trivial cancellation properties, as Goldfeld did \cite{Goldfeld}.  

We are also capable of constructing other shifted series of classical number theoretic interest. By modifying Hoffstein's Poincar\'{e} series, which will be described in the next section, we can get a continuation of
\beq
\sum_{a,c=1}^\infty \frac{\tau(4ac+h)}{a^{s+v}c^s}=\sum_{m=1}^\infty \frac{\tau(4m+h)\sigma_{-v}(m)}{m^s} \label{biggesty}
\eeq
where $\tau(n)$ is the $n$-th Fourier coefficient of the $\theta$-function and $\sigma_{-v}$ is the sum of positive divisors function. The author and collaborators used knowledge about the poles of the above function to asymptotically count the number of square discriminants, $h=b^2-4ac$, where $|a|,|b|,|c|<X$ as $X$ gets large. This is detailed in a work \cite{big}, currently in preparation, in which an earlier result due to Oh and Shah \cite{Oh} is substantially sharpened. Other work, such as Luo's more recent papers on shifted sums of half-integral weight forms and theta functions with cusp forms \cite{Luo1,Luo2}, respectively, further demonstrates that shifted convolution sums have become objects of general research interest in and of themselves. 

\section{Triple Shifted Sums}
This work is primarily dedicated to understanding the shifted convolution series of the form
\beq
T_{f_1,f_2,f_3}^\pm (s_1,s_2,s_3)=\sum_{m,h,n \geq 1}^\infty \frac{a_1(m-h)\ol{a_2(m)}a_3(h\pm n)}{m^{s_1+\frac{3}{2}k-1}h^{s_2}n^{s_3+\frac{k-1}{2}}} \label{refy2}
\eeq
where the $a_i$s are Fourier coefficients of even weight $k>0$ holomorphic forms, $f_1,f_2,f_3,$ for $\Gamma=SL_2(\mt{Z})$. It will be shown that the above series have meromorphic continuations to all of $\mt{C}^3$ and we will classify all of its polar lines. We achieve this by twice decomposing these sums as spectral expansions and, in the case of $T^+$, using Bochner's convexity theorem \cite{Bochner} to force a meromorphic continuation where we lack absolute convergence. We then use the information about the poles of these functions and the inverse Mellin transform to produce the main theorem of this work, which is the first result of which we are aware that uses these triple sums:
\newtheorem*{main}{Theorem \ref{main}}
\begin{main}
When $A_i(r)$ are the respective normalized coefficients of even, positive weight $k$ holomorphic cusp forms $f_i$, we have that for $X \gg 1$,
\bal
&\mkern-18mu \sum_{m,h,n\geq 1}^\infty \frac{A_1(m-h)\ol{A_2(m)}A_3(h\pm n)(1-\frac{h}{m})^{\frac{k}{2}}(1\pm \frac{n}{h})^{\frac{k}{2}}}{\sqrt{(m-h)(m)(h\pm n)} } e^{-(\frac{m+h+n}{X})}   \label{booyahz}
\\& \mkern300mu =T^{\pm}_{f_1,f_2,f_3}(1-\tkt,\tkt,\thf-\tkt)+\mc{O}_{f_1,f_2,f_3}(X^{-\hf+\vep}) \notag
\end{align}
where $\vep>0$ is arbitrarily small.
\end{main}

This result demonstrates substantial nontrivial cancellation of these shifted Fourier coefficients. In particular, we would like to obtain a similar asymptotic result for the uncoupled sum
\beq
\sum_{m,h,n \geq 1}^\infty \frac{A_1(m-h)\ol{A_2(m)}A_3(h\pm n)}{\sqrt{(m-h)(m)(h \pm n)}}e^{-\lt(\frac{m+h+n}{X}\rt)} \label{fail}
\eeq
though this has proven difficult. While it is possible to remove $(1+\frac{n}{h})^\kt$ by means of the integral formula for $(1+x)^{-\beta}$, convergence issues obstruct our ability to remove $(1-\frac{h}{m})^\kt$ and $(1-\frac{n}{h})^\kt$ at present.

Triple sums of this form \eqref{booyahs} and \eqref{fail} are of interest as, when $A_1(n)=A_2(n)=A_3(n) \in \mt{R}$ they seem to correspond to ``off-diagonal'' contributions of $L(\hf,f)^3$. Indeed, referring again to the approximate functional equation as a heuristic we see that, letting $W_f(m_1,m_2,m_3;X):=V_f\lt( \tfrac{m_1}{X}\rt)V_f\lt( \tfrac{m_2}{X}\rt)V_f\lt( \tfrac{m_3}{X}\rt)$, we have
\begin{align}
L&(\thf,f)^3=\sum_{m_i=1}^\infty  \frac{A(m_1)A(m_2)A(m_3)}{\sqrt{m_1m_2m_3}}W_f(m_1,m_2,m_3;X) \\
=&  2\sum_{m,h,n \geq 1}  \frac{A(m-h)A(m)A(h-n)}{\sqrt{(m-h)(m)(h-n)}}W_f(m-h,m,h-n;X) \notag  \\
&+2\sum_{m,h,n \geq 1}  \frac{A(m-h)A(m)A(h+n)}{\sqrt{(m-h)(m)(h+n)}}W_f(m-h,m,h+n;X) \notag \\
&+2\sum_{m,h \geq 1}  \frac{A(m-h)A(m)A(h)}{\sqrt{(m-h)mh}}W_f(m-h,m,h;X) + \sum_{m,h \geq 1}  \frac{A(m)^2A(h)}{\sqrt{m^2h}}W_f(m,m,h;X). \notag
\end{align}
The first two sums correspond to $T^\pm$ as we described above. The third sum is a simpler construction, one step removed from the construction of $T^\pm$ that is omitted from this work, and the last corresponds to $L(s,f\otimes f)L(s,f)$, a well-understood object. The ultimate goal is to get a formula for the asymptotics of higher moments of $L$-functions, approximating objects like
$$
\lt(\frac{1}{2\pi i } \rt)^3 \iiint\limits_{(2)(2)(2)} L(s_1,f)L(s_2,f)L(s_3,f) \Gamma(s_1-\thf)\Gamma(s_2-\thf)\Gamma(s_3-\thf)X^{s_1+s_2+s_3-\frac{3}{2}} \ ds_1ds_2ds_3
$$
or
$$
\frac{1}{2\pi i } \int\limits_{(2)} L(s,f)^3\Gamma(s-\thf)X^{s-\hf} \ ds 
$$
by decomposing them into pieces manageable by $T^\pm$. Integrals such of these are used to derive approximate functional equations, thus asymptotic formulas could in turn produce subconvexity estimates. 

The method for constructing $T^-$ requires very little new technology, but seems to have not been investigated until now. By making use of the well-studied Poincar\'{e} series
$$
P_h(z,s):=\sum_{\gamma \in \Gamma_\infty \backslash \Gamma} (\im \gamma z)^s e^{2\pi i h \gamma z}
$$
for $h \geq 1$, which is square integrable for sufficiently large $\re s$, we can expand the Petersson inner product $\LA P_h(*,s),(\im *)^k \ol{f_1}f_2 \RA$ either by untiling as in the Rankin-Selberg convolution  or by taking the spectral expansion of $P_h(z,s)$. The equivalence of these two expansions gives 
 \begin{align} \label{bigbird}
D^-_{f_1,f_2}(s;h)&=\sum_{m}  \frac{a_1(m-h)\ol{a_2(m)}}{m^{s+k-1}}  \\
&= \frac{(4\pi)^{k}h^{\hf-s}}{\Gamma(s+k-1)\Gamma(s)}  \lt( \sum_\ell \ol{\lambda_\ell(h)\rho_\ell(1)}\Gamma(s-\thf+it_\ell)\Gamma(s-\thf-it_\ell) \ol{\langle y^k \ol{f_1}f_2,\mu_\ell  \rangle} \rt. \notag \\
 &  \ \ + \lt.\frac{1}{2\pi i} \int_{(0)} \frac{\sigma_{2z}(h)h^{-z}\Gamma(s-\thf+z)\Gamma(s-\thf-z)}{2\zeta^*(1-2z)\zeta^*(1+2z)} \ol{\langle y^k\ol{f_1}f_2, E^*(*,\thf+z)\rangle} \ dz \rt) \notag
 \end{align}
when $\re s > \hf$. The discrete part of the spectrum arises from a linear combination of $\mu_\ell$s, which are orthonormal weight zero Hecke Maass eigenforms of the Laplacian for $\Gamma$, with Fourier-Whittaker expansion
\bal
\mu_\ell(z)&=\sum_{|m|\neq 0} \rho_\ell(\tfrac{m}{|m|})\lambda_{\ell}(|m|)2\sqrt{y} K_{it_\ell}(2\pi |m| y) e^{2\pi im x} \label{maass1} \\
& =\sum_{|m|\neq 0} \rho_\ell(\tfrac{m}{|m|})\lambda_{\ell}(|m|)|m|^{-\hf}W_{0,it_\ell}(4\pi |m| y) e^{2\pi im x} \notag
\end{align}
with eigenvalue $\qtr+t_\ell^2$, taking $t_\ell>0$ without loss of generality, and $\lambda_\ell(1)=1$. The continuous part of the spectrum is due to $E^*(z,s)$, the completed real-ananlytic Eisenstein series for $\Gamma$. Though the above formula for $D^-_{f_1,f_2}(s;h)$ changes as we move past the lines $\re s = \hf -r$, picking up residual terms due to the continuous part, an application of Stirling's approximation, Watson's  triple product formula \cite{Watson}, and the convexity bound on degree-eight $L$-functions gives us an absolutely convergent meromorphic expression for $D^-_{f_1,f_2}(s;h)$ for all $s \in \mt{C}$. 
By successively applying the raising operator  $\kt$-times to $\mu_\ell$, we can get the even weight $k$ Maass form
\beq
\mu_{\ell,k}(z)=   \sum_{|n|\neq 0} \rho_{\ell,k}(\tfrac{n}{|n|})\lambda_{\ell}(|n|)|n|^{-\hf} W_{\frac{kn}{2|n|},it_\ell}(4\pi |n| y)e^{2\pi i n x},
\eeq
which we can use in place of $f_2$ in the construction of $D^-_{f_1,f_2}$ above to get a similar spectral expansion and meromorphic continuation of
\beq
D^-_{f_3,\ell}(s;n)=\sum_{h}  \frac{a_3(h-n)\ol{\lambda_\ell(h)}}{h^{s+\frac{k-1}{2}}}. \label{fozie}
\eeq
By using the completed weight $k$ Eisenstein series we can also do the same for the series
\beq
D^-_{f_3,u}(s;n)=\sum_{h}  \frac{a_3(h-n)\sigma_{2u}(h)}{h^{s+\frac{k-1}{2}+u}}. \label{gonzo}
\eeq
Making use of these three shifted Dirichlet series, it is not difficult to see how they are combined to get a double spectral expansion of $T^-$ with an explicit convergent expression for all $\mt{C}^3$, though some work must be done to clarify the locations of poles and resolving issues of convergence. 

The construction of $T^+$ is significantly less straight-forward, as it requires meromorphic continuations of sums of the form
\beq
D^+_{f_3,\ell}(s;n)=\sum_{h}  \frac{a_3(h+n)\ol{\lambda_\ell(h)}}{h^{s+\frac{k-1}{2}}}, \ \ D^+_{f_3,u}(s;n)=\sum_{h}  \frac{a_3(h+n)\sigma_{2u}(h)}{h^{s+\frac{k-1}{2}+u}}. \label{refy1}
\eeq
To get a similar spectral expansion-type expression for these shifted Dirichlet series, it would appear that the natural object to consider would be $P_{-h}(z,s)$ as we considered $P_h(z,s)$ above, but this object isn't square integrable and so we cannot simply take its spectral expansion. 

Jeffrey Hoffstein constructed a square-integrable approximation of $P_{-h}(z,s)$: $P_{h,Y}(z,s)$ for large $Y \gg 1$ and small $\delta > 0$.  His intention was to meromorphically continue
\beq
D^+_{f_1,f_2}(s;n) = \sum_{h}  \frac{a_1(m+h)\ol{a_2(m)}}{m^{s+k-1}},
\eeq
as part of the greater aim of using it to derive subconvexity bounds for $L(\hf,f,\chi)$ in the conductor aspect. He has been largely successful, though I was eventually brought on as a collaborator in that project \cite{Jeff} to sort out many analytic details closely related to ones that emerge in this work.



Ultimately, all of the $D^+$ series noted above have meromorphic continuations to all $\mt{C}$, but they only have explicit spectral-type expansions when $s$ is sufficiently negative. For all other $s$ where the spectral and Dirichlet series expansion fail to be convergent, $D^+$ is given by a convergent contour integral expansion. Lacking the spectral expansion for all $s$ that we had in the $D^-$ case, we must use Bochner's theorem \cite{Bochner} to give a meromorphic continuation of
\beq
Z^+_{f_3,\ell}(s,w):=\sum_{h,n=1}^\infty \frac{a_3(h+n)\ol{\lambda_\ell(h)}}{h^{s+\frac{k-1}{2}}n^{w+\frac{k-1}{2}}}= \sum_{n=1}^\infty \frac{D^+_{f_1,\ell}(n;s)}{n^{s+\frac{k-1}{2}}}
\eeq
into the convex hull of the region where $Z^+_{f_3,\ell}(s,w)$ has a convergent expression, which happens to be all of $\mt{C}^2$. A similar argument is made to give a meromorphic continuation of the corresponding $Z^+_{f_3,u}(s,w)$ function, and combining these $Z^+$ functions with the expansion of $D^-_{f_1,f_2}$, we see that we get a meromorphic continuation of $T^+$ to all $\mt{C}^3$. 

From there, to get the bounds of the truncated sums we noted in the main theorem, we only need consider the inverse Mellin transform:
\beq 
\lt( \frac{1}{2\pi i }\rt)^3 \mkern-18mu \iiint\limits_{(2)(\kt+\hf)(2)}  \mkern-9mu T^\pm_{f_1,f_2,f_3}(s_1,s_2,s_3) \Gamma(s_1+\tkt-1)\Gamma(s_2-\tkt)\Gamma(s_3+\tfrac{k-1}{2})X^{s_1+s_2+s_3+\kt-\frac{3}{2}} ds_1 ds_2 ds_3. 
\eeq
Other estimates can be derived by taking variants of these integrals. 

\section{Future Research and Applications}

As mentioned in the previous section, the ultimate goal is to use analysis of these shifted sums and related objects to obtain precise results for third moments of $GL(2)$ $L$-functions. Conrey and Iwaniec \cite{CI} were able give a bound, in terms of the conductor of the real primitive character $\chi$, of $L^3(\hf,f,\chi)$ averaged over a family of holomorphic cusp forms of fixed weight and level. This in turn produced a subconvexity result for $L(\hf,f\chi)$ in the conductor aspect. An interesting application of these triple sums would be to compute an explicit asymptotic estimate of these averages; to produce a main term and an error term with a power savings and thus derive an improved subconvexity result. 

In a related area, Hoffstein \cite{Jeff} generalized the construction of $D^+_{f_1,f_2}(s;h)$ and $Z^+_{f_1,f_2}(s,w)$ for $f_i$ on different levels with added $\ell_i$ parameters to meromorphically continue the multiple Dirichlet series
\beq
Z_Q(s,w)= \sum_{ h, m_2 \ge 1 \atop m_1\ell_1= m_2\ell_2 + {hQ}} 
{a(m_1) \bar b(m_2)
\over
 (\ell_2m_2)^{s+k-1}(hQ)^{w + (k-1)/2}}.
\eeq
He then used an amplification argument, based on the work of Blomer \cite{Blomer} which was itself based off of the work of Duke, Friedlander and Iwaniec \cite{DFI}, to give an upper bound of $L(\hf,f,\chi)$ in the conductor aspect in terms of truncated sums computed from $Z_Q(s,w)$. Using this work as a guide, it should be possible to similarly generalize $T^\pm$ to higher levels with amplification parameters. From there, the goal is to construct an analogous amplification argument to produce some other subconvexity bound from triple shifted sums such as $T^\pm$.
\beq
L(\thf,f_1) \sum_{g \in \mc{F}} |L(\thf,f_1\otimes g)|^2
\eeq
where $\mc{F}$ is a Hecke basis of holomorphic cusp forms of weight $k$. The Petersson trace formula combined with the approximate functional equations of the $L$-functions above will yield truncated sums that correspond nicely with those derived above from $T^\pm$. By adding a weight term to the Poincar\'{e} series in the construction of $D^+$ or continuing to make use of raising operators, it would not be a difficult thing to generalize our construction to allow the holomorphic cusp forms to be of different weights. From there it should be possible to achieve a subconvexity bound for Rankin-Selberg $L$-functions in the weight aspect.

  Along similar lines, proceeding analogously from \eqref{biggesty} it is possible to analyze shifted series of half-integral weight terms via small variants of Hoffstein's Poincar\'{e} series. From there it seems reachable to derive a meromorphic continuation of
\beq
\sum_{m=1}^\infty \frac{a(m+h)\tau_k(m)}{m^s}, 
\eeq
where $\tau_k(n)$ is the Fourier coefficient of the $k$-th power of $\theta$, and use this to improve upon the bounds given by Luo \cite{Luo1}.

\chapter{The Spectral Expansion of \lowercase{$\uppercase{D}_{f,\ell}(s;h)$}}

\section{The Truncated Poincar\'{e} Series} \label{trunk}
In this chapter we will give a meromorphic continuation of $D_{f,\ell}^+(s;h)$ as in \eqref{refy1} by means of a spectral expansion, which is itself the second spectral expansion of the triple shifted sum $T^+$ in \eqref{refy2}. In this section, we begin by deriving a spectral expansion and meromorphic continuation of an approximation of the shifted convolution sum by means of a truncated Poincar\'{e} series. 

Let $\mu_{\ell,k}(z)$ be a weight $k$ Maass form corresponding to $\mu_\ell(z)$ in \eqref{maass1} via the isometry given by successive application of raising and lowering operators, that is
\beq
\mu_{\ell,k}(z):=\lt( \prod_{k>k'\geq 0, 2|k'}R_{k'}\rt) \mu_\ell(z) \label{isometry}
\eeq
where
\beq
R_{k'}:=iy\frac{\pd}{\pd x} +y\frac{\pd}{\pd y} +\frac{k'}{2}
\eeq
which means that
  \beq
\mu_{\ell,k}(z)=   \sum_{|n|\neq 0} \rho_{\ell,k}(\tfrac{n}{|n|})\lambda_{\ell}(|n|)|n|^{-\hf} W_{\frac{kn}{2|n|},it_\ell}(4\pi |n| y)e^{2\pi i n x}. \label{maassk}
  \eeq
We recall that
  \beq
  \rho_{\ell,k}(-1)=\epsilon_{\ell,k} \frac{\Gamma(\thf+it_\ell+\tkt)}{\Gamma(\thf+it_\ell-\tkt)}\rho_{\ell,k}(1) \mbox{ \ \ \ and \ \ \ } \rho_{\ell,k}(1)=(-1)^{\kt}\rho_{\ell}(1) \label{thatthing}
  \eeq
  where $\epsilon_{\ell,k}=\pm 1$ depends on the parity of $k$ and on whether $\mu_\ell$ is even or odd. We observe that the ratio of gamma factors is real-valued as $t_\ell \in \mt{R} \cup i\mt{R}$. 
Indeed, Selberg's eigenvalue conjecture gives us that $t_\ell \in \mt{R}_{>0}$ in the case when $\Gamma=SL_2(\mt{Z})$. We also note that 
\beq
\rho_\ell(1) \sim e^{\pif|t_\ell|}\log(1+|t_\ell|) 
\eeq
as $t_\ell \to \infty$. 
 
Let $h \in \mt{N}$. For $z \in \mt{H}$ and $s \in \mt{C}$ with $\re s>1$, consider the Poincar\'{e} series:
\beq
P_{-h}(z;s): = \sum_{\gamma \in \Gamma_\infty \backslash \Gamma} (\im(\gamma z))^se^{-2\pi i h \gamma z},
\eeq
which is the usual real-analytic Poincar\'{e} series except for the sign of the exponent. For fixed $z$ this series is invariant under the action of $\gamma \in \Gamma$ and is locally uniformly convergent and thus analytic for $s$ in this region. However, as noted in the introduction, it grows exponentially in $y$, which complicates efforts to determine its spectral properties. For this reason, for any fixed $Y \gg 1$ and $\delta >0$ we use the modified series $P_{h,Y}(z;s;\delta)$ given by 
\beq
P_{h,Y}(z;s;\delta): = \sum_{\gamma \in \Gamma_\infty \backslash \Gamma} (\im \gamma z )^s e^{-2\pi i h \re \gamma z + (2 \pi h \im \gamma z)(1-\delta)}\psi_Y(\im \gamma z).
\label{poin1}
\eeq
where $\psi_Y$ is a characteristic function on $[-Y,Y]$. The spectral properties of this modified Poincar\'{e} series have been well studied by Hoffstein and the author \cite{Jeff}, and it will be used here to similar effect. 
We also define the series  $P^{(2)}_{h,Y}(z;s;\delta)$:
\bal
& P_{h,Y}^{(2)}(z;s;\delta): =\label{poin2} \\
&=\sum_{h > m \geq 1 \atop \gamma \in \Gamma_\infty \backslash \Gamma} (\im \gamma z)^{s+\kt} a(h-m)\ol{\rho_{\ell,k}(-1)\lambda_\ell(m)}m^{-\hf}W_{-\kt,it_\ell}(4\pi m \im  \gamma z) e^{2\pi(m-h \delta) \im \gamma z} \psi_Y(\im \gamma z). \notag
\end{align}

We see that $P_{h,Y}(z;s;\delta)$ consists of finitely many terms and is in $L^2(\Gamma \backslash \mt{H})$, the Hilbert space of functions $g: \mt{H} \rightarrow \mt{C}$ satisfying 
\beq 
g(\gamma z)=g(z) \ \ \mbox{and} \ \ \int_{\Gamma \backslash \mt{H}} |g(z)|^2 \frac{dx dy}{y^2} < \infty
\eeq
for all $s \in \mt{C}$. 
 Similarly, we see that $ P_{h,Y}^{(2)}(z;s;\delta) \in L^2(\Gamma \backslash \mt{H})$. The Petersson inner product of two functions $F,G$ in $L^2(\Gamma \backslash \mt{H})$ is defined by 
\[
\LA F,G \RA =\int \int_{\Gamma \backslash \mt{H}} F(z)\ol{G(z)} \frac{dxdy}{y^2}. 
\] 
 Let $V_{f,\ell}(z)$ denote the $L^2(\Gamma \backslash \mt{H})$ function $y^{k/2}\ol{f(z)}\mu_{\ell,k}(z)$. which is rapidly decreasing as $y \rightarrow \infty$. Our approach to studying \eqref{refy1} will be to compute $ \LA P_{h,Y}(*;s;\delta),V_{f,\ell} \RA- \LA  P_{h,Y}^{(2)}(*;s;\delta),1 \RA$ in different ways, and to then let $Y \rightarrow \infty$ and $\delta \to 0$. 

To begin we unfold the respective inner products:
\begin{subequations}\label{tag1}
\begin{align}
\mc{I}_{f,\ell, Y,\delta}(s;h)
&:= \LA P_{h,Y}(*;s;\delta),V_{f,\ell} \RA- \LA  P_{h,Y}^{(2)}(*;s;\delta),1 \RA \\
&=\iint\limits_{\Gamma \backslash \mt{H}} \lt( P_{h,Y}(z;s;\delta)f(z)\ol{\mu_{\ell,k}(z)} y^{k/2} -P_{h,Y}^{(2)}(z;s+k/2;\delta)\rt) \frac{dx dy}{y^2}  \notag \\
&= \sum_{n \geq 1} \sum_{|m|\neq 0} \frac{a(n)\ol{\rho_{\ell,k}(\tfrac{m}{|m|})\lambda_{\ell}(|m|)}}{|m|^{\hf}} \int_0^1 e^{2\pi i x(n-m-h)} \ dx \notag \\
& \ \ \ \times \int_{Y^{-1}}^Y e^{-2\pi y(n-h(1-\delta))} y^{s+\kt-1}W_{\frac{km}{2|m|},it_\ell}(4\pi|m|y)\frac{dy}{y} \notag \\
&\notag \ \ \ - \ol{\rho_{\ell,k}(-1)} \sum_{h>m\geq 1}  \frac{a(h-m)\ol{\lambda_\ell(m)}}{m^{\hf}}  \int_{Y^{-1}}^Y y^{s+\kt-1}W_{-\kt,it_\ell}(4\pi m y)e^{-2\pi y(-m+h \delta)} \frac{dy}{y} \\ 
&=\frac{\ol{\rho_{\ell,k}(1)}}{(2\pi)^{s+\kt-1}} \sum_{m> 0} \frac{a(m+h)\ol{\lambda_\ell(m)}}{m^{s+\frac{k-1}{2}}}  \int_{Y^{-1}2\pi m}^{Y2\pi m} e^{-y(1+h\delta/m)}W_{\kt,it_\ell}(2y)y^{s+\kt-1} \frac{dy}{y}. \label{intend}
\end{align}
\end{subequations}
Using well-known bounds for the Whittaker function, 
\beq
W_{\lambda,\gamma+it}(2y) \ll_{\lambda,\gamma,t}  e^{-y}y^\lambda \ \ \mbox{for } y \gg 1 
\label{bessel}
\eeq
\beq
W_{\lambda,\gamma+it}(2y) \ll_{\lambda,\gamma,t}  \lt\{
\begin{array}{ll} y^{\frac{1}{2}-|\gamma |} & \mbox{if } \gamma \neq 0 \\
y^{\frac{1}{2}}|\log y| & \mbox{if } \gamma = 0 \end{array} \rt. \ \ \mbox{for } y \ll 1,
\eeq
where $\lambda,\gamma,t \in \mt{R}$, we see that for $\re s >\hf-\kt+|\im t_\ell|$, and all $Y > 1$ and $\delta >0$, the integral in (\ref{intend}) converges absolutely and satisfies the upper bound
\beq
 \int_{Y^{-1}2\pi m}^{Y2\pi m} e^{-y(1+h\delta/m)}W_{\kt,it_\ell}(2y)y^{s+\kt-1} \frac{dy}{y} 
\ll_{k,t_\ell} \Gamma(\re s+k-1).
\eeq
Using the known formula 7.621(3) in \cite{GR}   for $\int_0^\infty t^{\alpha}e^{-st} W_{\lambda,\mu}(q t) \ dt$ when $\re (\alpha \pm \mu+\tfrac{3}{2} )> 0$, $\re s > -q/2$ and $q>0$, we have that the integral in (\ref{intend}) converges to 
\begin{align}
\label{side1}  \frac{\Gamma(s+(k-1)/2+it_\ell)\Gamma(s+(k-1)/2-it_\ell)}{2^{s+\kt-1}\Gamma(s)}
\end{align}
as $Y\rightarrow \infty$ and $\delta \rightarrow 0$. This gives us that 
\beq
\lim_{{Y\rightarrow \infty} \atop {\delta\rightarrow 0}} \mc{I}_{f,\ell, Y,\delta}(s;h) =
 \frac{\ol{\rho_{\ell,k}(1)}}{\mcV  (4\pi)^{s+\kt-1}}  \frac{\Gamma(s+\frac{k-1}{2}+it_\ell)\Gamma(s+\frac{k-1}{2}-it_\ell)}{\Gamma(s)} \sum_{m> 0} \frac{a(m+h)\ol{\lambda_\ell(m)}}{m^{s+\frac{k-1}{2}}}. 
\eeq 
Thus we obtain a shifted convolution sum of the Fourier coefficients of the holomorphic form and the Maass form. Next we intend to give an analytic continuation of this sum in terms of a spectral expansion of the Poincar\'{e} series we used to define it. 

Since $P_{h,Y}(z;s;\delta) \in L^2(\Gamma \backslash \mt{H})$, it has a spectral expansion, 
\begin{align}
P_{h,Y}(z;s;\delta) &= \label{side2} \\
&\sum_{j=1}^\infty \langle P_{h,Y}(*;s;\delta),u_j \rangle u_j(z)+\frac{1}{4\pi} \int_{-\infty}^\infty\langle P_{h,Y}(z;s;\delta), E(*,\thf+it) \rangle E(z,\thf+it) \ dt, \notag
\end{align}
 where $u_j \in L^2(\Gamma \backslash \mt{H})$ are a basis of weight zero Maass forms which are orthonormal with respect to the Petersson inner product, simultaneously diagonalized with respect to the Hecke operators and the Laplacian. The function $E(z,s)$ is the real-analytic Eisenstein series
\beq
E(z,s)=\sum_{\gamma \in \Gamma_\infty \backslash \Gamma} (\im(\gamma z))^s.
\eeq
\begin{remark}
Since every Maass form $u_j$ has two associated paramters, $\pm t_j$, and we will often find ourselves summing over these parameters instead of the Maass forms themselves, we extend the indices of these parameters to all $\mt{Z}_{\neq 0}$ by the convention that $t_{-j}=-t_j$. Since we are working in $SL_2(\mt{Z})$, $t_j \neq 0$ and so this convention presents no complications. 
\end{remark}
The spectral expansion of $P_{h,Y}$  yields an  expansion of $\mc{I}$:
 \begin{align}
\mc{I}_{f,\ell, Y,\delta}(s;h)
&=\LA P_{h,Y}(*;s;\delta),V_{f,\ell} \RA- \LA  P_{h,Y}^{(2)}(*;s;\delta),1\RA \label{spec2} \\
&=\sum_j \ol{\LA V_{f,\ell}, u_j \RA} \langle P_{h,Y}(*;s;\delta),u_j \rangle \notag \\
&\ +\frac{1}{4\pi} \int_{-\infty}^\infty\ol{\LA  V_{f,\ell},E(*,1/2+it)\RA}\langle P_{h,Y}(z;s;\delta), E(*,1/2+it) \rangle \ dt \notag \\
& \ -\LA  P_{h,Y}^{(2)}(*;s;\delta),1 \RA.  \notag
\end{align}
We then observe that \bal
\langle P_{h,Y}(*;s;\delta),u_j \rangle &=  \int_{Y^{-1}}^Y \int_0^1 y^{s+\hf} e^{-2\pi i h x} e^{2\pi h y(1-\delta) }\ol{u_j(z)} \frac{dxdy}{y^2} \notag \\ 
&=2\ol{\rho_j(-1)\lambda_j(h)}\int_{Y^{-1}}^Y y^{s-\hf} e^{2 \pi h y (1-\delta)}K_{it_j}(2 \pi h y) \frac{dy}{y} \notag \\
&=\frac{2\ol{\rho_j(-1)\lambda_j(h)}}{\mcV (2\pi h)^{s-\hf}}\int_{2\pi h Y^{-1}}^{2 \pi h Y} y^{s-\hf} e^{ y (1-\delta)}K_{it_j}(y) \frac{dy}{y} \label{oldint1},
\end{align} 
and similarly, using the Fourier expansion of $E(z,\thf+it)$ we get
\bal
\langle   P_{h,Y}(*;s;\delta),&E(*,\thf+it)\rangle=   \notag \\
& \frac{2h^{-it}\sigma_{2it}(h)}{\mcV \zeta^*(1-2it)(2\pi h)^{s-\hf}}\int_{2\pi h Y^{-1}}^{2 \pi h Y} y^{s-\hf} e^{y (1-\delta)}K_{it}( y) \frac{dy}{y}.  \label{oldint2}
\end{align}
The spectral expansion of $P_{h,Y}$ was well studied in \cite{Jeff} and the integrals in \eqref{oldint1} and \eqref{oldint2} were meromorphically continued in $s$ in the limit as $Y \rightarrow \infty$ and $\delta \to 0$. However  we must also similarly continue a generalization of this integral that appears in $\LA  P_{h,Y}^{(2)}(*;s;\delta),1 \RA$,
\begin{align}
& \LA  P_{h,Y}^{(2)}(*;s;\delta),1 \RA \notag\\
&\ \ \ = \frac{\ol{\rho_{\ell,k}(-1)}}{\mcV (2\pi)^{s+\kt-1}} \sum_{h>m\geq 1}  \frac{a(h-m)\ol{\lambda_\ell(m)}}{m^{s+(k-1)/2}}  \int_{2\pi m Y^{-1}}^{2\pi m Y} y^{s+\kt-1}W_{-\kt,it_\ell}(2 y)e^{y(1-\frac{h \delta}{m})} \frac{dy}{y}, \label{mhere}
\end{align}
and this is done in the next section.
\section{The $M$-functions}

Let
\beq
M_{Y,h,k}(s,z/i,\delta):=\int_{Y^{-1}2\pi h}^{Y2\pi h}  y^{s-1}e^{y(1-\delta)}W_{{\frac{k}{2}},z}(2y) \frac{dy}{y}
\eeq 
where $k \in \mt{R}, z\in \mt{C}$, $Y\gg 1$ and $\delta >0$. This is a generalization of the integral that appears in \eqref{mhere} as well as those that appear in \eqref{oldint1} and \eqref{oldint2} 
as $W_{0,z}(2y)=\lt(\frac{2y}{\pi}\rt)^{\hf}K_{z}(y)$. Now let
\beq
M_k(s,z/i,\delta):=\int_{0}^{\infty}  y^{s-1}e^{y(1-\delta)}W_{\frac{k}{2},z}(2y) \frac{dy}{y}.
\eeq
which we expect to be the limit of $M_{Y,k}$ as $Y \to \infty$. 

As has been noted, this is a generalization of a similarly defined integral in \cite{Jeff} and so our results should agree with that work if we let $k=0$. Since our short-term aim is to remove the dependence on $Y$ and $\delta$ in $M_{Y,k}(s,z/i,\delta)$, in this section we investigate how $M_k(s,z/i,\delta)$ and $M_{Y,k}(s,z/i,\delta)$ are related as $Y \rightarrow \infty$ and we give growth properties and residues of a meromorphic continuation of $M_k(s,z/i,\delta)$. 

\begin{remark}
This integral also arrises with some modification in \cite{big}, as well as elsewhere and so we let $k \in \mt{R}$, though that generality is not required for the purposes of this work.
\end{remark}

We begin the following proposition, which gives us the degree to which $M_{k,Y}$ is approximated by $M_k$ as $Y \to \infty$. 
\begin{proposition} \label{ybound}
For fixed $\vep>0, Y \gg 1, $ $1 > \delta >0$, and $A \in \mt{Z}_{\geq 0}$, we have that for $\re s>\hf+|\re z| +\vep$
\begin{align}
|M_{Y,h,k}(s,z/i,\delta)-M_k(s,z/i,\delta)| \ll \frac{e^{-Y2\pi h \delta}(Yh)^{\re s +\frac{k}{2}+A+\vep-2}}{\delta(1+|\im z|)^A}+\frac{(Y^{-1}h)^{\re s-\hf-|\re z|-\vep}}{(1+|\im z|)^A} \label{ygrowth}
\end{align} where the implied constant is dependent on $A$, $k$, $\re s$, $\re z$  and $\vep$. 
\end{proposition}
\begin{proof}To prove this we require the following lemma:
\begin{lemma}
For any integer $A \geq 0$ and arbitrarily small $\vep>0$, we have for $y \gg 1$ that
\beq 
W_{{\frac{k}{2}},z}(2y) \ll \frac{e^{-y}y^{ A  +\frac{k}{2}+\vep}}{(1+|\im z|)^A}. \label{1part}
\eeq
Similarly, when $y \ll 1$ we have 
\beq 
W_{{\frac{k}{2}},z}(2y) \ll \frac{y^{\hf-|\re z|}}{(1+|\im z|)^ A}. \label{2part}
\eeq
The implied constants for both bounds are independent of $y$ and $\im z$. 
\end{lemma}

\begin{proof} Let $\mu=w+it$ with $w,t \in \mt{R}$, and also let $\lambda \in \mt{R}$. It is well-known (see 9.227 in \cite{GR}) that when $\lambda$ and $\mu$ are fixed and $y\to \infty$ that 
\beq
W_{\lambda,\mu}(y) \ll y^\lambda e^{-\frac{y}{2}} \label{woot3}
\eeq
and (see 3.6.2 in \cite{GoHu}) if $y \to 0$ then
\beq
W_{\lambda,\mu}(y) \ll \lt\{ 
\begin{array}{ll}
y^{\hf-|w|} & \mbox{if } w\neq 0\\
y^{\hf}\log(y)& \mbox{if } w= 0.
\end{array} \label{woot4}
\rt.
\eeq
For the purposes of this lemma, we need to verify that comparable bounds can be made which are uniform in $t$. 

From 9.223 in \cite{GR}, we know that when $y>0$ we have that
\beq 
W_{\lambda,\mu}(y) = \frac{e^{-\frac{y}{2}}}{2\pi i }\int_{C} \frac{\Gamma(u-\lambda)\Gamma(-u-\mu+\hf)\Gamma(-u+\mu+\hf)}{\Gamma(-\lambda+\mu+\hf)\Gamma(-\lambda-\mu+\hf)}y^u \ du \label{above}
\eeq
where $C$ is a path from $-i\infty$ to $i\infty$ chosen in such a way that the poles of the function $\Gamma(u-\lambda)$ are to the left of $C$ and the poles of the functions $\Gamma(-u-\mu+\hf)$ and $\Gamma(-u+\mu+\hf)$ are to the right. We want to be sure that none of these poles overlap, so for now we will let $t \neq 0$. If we shift the line of integration of $u$ in \eqref{above} right  to $(\lambda+\vep)$ for some arbitrarily small $\vep>0$, such that the line does not intersect any poles of the integrand, then we pass over the simple poles at $u=r_1+\hf+ \mu$ where $0\leq r_1 \leq \lfloor \lambda-\hf+\vep- w \rfloor$ and at  $u=r_2+\hf- \mu$ where $0\leq r_2 \leq \lfloor \lambda-\hf+\vep+ w \rfloor$. This gives us that
\begin{align}
W_{\lambda,\mu}(y) =& \frac{e^{-\frac{y}{2}}}{2\pi i }\int_{(\lambda+\vep)} \frac{\Gamma(u-\lambda)\Gamma(-u-\mu+\hf)\Gamma(-u+\mu+\hf)}{\Gamma(-\lambda+\mu+\hf)\Gamma(-\lambda-\mu+\hf)}y^u \ du \label{above2} \\
& + \sum_{r=0}^{\lfloor \lambda-\hf+\vep-w \rfloor} e^{-\frac{y}{2}}y^{r+\hf+\mu} \frac{(-1)^r\Gamma(r+\hf+\mu-\lambda)\Gamma(-2\mu-r)}{r!\Gamma(-\lambda+\mu+\hf)\Gamma(-\lambda-\mu+\hf)} \notag \\ 
& +   \sum_{r=0}^{\lfloor \lambda-\hf+\vep+w \rfloor} e^{-\frac{y}{2}}y^{r+\hf-\mu} \frac{(-1)^r\Gamma(r+\hf-\mu-\lambda)\Gamma(2\mu-r)}{r!\Gamma(-\lambda+\mu+\hf)\Gamma(-\lambda-\mu+\hf)}. \notag
\end{align}
 From Stirling's approximation formula, we observe that if $y\geq 1$
\begin{align}
 & \sum_{r=0}^{\lfloor \lambda-\hf+\vep-w \rfloor}  e^{-\frac{y}{2}}y^{r+\hf+\mu} \frac{(-1)^r\Gamma(r+\hf+\mu-\lambda)\Gamma(-2\mu-r)}{r!\Gamma(-\lambda+\mu+\hf)\Gamma(-\lambda-\mu+\hf)} \\ 
 &\notag +   \sum_{r=0}^{\lfloor \lambda-\hf+\vep+w \rfloor} e^{-\frac{y}{2}}y^{r+\hf-\mu} \frac{(-1)^r\Gamma(r+\hf-\mu-\lambda)\Gamma(2\mu-r)}{r!\Gamma(-\lambda+\mu+\hf)\Gamma(-\lambda-\mu+\hf)} \ll_{\lambda,w,\vep} e^{-\frac{y}{2}}y^{\lambda+\vep} e^{-\frac{\pi}{2}|t|}(1+|t|)^{\lambda-\hf+|w|},
\end{align}
and if $y<1$
\begin{align}
 & \sum_{r=0}^{\lfloor \lambda-\hf+\vep-w \rfloor}  e^{-\frac{y}{2}}y^{r+\hf+\mu} \frac{(-1)^r\Gamma(r+\hf+\mu-\lambda)\Gamma(-2\mu-r)}{r!\Gamma(-\lambda+\mu+\hf)\Gamma(-\lambda-\mu+\hf)}  \\ 
 & \notag+  \mkern-9mu \sum_{r=0}^{\lfloor \lambda-\hf+\vep+w \rfloor} e^{-\frac{y}{2}}y^{r+\hf-\mu} \frac{(-1)^r\Gamma(r+\hf-\mu-\lambda)\Gamma(2\mu-r)}{r!\Gamma(-\lambda+\mu+\hf)\Gamma(-\lambda-\mu+\hf)} \ll_{\lambda,w,\vep} e^{-\frac{y}{2}}y^{\hf-|w|} e^{-\frac{\pi}{2}|t|}(1+|t|)^{\lambda-\hf+|w|}.
\end{align}
Again using Stirling's approximation, we get the bound 
\begin{align}
& \int_{(\lambda+\vep)} \frac{\Gamma(u-\lambda)\Gamma(-u-\mu+\hf)\Gamma(-u+\mu+\hf)}{\Gamma(-\lambda+\mu+\hf)\Gamma(-\lambda-\mu+\hf)}y^u \ du \label{woot1} \\
&  \ll_{\lambda,w,\vep} y^{\lambda+\vep} \int_{-\infty}^{\infty} (1+|v|)^{\vep-\hf}(1+|t-v|)^{-\vep}(1+|t+v|)^{-\vep}e^{-\frac{\pi}{2}(2\max(|t|,|v|)-2|t|+|v|)} \ dv \ll_{\lambda,w,\vep} y^{\lambda+\vep}. \notag
\end{align}
The case where $t=0$ is covered by \eqref{woot3} and \eqref{woot4}, and so we obtain the following bounds which are  uniform in $\im z$: if $y>1$ then
\beq
W_{\kt,z}(2y) \ll_{\kt,\re z,\vep} e^{-y}y^{\kt+\epsilon} \label{wgrow1}
\eeq
and if $y <1$ then 
\beq
W_{\kt,z}(2y) \ll_{\kt,\re z} \lt\{
\begin{array}{ll}
y^{\hf-|\re z|} & \mbox{if } w\neq 0\\
y^{\hf}\log(y)& \mbox{if } w= 0.
\end{array}  \label{wgrow2}
\rt.
\eeq
By adding the recurrence relations 9.234(1) and 9.234(2) in \cite{GR}:
\beq
W_{\lambda,\mu}(y) = \sqrt{y} W_{\lambda-\frac{1}{2},\mu-\hf}(y)+\lt(\hf-\lambda+\mu\rt)W_{\lambda-1,\mu}(y) 
\eeq
and
\beq
W_{\lambda,\mu}(y) = \sqrt{y} W_{\lambda-\frac{1}{2},\mu+\hf}(y)+\lt(\hf-\lambda-\mu\rt)W_{\lambda-1,\mu}(y),
\eeq
and then changing $\lambda \to \lambda +1$, we get the formula
\beq
W_{\lambda,\mu}(y)=\frac{\sqrt{y}}{2\mu}\lt(W_{\lambda+\hf,\mu+\hf}(y)- W_{\lambda+\hf,\mu-\hf}(y) \rt).
\label{rec1}
\eeq
When $y$ is large, combining  (\ref{wgrow1}) with the recursive relation (\ref{rec1}) $ A $ times on $W_{\kt,z}(2y)$ for some integer $A \geq 0$, we get \eqref{1part}.
Similarly when $y$ is small, using (\ref{rec1})  $A$ times with (\ref{wgrow2}) we get \eqref{2part} . \end{proof}
From the previous lemma, we can derive that for $Y$ sufficiently large
\begin{align}
\int_{Y2\pi h}^\infty y^{s-1} e^{y(1-\delta)}W_{\kt,z}(2y) \frac{dy}{y} \notag 
& \ll 
 \int_{Y2\pi h}^\infty y^{\re s -1}e^{y(1-\delta)}\frac{e^{-y}y^{A+\frac{k}{2}+\vep}}{(1+|\im z|)^A} \frac{dy}{y} \notag \\ 
 & \ll \frac{1}{(1+|\im z|)^A} \int_{Y2\pi h}^\infty y^{\re s+\kt+A+\vep-1}e^{-y\delta} \frac{dy}{y} \notag\\
 & =  \frac{1}{(1+|\im z|)^A\delta^{\re s + \kt +A+\vep -1}} \Gamma\lt[\re s +A+\vep+\tkt -1, Y2\pi h \delta\rt] \notag \\
 & \ll \frac{e^{-Y2\pi h \delta}(Yh)^{\re s+\kt+A+\vep-2}}{\delta(1+|\im z|)^A} \label{top1}
\end{align} 
using the fact that $\Gamma[r,x]:=\int_x^\infty e^{-y} y^r \frac{dy}{y} \ll e^{-x}x^{r-1}$ as $x \rightarrow \infty$. The implied constant is only dependent on $A,\re s,k,\vep$ and $\re z$. Similarly, when $\re s>\hf+|\re z|+\vep$ we have
\begin{align}
\int_0^{Y^{-1}2\pi h} y^{s-1} e^{y(1-\delta)}W_{\kt,z}(2y) \frac{dy}{y} \notag 
& \ll 
 \int_0^{Y^{-1}2\pi h} y^{\re s -1}e^{y(1-\delta)}\frac{y^{\hf-|\re z|}}{(1+|\im z|)^A} \frac{dy}{y} \notag \\ 
 & \ll \frac{1}{(1+|\im z|)^A} \int_0^{Y^{-1}2\pi h} y^{\re s-\hf-|\re z|}e^{y(1-\delta)} \frac{dy}{y} \notag \\
 & \ll \frac{(Y^{-1} h)^{\re s-\hf-|\re z|}}{(1+|\im z|)^A}. \label{bottom1}
\end{align} 
Here the implied constant is only dependent on  $A,k,\vep$ and $\re z$. Putting (\ref{top1}) and (\ref{bottom1}) together, we prove the proposition.
 \end{proof}
 
Now we direct our attention to $M_k(s,z/i,\delta)$ and prove it has a meromorphic continuation for all $(s,z) \in \mt{C}^2$. The following proposition locates the poles of $M_k(s,z/i,\delta)$ in both $s$ and $z$ and gives their residues and gives three different growth estimates in regions where $M_k(s,z/i,\delta)$ is analytic, particularly when $\re s \leq 1-\kt$.
\begin{proposition} \label{props} Fix small $\vep>0$ and $\delta>0,$  and let $k \in \mt{R}$. Furthermore let $s=\sigma+ir$ where $\sigma,r \in \mt{R}$ and $\im z =t$. The function $M_k(s,z/i,\delta)$ has a meromorphic continuation to all $(s,z)\in\mt{C}^2$ with simple polar lines at the points $s-\hf\pm z\in \mt{Z}_{\leq 0}$ when $\hf - \kt\pm z \notin \mt{Z}_{\leq 0}$. For fixed $z\notin \hf \mt{Z}$, the residues at these points are given by
\bal
 \label{resfin} \ \ \ \stackrel[s=\hf-\ell\pm z]{}{\mbox{\emph{Res}} }M_k(s,z/i,\delta) =&\frac{(-1)^\ell2^{\hf+\ell\mp z}\Gamma(\hf \mp z-\kt+\ell)\Gamma(\pm 2z-\ell)}{\ell! \Gamma(\hf-\kt + z)\Gamma(\hf-\kt-z)}\\
& \notag +\mc{O}_{\ell,\re z, b}\lt((1+|t|)^{\ell+\kt-\hf-\re z -b}e^{-\frac{\pi}{2}|t|}\delta  \rt) 
\end{align}
where $\ell \in \mt{Z}_{\geq 0}$ and $b< \min(-1,\hf-\sigma-\re z,-2\re z)$. If $\ell\pm 2z \in \mt{Z}_{\geq 0}$ then $ M_k(s,z/i,\delta)$ has a double pole at $s=\hf-\ell \mp z$. Otherwise the poles are simple and are as described above. The Laurent series around these double poles are of the form
\begin{align} \label{laurent}
M(s,z/i,\delta) = \frac{c_2^\pm(\ell,z,k)+\Or_{\ell,z,k}(\delta)}{(s-\thf+\ell\mp z)^{2}}+\frac{c_1^\pm(\ell,z,k)+\Or_{\ell,z,k}(\delta)}{(s-\thf+\ell\mp z)}+\Or_{\ell,z,k}(1)+\Or_{\ell,z,k}(\delta^{1-\vep}).
\end{align}

If we restrict $z$ to the region $0<|\re z| <\vep$, we also have poles in $z$ of the form
\begin{align}
\resi{z=\pm (s+m-\hf)} M_k(s,z/i,\delta)=&\mp \frac{2^{1-s}(-1)^m  \Gamma(2s+m-1)\Gamma(1-s-\tkt)}{m! \Gamma(1-s-m-\tkt)\Gamma(s+m-\kt)}\notag\\
& \ \ \ +\mc{O}_{\sigma,m,k,b}\lt(\frac{\Gamma(2s+m-1)}{\Gamma(s+m-\kt)}(1+|r|)^{1-2\sigma-b} \delta\rt). \label{zeepoles}
\end{align}
when $s$ is near the line $\sigma = \hf-m$. These residues have a meromorphic continuation to 
$\sigma < \hf -m+\vep$ that agrees with the representation above.
  
  For $s$ and $z$ at least a distance of $\vep>0$ from the poles, there exists $A  \in \mt{R}-\mt{Z}$, independent of $\delta$, $r$, and $t$, such that $A> 1+|\sigma| + |\re z|+|\tkt|$ and
\beq
M_k(s,z/i,\delta) \ll_{A,\vep} \delta^{-A}(1+|t|)^{2\sigma - 2 - 2A+k}(1+|r|)^{9A}e^{-\frac{\pi}{2}|r|}.
\label{prop1}
\eeq
 For $\sigma < 1-\kt-\vep_0$ and $s$ at least a distance of $\vep$ away from the poles of $M_k(s,z/i,\delta)$ and $\delta(1+|t|)^2 \leq 1$ we have 
 \bal 
 \label{ll1} M_k(s,z/i,\delta)&=\frac{2^{1-s}\Gamma(s-\frac{1}{2}-z)\Gamma(s-\hf+z)\Gamma(1-s-\kt)}{\Gamma(\frac{1}{2}-\frac{k}{2}+z)\Gamma(\hf-\kt-z)} \\
& \notag + \mc{O}_{A,b,\vep_0} \lt( (1+|t|)^{2\sigma -2+k+2\ep}(1+|r|)^{9A-2b}e^{-\frac{\pi}{2}|r|}\delta^{\vep_0} \rt) 
 \end{align}
while for $\delta(1+|t|)^2>1$ we have
 \beq
M_k(s,z/i,\delta) \ll_{A,\vep} (1+|t|)^{2\sigma - 2 +k}(1+|r|)^{9A}e^{-\frac{\pi}{2}|r|}.
\label{ll2}
\eeq
When $\re z =0$ and $|t|,|r| \gg 1$, $|s\pm z -\hf -m| =\ep>0$,  for $\ep$ small, we have
\beq \label{nearp1}
M_k(s,z/i,\delta) \ll_{m,A,b} \ep^{-1}\delta^{-A}(1+|r|)^{11A-4b}e^{-\frac{\pi}{2}|r|}.
\eeq

\end{proposition} 
 
\begin{proof} From 7.621(3) in \cite{GR} we have that 
\begin{align}
\int_0^\infty e^{-st}t^\alpha W_{\lambda,\mu}(qt) \ dt &
= \frac{\Gamma(\alpha+\mu+\frac{3}{2})\Gamma(\alpha-\mu+\frac{3}{2})q^{\mu+\hf}}{\Gamma(\alpha-\lambda+2)}\lt(s+\hf q\rt)^{-\alpha-\mu-\frac{3}{2}} \\
& 
\times F\lt(\alpha+\mu+\frac{3}{2}, \ \mu-\lambda+\hf; \ \alpha-\lambda+2;  \ \frac{2s-q}{2s+q} \rt) \notag
\end{align}
when $\re\lt(\alpha \pm \mu +\frac{3}{2}\rt)>0, \ \sigma> -\frac{q}{2},$ and $q>0$, and $F$ is a hypergeometric function. Thus when $\re(s-\hf \pm z)>0$ and $\delta>0$,
\begin{align}
\label{gr1}
&M_k(s,z/i,\delta)
=\int_{0}^{\infty}  y^{s-1}e^{y(1-\delta)}W_{\frac{k}{2},z}(2y) \frac{dy}{y} \\
&=\frac{2^{z+\hf}\Gamma(s-\hf+z)\Gamma(s-\hf-z)}{\Gamma(s-\kt)\delta^{s-\hf+z}} \times F\lt(s-\hf+z, \ z-\kt+\hf; \ s-\kt; \ 1-\frac{2}{\delta} \rt). \notag 
\end{align} 
From 9.113 in \cite{GR} we know that when $|\arg(-\eta)|<\pi$ that
\beq \label{bunk}
F(\alpha, \ \beta; \ \gamma; \ \eta)=\frac{\Gamma(\gamma)}{\Gamma(\alpha)\Gamma(\beta)}\times \frac{1}{2\pi i}  \int_C \frac{\Gamma(\alpha+t)\Gamma(\beta+t)\Gamma(-t)}{\Gamma(\gamma+t)}(-\eta)^t \ dt
\eeq
where $C$ is a path of integration chosen such that the poles of $\Gamma(\alpha+t),\Gamma(\beta+t) $ lie to the left of $C$ and the poles of $\Gamma(-t)$ lie to the right. Thus 
\begin{align}
\label{gr2}
& F\lt(s-\hf+z, \ z-\kt+\hf; \ s-\kt ; \  1-\frac{2}{\delta} \rt)  \\ & =\frac{\Gamma(s-\kt)}{\Gamma(s-\hf+z)\Gamma(z-\kt+\hf)} \\
& \ \ \times \frac{1}{2\pi i} \int_C \frac{\Gamma(s-\hf+z+u)\Gamma(z-\kt+\hf+u)\Gamma(-u)}{\Gamma(s-\kt+u)}\lt(\frac{2}{\delta}-1\rt)^u \ du. \notag
\end{align}
 Putting (\ref{gr1}) and (\ref{gr2}) together we get
\begin{align}
M_k(s,z/i,\delta)&=\frac{\Gamma(s-\frac{1}{2}-z)2^{\frac{1}{2}+z}}{\Gamma(z-\frac{k}{2}+\frac{1}{2})\delta^{s-\frac{1}{2}+z}} \label {grow13}\\
& \times \frac{1}{2 \pi i} \int_{C} \frac{\Gamma(s-\frac{1}{2}+z+u)\Gamma(z-\frac{k}{2}+\frac{1}{2}+u)\Gamma(-u)}{\Gamma(s-\frac{k}{2}+u)}\lt(\frac{2}{\delta}-1\rt)^u \ du \notag
\end{align}
where $C$ is a curve where the poles of $\Gamma(s-\frac{1}{2}+z+u),\Gamma(z-\frac{k}{2}+\frac{1}{2}+u)$ are left of $C$ and the poles of $\Gamma(-u)$ are to the right. This immediately gives us a meromorphic continuation of $M_k(s,z/i,\delta)$ to all $(s,z)\in\mt{C}^2$ except potentially along the lines $s-\hf + z \in \mt{Z}_{\leq 0}$ and $z-\kt+\hf \in  \mt{Z}_{\leq 0}$, as the curve $C$ would be undefined. We also note the potential polar lines when $s-\hf-z \in \mt{Z}_{\leq 0}$. Suppose $s$ and $z$ are not along any of these lines. Define $R$  such that $\sigma+\re z = \frac{1}{2}-R$ and let $A>1+|\sigma|+|\re z|+|\kt|$ such that $A$ is some $\vep>0$ distance away from integers. Knowing the residues of the gamma function, we can straighten $C$ to the line $\re u=A$ to get
\beq
M_k(s,z/i,\delta)=\frac{\Gamma(s-\frac{1}{2}-z)2^{\frac{1}{2}+z}}{\Gamma(z-\frac{k}{2}+\frac{1}{2})\delta^{s-\frac{1}{2}+z}}N_1(s,z,\delta) \label{blip2}
\eeq
where
\begin{align}
N_1(s,z,\delta)&=\sum_{0\leq \ell <A} R(s,z,\ell)\label{blip}\\
&\notag+\frac{1}{2\pi i}\int_{(A)} \frac{\Gamma(s-\frac{1}{2}+u+z)\Gamma(\frac{1}{2}-\frac{k}{2}+u+z)\Gamma(-u)}{\Gamma(s-\frac{k}{2}+u)}\lt(\frac{2}{\delta}-1\rt)^{u} \ du
\end{align}
and
\beq
R(s,z,\ell)=\frac{\Gamma(s-\frac{1}{2}+z+\ell)\Gamma(z-\frac{k}{2}+\frac{1}{2}+\ell)(-1)^{\ell}}{\ell! \Gamma(s-\frac{k}{2}+\ell)}\lt(\frac{2}{\delta}-1\rt)^\ell.
\label{pole}
\eeq
We see that the integral part of \eqref{blip} contributes no poles in $s$ or $z$ as $\sigma - \frac{1}{2}+\re z+\re u > -A+A=0$ and $\re z -\tfrac{k-1}{2}+\re u > -A+A=0$. Noting the cancellation of simple poles of $\Gamma(z-\kt+\hf+\ell)$ by the $\Gamma(z-\kt+\hf)$ in the denominator of \eqref{blip2},  the only poles of $M_k(s,z/i,\delta)$ are at $s=\frac{1}{2}\pm z-m$ for $m \in \mt{Z}_{\geq 0}$.

Still assuming that $s-\hf\pm z \notin \mt{Z}_{\leq 0}$ and $A+R \notin \mt{N}_{0}$, shift the line of integration to  $\re u=A+R$. Performing the change of variables $u \rightarrow u-(s-\frac{1}{2}+z)$ we get 
\beq
M_k(s,z/i,\delta)=\frac{\Gamma(s-\frac{1}{2}-z)2^{\frac{1}{2}+z}}{\Gamma(z-\frac{k}{2}+\frac{1}{2})\delta^{s-\frac{1}{2}+z}}N_2(s,z,\delta)
\eeq
where 
\begin{align}
N_2(s,z,\delta)=&\sum_{0\leq \ell <A+R} R(s,z,\ell) \label{awed}\\
&\notag +\frac{1}{2\pi i}\int_{(A)} \frac{\Gamma(u)\Gamma(1-s-\frac{k}{2}+u)\Gamma(s-\frac{1}{2}+z-u)}{\Gamma(\frac{1}{2}-\frac{k}{2}-z+u)}\lt(\frac{2}{\delta}-1\rt)^{u-s+\frac{1}{2}-z} \ du.
\end{align}
Let $s=\sigma+ir$ with $\sigma,r\in\mt{R}$, $\im z =t$ and $u=A+iv$ with $v \in \mt{R}$. Stirling's Approximation gives us that for $s$ at least $\vep>0$ from the poles of $M_k(s,z,\delta)$ we have
\beq \frac{\Gamma(s-\frac{1}{2}-z)2^{\frac{1}{2}+z}}{\Gamma(z-\frac{k}{2}+\frac{1}{2})\delta^{s-\frac{1}{2}+z}} \ll_{A,\vep} \delta^R(1+|r-t|)^{\sigma-\re z-1}(1+|t|)^{\kt-\re z}e^{-\frac{\pi}{2}(|r-t|-|t|)}
\label{grow5}
\eeq
and
\beq
R(s,z,\ell) \ll_{A,\vep} \delta^{-\ell} \frac{(1+|r+t|)^{\sigma+\ell+\re z -1}(1+|t|)^{-\frac{k}{2}+\re z+\ell}}{(1+|r|)^{-\kt+\sigma+\ell-\hf}}e^{-\frac{\pi}{2}(|r+t|+|t|-|r|)}.
\eeq

Thus
\begin{align}
\label{grow3}
&\frac{\Gamma(s-\frac{1}{2}-z)2^{\frac{1}{2}+z}}{\Gamma(z-\frac{k}{2}+\frac{1}{2})\delta^{s-\frac{1}{2}+z}}  R(s,z,\ell) \\
&\ll_{A,\vep} \notag \delta^{R-\ell}  \frac{(1+|r+t|)^{\sigma+\ell+\re z -1}(1+|r-t|)^{\sigma-\re z -1}(1+|t|)^{\ell}}{(1+|r|)^{\sigma+\ell-\hf-\kt}}e^{-\frac{\pi}{2}(2\max(|r|,|t|)-|r|)}.
\end{align}
Given that  $\ell \leq \lfloor A+R \rfloor$, we see that 
\begin{align}
\label{grow4}
\frac{\Gamma(s-\frac{1}{2}-z)2^{\frac{1}{2}+z}}{\Gamma(z-\frac{k}{2}+\frac{1}{2})\delta^{s-\frac{1}{2}+z}} &\sum_{0\leq 0 <A+R} R(s,z,\ell) \\
& \notag \ll_{A,\vep} \delta^{-A}(1+|t|)^{4A}(1+|r|)^{3A} e^{-\frac{\pi}{2}(2\max(|r|,|t|)-|r|)}.
\end{align}
Now, using Stirling's Approximation for $\re u =A$ we are able to note that
\begin{align}
&\frac{\Gamma(u)\Gamma(1-s-\frac{k}{2}+u)\Gamma(s-\frac{1}{2}+z-u)}{\Gamma(\frac{1}{2}-\frac{k}{2}-z+u)}\lt(\frac{2}{\delta}-1\rt)^{u-s+\frac{1}{2}-z} \label{grow6}\\
& \ll \delta^{-A-R} \frac{(1+|v|)^{A-\frac{1}{2}}(1+|r-v|)^{A-\sigma+\hf-\kt}(1+|r+t-v|)^{\sigma-A-1+\re z}}{(1+|t-v|)^{A-\kt-\re z}e^{\frac{\pi}{2}(|v|+|r-v|+|r+t-v|-|t-v|)}}. \notag
\end{align} 
Combining (\ref{grow5}) and (\ref{grow6}) we get an exponent of $|v|+|r-v|+|r+t-v|-|t-v|+|r-t|-|t|$. 
\begin{lemma}\label{dumblemma}
For $v,r,t \in \mt{R}$, 
\beq
 |v|+|r-v|+|r+t-v|-|t-v|+|r-t|-|t|\geq \max(|r|,|v|-|r|).
 \eeq
\end{lemma}
\begin{proof}
We can explicitly compute that
\begin{align}
M_1:&=|v|+|r-v|+|r+t-v|+|r-t|\notag \\
&=\max(|3r-v|,|r+2t-v|,3|r-v|,|r+2t-3v|,|r+v|,|v+2t-r|,|2t-r-v|)
\end{align}
and
\beq
M_2:=|t-v|+|t|=\max(|2t-v|,|v|).
\eeq Furthermore, we see that
\beq
M_2+|v|=\max(|2v|,|2t-2v|,|2t|)
\eeq
and
\beq
M_2+|r|=\max(|2t-v+r|,|2t-v-r|,|v-r|,|v+r|)\leq M_1,
\eeq
so $M_1-M_2 \geq |r|$. We also see that $M_2+|v|=|2t-2v|$ if and only if $v$ and $t$ have opposite signs, in which case 
\beq
M_1+|r|\geq |r+2t-3v|+|r| \geq |2t-3v|\geq |2t-2v| =M_2+|v|.
\eeq
If $M_2+|v|=|2t|$, then $v$ and $t$ have the same sign so
\beq
M_1+|r|\geq |v+2t-r|+|r| \geq |v+2t| \geq |2t| =M_2+|v|.
\eeq
In the case when $M_2+|v|=|2v|$ then $v$ and $t$ again have the same sign. If also $|v| \leq |2t|$ then 
\beq
M_1+|r| \geq |r|+|r+v+2t|\geq |v+2t| \geq |v|+|2t| \geq |2v|=M_2+|v|,
\eeq
and if instead $|v| \geq |2t|$ then
\beq
M_1+|r| \geq |r+2t-3v| +|r| \geq |2t-3v| \geq 3|v|-2|t| \geq 3|v|-|v| \geq |2v| =M_2+|v|.
\eeq
This all means that $M_1-M_2 \geq |v|-|r|$ in every case, proving the lemma.  
\end{proof}
When $|t| \gg |r|$ we see that
\begin{align}
&\frac{\Gamma(s-\frac{1}{2}-z)2^{\frac{1}{2}+z}}{\Gamma(z-\frac{k}{2}+\frac{1}{2})\delta^{s-\frac{1}{2}+z}}  \int_{(A)} \frac{\Gamma(u)\Gamma(1-s-\frac{k}{2}+u)\Gamma(s-\frac{1}{2}+z-u)}{\Gamma(\frac{1}{2}-\frac{k}{2}-z+u)}\lt(\frac{2}{\delta}-1\rt)^{u-s+\frac{1}{2}-z} \ du \notag \\
& \ \ \  \ll_A \delta^{-A}(1+|t|)^{2\sigma-2-2A+k} e^{-\frac{\pi}{2}|r|}\int_{A-2ir}^{A+2ir} (1+|v|)^{A-\hf}(1+|v-r|)^{A-\kt+\hf-\sigma} \ du \notag \\
& \ \ \ \ \ + \delta^{-A} (1+|t|)^{\sigma+\kt-2\re z-1}e^{-\frac{\pi}{2}|r|} \int\limits_{\re u =A \atop |v|>|2r|} \frac{(1+|v|)^{2A-\sigma-\kt}}{(1+|t-v|)^{2A-2\re z -\sigma-\kt +1 }} e^{-\frac{\pi}{2}(|v|-2|r|)} \  du \notag \\
& \ \ \ \ll_A \delta^{-A} (1+|t|)^{2\sigma+k-2-2A} (1+|r|)^{3A} e^{-\frac{\pi}{2}|r|} \label{grow14}.
\end{align}
By a similar argument, when $|t| \sim |r|$ or $|r| \gg |t|$  and $s$ and $z$ are at least $\vep$ away from the poles of $M_k(s,z/i,\delta)$, then
\begin{align}
&\frac{\Gamma(s-\frac{1}{2}-z)2^{\frac{1}{2}+z}}{\Gamma(z-\frac{k}{2}+\frac{1}{2})\delta^{s-\frac{1}{2}+z}}  \int_{(A)} \frac{\Gamma(u)\Gamma(1-s-\frac{k}{2}+u)\Gamma(s-\frac{1}{2}+z-u)}{\Gamma(\frac{1}{2}-\frac{k}{2}-z+u)}\lt(\frac{2}{\delta}-1\rt)^{u-s+\frac{1}{2}-z} \ du \notag \\
& \ \ \ \ll_{A,\vep} \delta^{-A}(1+|r|)^{4A}(1+|t|)^{\kt-\re z}e^{-\frac{\pi}{2} |r|} \label{grow24}
\end{align}
and since when $|t| \leq |r|$, 
\beq
(1+|t|)^{2\sigma+k-2-2A}(1+|r|)^{9A} \geq (1+|r|)^{4A}(1+|t|)^{\kt-\re \, z},
\eeq
  we get (\ref{prop1}) in the proposition by combining \eqref{grow4} with \eqref{grow14} and \eqref{grow24}.

Assume for the present that $s+\kt-1 \notin \mt{Z}$. If we go back to the contour integral representation in \eqref{grow13} and straighten the line of integration to $\re u=b$ where $b<\min(-2\re z,R,-1)$ and this line does not cross any of the poles of the integrand. We move past the poles at $u=\frac{1}{2}-s-\ell-z$ for $0\leq \ell \leq \lf R-b\rf $ and $u=\frac{k-1}{2}-\ell-z$ for $0 \leq \ell \leq \lf\frac{k-1}{2}-\re z -b \rf$ getting that  
\begin{subequations}\label{thelot}
\beq
M_k(s,z/i,\delta)=C+D+E,
\eeq
 where 
\begin{align}
C&=\frac{2^{\frac{1}{2}+z}\Gamma(s-\frac{1}{2}-z)}{\Gamma(\frac{1}{2}-\frac{k}{2}+z)} \label{cee} \\
&\times \notag \sum_{\ell=0}^{\lf R-b \rf} \frac{(-1)^\ell\Gamma(1-s-\frac{k}{2}-\ell)\Gamma(s-\frac{1}{2}+z+\ell)}{\ell! \Gamma(\frac{1}{2}-\frac{k}{2}-z-\ell)}(2-\delta)^{\frac{1}{2}-s-z}\lt(\frac{2}{\delta}-1\rt)^{-\ell} 
\end{align}
\bal
D&=\label{cee2}\frac{2^{\frac{1}{2}+z}\Gamma(s-\frac{1}{2}-z)}{\delta^{s-\frac{1}{2}+z}\Gamma(\frac{1}{2}-\frac{k}{2}+z)} \\
&\times \notag \sum_{\ell=0}^{\lf \frac{k-1}{2}-\re z-b\rf} \frac{(-1)^\ell\Gamma(s+\frac{k}{2}-1-\ell)\Gamma(\frac{1}{2}-\frac{k}{2}+\ell+z)}{\ell! \Gamma(s-\frac{1}{2}-\ell-z)}\lt(\frac{2}{\delta}-1\rt)^{\frac{k-1}{2}-\ell-z}
\end{align}
\bal
E&=\label{cee3} \frac{2^{\frac{1}{2}+z}\Gamma(s-\frac{1}{2}-z)}{\delta^{s-\frac{1}{2}+z}\Gamma(\frac{1}{2}-\frac{k}{2}+z)} \\
& \notag \times \frac{1}{2\pi i} \int_{(b)}\frac{\Gamma(s-\frac{1}{2}+z+u)\Gamma(z-\frac{k}{2}+\frac{1}{2}+u)\Gamma(-u)}{\Gamma(s-\frac{k}{2}+u)} \lt(\frac{2}{\delta}-1\rt)^u \ du. 
\end{align}
\end{subequations}
 It's worth noting that $D=0$ if $\tfrac{k-1}{2}-\re z-b <0$. We also note that $C$ and $D$ both appear to have poles for when $s+\kt-1 \in \mt{Z}$ but either recalling the poles of $M_k(s,z/i,\delta)$ given by the formula (\ref{blip}), or by explicitly computing residues of both $C$ and $D$ and observing that they cancel, we see that these poles do not exist, thus we may also allow such values of $s$. In any case, from the above representation of $M_k(s,z/i,\delta)$ we can compute the residues of the poles at $s=\hf\pm z-\ell$, in particular, it is easy to see that $D$ does not contribute residues for any of these poles. So for fixed  $z\notin \hf\mt{Z} $, 
\bal
& \label{res1} \ \ \ \stackrel[s=\hf-m+z]{}{\res} M_k(s,z/i,\delta) = \frac{2^{\hf+z}(-1)^m}{m! \Gamma(\hf-\kt+z)}\\
& \times \lt[\sum_{\ell=0}^{\lfloor m-\re 2z -b \rfloor} \frac{(-1)^\ell\Gamma(\hf+m-z-\kt-\ell)\Gamma(-m+2z+\ell)(2-\delta)^{m-2z}(\frac{2}{\delta}-1)^{-\ell}}{\ell! \Gamma(\hf-\kt-z-\ell)} \rt. \notag\\
& \ \ \  \lt. + \frac{1}{2\pi i} \int_{(b)} \frac{\Gamma(-m+2z+u)\Gamma(\hf+z-\kt+u)\Gamma(-u)\delta^{m-2z}(\frac{2}{\delta}-1)^{u}}{\Gamma(\hf-m+z-\kt+u)} \ du \rt] \notag
\end{align}
 and
\bal
\label{res2}  \stackrel[s=\hf-m-z]{}{\res} M_k(s,z/i,\delta) =& \frac{2^{\hf+z}(-1)^m\Gamma(-m-2z)(2-\delta)^m}{\Gamma(\hf-\kt+z)} \\
& \times \sum_{\ell =0}^{\lf m-b \rf} \frac{\Gamma(\hf+m+z-\kt-\ell)(\frac{2}{\delta}-1)^{-\ell}}{\ell! (m-\ell)!\Gamma(\hf-\kt-z-\ell)}. \notag
\end{align}
Again recall $\im z =t$. It is easy to see that the dominant term with respect to $\delta$ in the sum of (\ref{res1}) comes from $\ell=0$, and letting $u=b+iv$, Stirling's Approximation gives us that
\bal
& \label{res3} \ \ \ \stackrel[s=\hf-m+z]{}{\res} M_k(s,z/i,\delta) =\frac{2^{\hf+m-z}(-1)^m\Gamma(\hf+m-z-\kt)\Gamma(-m+2z)}{m! \Gamma(\hf-\kt+z)\Gamma(\hf-\kt-z)}\\
& \notag +\mc{O}_{m,\re z, b}\lt((1+|t |)^{m+\kt-\hf-\re z-b}e^{-\frac{\pi}{2}|t |}\delta  \rt) \\
&\notag +\mc{O}_{m,\re z, b}\lt(\int_{-\infty}^\infty \delta^{m-b}\frac{(1+|2t +v|)^{b-m-\hf}(1+|t +v|)^{m}(1+|t |)^{\kt-\re z}}{(1+|v|)^{\hf+b}} e^{-\frac{\pi}{2}(|2 t  +v|+|v|-|t |)} \ dv \rt).
\end{align}
We see that since $|2t+v|+|v|-|t|\geq |v|-|t|$ and $|t|$ everywhere, an argument like that given in \eqref{grow14} gives that
\bal
 \label{res4} \ \ \ \stackrel[s=\hf-m+z]{}{\res} M_k(s,z/i,\delta) =&\frac{2^{\hf+m-z}(-1)^m\Gamma(\hf+m-z-\kt)\Gamma(-m+2z)}{m! \Gamma(\hf-\kt+z)\Gamma(\hf-\kt-z)}\\
& \notag +\mc{O}_{m,\re z, b}\lt((1+|t |)^{m+\kt-\hf-\re z-b}e^{-\frac{\pi}{2}|t |}\delta  \rt).
\end{align}
In (\ref{res2}), $\ell=0$ also gives the dominant term with respect to $\delta$. Stirling's Approximation easily gives us that
\bal 
\label{res5} \stackrel[s=\hf-m-z]{}{\res} M_k(s,z/i,\delta) =& \frac{2^{\hf+m+z}(-1)^m\Gamma(-m-2z)\Gamma(\hf+m+z-\kt)}{m!\Gamma(\hf-\kt+z)\Gamma(\hf-\kt-z)} \\
& +\mc{O}_{m,\re z,b}\lt(  (1+|t|)^{\kt-\hf-\re z}e^{-\frac{\pi}{2}|t|} \delta \rt). \notag
\end{align}
Putting together (\ref{res4}) and (\ref{res5}) we get (\ref{resfin}) in the proposition.

Now suppose that $z \in \hf\mt{Z}$ and $m \in \mt{Z}_{\geq 0}$. If $m+2z \in \mt{Z}_{\geq 0}$ then $ M_k(s,z/i,\delta)$ has a double pole at $s=\hf-m-z$. If $m-2z \in \mt{Z}_{\geq 0}$ then $ M_k(s,z/i,\delta)$ has a double pole at $s=\hf-m+z$. Otherwise the poles are simple as described above. When double poles occur it is not difficult to show from \eqref{cee} that the $n$th coefficient
of the Laurent series around these double poles at $s=\hf-m\pm z$ are of the form $c_n^\pm(m,z,k)+\mc{O}_{m,z,k}(\delta^{1-\vep})$ where there $c_n^\pm$ can be computed explicitly. This gives us \eqref{laurent}. 

Now we consider the poles of $M_k(s,z/i,\delta)$ at $z=\pm(\hf-s-m)$, when $0<|\hf-s-m| < \vep$. Using \eqref{thelot} we compute the residues: 
\bal
R_1(s,m,\delta): = \resi{z=\hf-s-m} M_k(s,z/i,\delta) =& \frac{2^{1-s-m}(-1)^m(2-\delta)^m\Gamma(2s+m-1)}{\Gamma(1-s-m-\kt)}\notag \\
& \ \ \ \  \times \sum_{\ell =0}^{\lf m-b \rf} \frac{\Gamma(1-s-\tkt-\ell)(\tfrac{2}{\delta}-1)^{-\ell}}{\ell! (m-\ell)!\Gamma(s+m-\kt-\ell)}
\end{align}
and
\begin{align}
\label{therealthing} R_2(s,m,\delta): =& \resi{z=s+m-\hf} M(s,z/i,\delta)=-\frac{2^{s+m}(-1)^m}{m!\Gamma(s+m-\tkt)} \\&\notag \times  \lt[ \sum_{\ell=0}^{\lf 1-2\sigma-m-b \rf} \frac{(-1)^\ell\Gamma(1-s-\kt-\ell)\Gamma(2s+m-1+\ell)(2-\delta)^{1-m-2s}\lt(\frac{2}{\delta}-1 \rt)^{-\ell}}{\ell! \Gamma(1-s-m-\ell-\tkt)} \rt. \notag  \\
& \lt. \ \ \ \ + \frac{1}{2\pi i} \int_{(b)} \frac{\Gamma(2s+m-1+u)\Gamma(s+m-\kt+u)\Gamma(-u)\delta^{1-2s-m}\lt( \frac{2}{\delta}-1 \rt)^u}{\Gamma(s-\tkt+u)} \ du  \rt] \notag.
\end{align}
We see that $R_1(s,m,\delta)$ easily has a meromorphic continuation to all $s \in \mt{C}$ to the left $\sigma = \hf-m+\vep$. Stirling's approximation gives us that
\begin{align}
R_1(s,m,\delta) =& \frac{2^{1-s}(-1)^m  \Gamma(2s+m-1)\Gamma(1-s-\tkt)}{m! \Gamma(1-s-m-\tkt)\Gamma(s+m-\kt)} \notag \\
& \ \ \ + \mc{O}_{\sigma,m,k,b}\lt( \frac{\Gamma(2s+m-1)}{\Gamma(s+m-\kt)}(1+|r|)^{m} \delta  \rt).
\end{align}
Comparatively, $R_2(s,m,\delta)$ has a less obvious continuation: while the finite sum has an easy meromorphic continuation to all $s$, the integral component requires a bit more attention. Consider the function given by the integral
\beq
J(s,m,\delta):=\frac{1}{2\pi i} \int_{(b)} \frac{\Gamma(2s+m-1+u)\Gamma(s+m-\kt+u)\Gamma(-u)\delta^{1-2s-m}\lt( \frac{2}{\delta}-1 \rt)^u}{\Gamma(s-\tkt+u)} \ du 
\eeq
when $\sigma $ is near $\hf -m$, specifically to the right of $2 \sigma + 2m +\{ b \}-1=0$. We see that the integral itself, which is not to be confused with the meromorphic continuation of $J(s,m,\delta)$, is absolutely convergent for all fixed $s$ to the left of $\sigma = \hf-m+\vep$, except for the lines $2 \sigma + 2m +\{b\}-1=-r$ for $r \geq 0$. We can get a meromorphic continuation of $J(s,m,\delta)$ by moving $s$ left up to the right side of each subsequent line, shifting the line of integration to the right past the now nearby pole, then moving $s$ further left and shifting the line of integration back. Repeating as we move $s$ past each line we get 
\begin{align}
J(s,m,\delta)=&\frac{1}{2\pi i} \int_{(b)} \frac{\Gamma(2s+m-1+u)\Gamma(s+m-\kt+u)\Gamma(-u)\delta^{1-2s-m}\lt( \frac{2}{\delta}-1 \rt)^u}{\Gamma(s-\tkt+u)} \ du \notag \\
& + \sum_{\ell=m-\lf b \rf}^{m-\lf b \rf+r} \frac{(-1)^\ell \Gamma(1-s-\ell-\tkt)\Gamma(2s+m+\ell-1)(2-\delta)^{1-2s-m} (\frac{2}{\delta}-1)^{-\ell}}{\ell!\Gamma(1-s-\ell-m-\tkt)}
\end{align}
when $-\frac{1}{2}-m-\frac{\{b\}}{2}-\tfrac{r}{2} < \sigma < \tfrac{1}{2} -m-\frac{\{b\}}{2} -\frac{r}{2}$ for $r \geq 0.$ Thus $J(s,m,\delta)$ has a meromorphic continuation to the left of $\sigma = \hf -m +\vep$.  From this we have that $R_2(s,m,\delta)$ indeed has an meromorphic continuation to the same region and is accurately given by the formula in \eqref{therealthing}. Another application of Stirling's gives us that 
\begin{align}
R_2(s,m,\delta)=&-\frac{(-1)^m2^{1-s} \Gamma(1-s-\tkt)\Gamma(2s+m-1)}{m! \Gamma(s+m-\tkt)\Gamma(1-s-m-\tkt)} \notag \\
&\ \ \ +\mc{O}_{\sigma,m,k,b}\lt(\frac{\Gamma(2s+m-1)}{\Gamma(s+m-\kt)}(1+|s|)^{1-2\sigma-b} \delta \rt).
\end{align}
Thus we have the result corresponding to \eqref{zeepoles} in the proposition.

Now again recall that $A>1+|\sigma|+|\re z| +|\kt|$. We want to use these representations of $C,D,E$ in \eqref{thelot} to produce subtler growth estimates as $\delta \rightarrow 0$ and $\im z \to \infty$. Suppose $\sigma <1-\kt$, then for $C$ in \eqref{cee}, we note that the dominant growth term with respect to $\delta$ arises when $\ell=0$. With this in mind we have that when $s$ and $z$ are at least $\vep>0$ away from the poles at $s\pm z -\hf \in \mt{Z}_{\leq 0}$, 
\bal
C&=\frac{2^{1-s}\Gamma(s-\frac{1}{2}-z)\Gamma(s-\hf+z)\Gamma(1-s-\kt)}{\Gamma(\frac{1}{2}-\frac{k}{2}+z)\Gamma(\hf-\kt-z)} \\
&+\sum_{\ell=1}^{\lf R-b \rf} \mc{O}_{A,b,
\vep}\lt( \frac{(1+|r-t|)^{\sigma-1-\re z} (1+|t|)^{\ell+k}(1+|r|)^{\hf-\sigma -\kt-\ell}(1+|r+t|)^{\sigma+\ell-1+\re z}}{e^{\frac{\pi}{2}(|r-t|+|r|+|r+t|-2|t|)}} \delta^ \ell \rt).
 \notag
\end{align}
Since $|r-t|+|r|+|r+t|-2|t| =2\max(|r|-|t|,0)+|r|$ we get
 \bal
\label{growc} C&=\frac{2^{1-s}\Gamma(s-\frac{1}{2}-z)\Gamma(s-\hf+z)\Gamma(1-s-\kt)}{\Gamma(\frac{1}{2}-\frac{k}{2}+z)\Gamma(\hf-\kt-z)} \\
&+\mc{O}_{A,b,\vep}\lt( (1+|r|)^{9A-2b} (1+|t|)^{2\sigma-2+k}e^{\frac{\pi}{2}|r|} \max_{\ell=1,\lf R-b \rf} \lt(\delta^ \ell(1+|t|)^{2\ell}\rt) \rt). 
 \notag
\end{align}
If $\sigma < 1 -\kt$ then the exponent of $\delta$ is always greater than 0 for every term in $D$ in \eqref{cee2}, unless $\tfrac{k-1}{2}-\re z -b < 0$, where $D=0$. Using Stirling's approximation we get that
\bal
D \ll_{A,b,\vep} \max_{\ell=0,\lf \frac{k-1}{2}-\re z-b\rf} \lt( \delta^{1-\sigma-\kt+\ell} (1+|r-t|)^\ell (1+|r|)^{\sigma+\kt-\frac{3}{2}-\ell}(1+|t|)^\ell e^{-\frac{\pi}{2}|r|}   \rt).
\end{align}
which gives us
\beq 
\label{growd}  D \ll_{A,b,\vep} \delta^{1-\sigma-\kt}  (1+|r|)^{9A-2b}e^{-\frac{\pi}{2}|r|}\max_{\ell=0,\lf \frac{k-1}{2}-\re z-b\rf}\lt(\delta^\ell (1+|t|)^{2\ell}    \rt).
\eeq

Now we deal with $E$ in \eqref{cee2}. By the change of variables $u \to u-s+\hf-z$, Stirling's approximation and Lemma \ref{dumblemma} to get 
\bal
E \ll_{A,b,\vep} & \  \delta^{R-b} (1+|t|)^{-1-2b+k-2\re z}(1+|r|)^{9A-2b}e^{-\frac{\pi}{2}|r|} \label{growe},
\end{align}
by following an argument analogous to that given following \eqref{grow14}.

Consider now the two cases: $\delta(1+|t|)^2 \leq 1$ and $\delta(1+|t|)^2 >1$.  In the prior case, for small enough $\vep>0$
\begin{align*}
&\max_{\ell =1,\lf R-b \rf} (\delta^\ell(1+|t|)^{2\ell})=\delta(1+|t|)^2 \leq \delta^\vep(1+|t|)^{2\vep} \\
&\max_{\ell =\lf \frac{k-1}{2}-\re z-b \rf,0} (\delta^\ell(1+|t|)^{2\ell})=1
\end{align*}
so the error term in  (\ref{growc}) is
\begin{align}
\mc{O}_{A,b,\vep}\lt( (1+|r|)^{9A-2b} (1+|t|)^{2\sigma-2+k+2\vep}e^{-\frac{\pi}{2}|r|}\delta^\vep \rt). \label{err1}
\end{align}
From \eqref{growd}
\beq 
 D \ll_{A,b,\vep} \delta^{1-\sigma-\kt}  (1+|r|)^{9A-2b}e^{-\frac{\pi}{2}|r|}. \label{err2},
\eeq
supposing $D \neq 0$. Finally from (\ref{growe}) we have
\begin{align}
&\notag  \ll (1+|t|)^{2\sigma -2 +k+2\ep}(1+|r|)^{9A-2b}e^{-\frac{\pi}{2}|r|}\delta^R.
\end{align}
Putting all of this together we get (\ref{ll1}) in the proposition.

In the second case, $\delta(1+|t|)^2 >1$, we return to the first estimate in (\ref{prop1}). Substituting $\delta^{-A} < (1+|t|)^{2A}$ we obtain  (\ref{ll2}) in the proposition.

 Let $\re  z =0$, and $|t|,|r| \gg 1$ and suppose that $|s+z-\hf-m|=\vep>0$, which forces $t \sim -r$. Then using \eqref{thelot}, Stirling's approximation, \eqref{growd}, \eqref{growe} and the bound
$$
\Gamma(s-\thf+z+\ell) \ll_m \ep^{-1},
$$
when $\ell \leq m$, we have that 
\beq \label{nearp}
M_k(s,z/i,\delta) \ll_{m,A,b}  \ep^{-1} \delta^{-A} (1+|r|)^{11A-4b}e^{-\frac{\pi}{2}|r|}.
\eeq
Similarly suppose that $|s-z-\hf-m|=\ep>0$, which forces $t \sim r$. Again using \eqref{thelot}, Stirling's approximation, and that
$$
\Gamma(s-\thf-z) \ll_m \ep^{-1},
$$
in this region, we have that
\beq
M_k(s,z/i,\delta) \ll_{m,A,b} \ep^{-1} \delta^{-A} (1+|r|)^{11A-4b}e^{-\frac{\pi}{2}|r|}
\eeq
which together with \eqref{nearp} gives \eqref{nearp1}, completing the proof of the proposition.
\end{proof}

\section{The Limit As $Y\rightarrow \infty$}
We can now analyze the uniformity of the convergence of the spectral decomposition of $\mc{I}_{f,\ell, Y,\delta}(s;h)$ given in (\ref{spec2}). We begin by letting $\re s > 1$. Recall that by (\ref{spec2}) we may write
\beq 
\mc{I}_{f,\ell, Y,\delta}(s;h)=\mc{I}_{f,\ell, Y,\delta}^{\mbox{\tiny cusp}}(s;h)+\mc{I}_{f,\ell, Y,\delta}^{\mbox{\tiny con}}(s;h)-\LA  P_{h,Y}^{(2)}(*;s;\delta),1 \RA
\eeq 
where
\bal
\mc{I}_{f,\ell, Y,\delta}^{\mbox{\tiny cusp}}(s;h) & = \sum_{j=1}^\infty \frac{\ol{2\rho_j(-1)\lambda_j(h)}}{\mcV (2\pi h)^{s-\hf}} \label{sum1}\\
& \ \times \lt(\int_{2\pi h Y^{-1}}^{2 \pi h Y} y^{s-\hf} e^{ y (1-\delta)}K_{it_j}(y) \frac{dy}{y}\rt) \ol{\LA V_{f,\ell}, u_j \RA} \notag \\
& =  \sum_j\sqrt{2\pi}\frac{\ol{\rho_j(-h)}}{\mcV (2\pi h)^{s-\hf}}M_{Y,h,0}(s,t_j,\delta) \ol{\LA V_{f,\ell}, u_j \RA} 
\notag
\end{align}
and
\bal
\mc{I}_{f,\ell, Y,\delta}^{\mbox{\tiny con}}(s;h) & = \frac{1}{4\pi}\int_{-\infty}^\infty  \frac{2h^{-it}\sigma_{2it}(h)}{\mcV \zeta^*(1-2it)(2\pi h)^{s-\hf}} \label{integ1} \\ 
& \ \times \notag \lt( \int_{2\pi h Y^{-1}}^{2 \pi h Y} y^{s-\hf} e^{y (1-\delta)}K_{it}( y) \frac{dy}{y}\rt)\ol{\LA V_{f,\ell}, E(*,1/2+it) \RA}  \ dt  \\
& =\frac{1}{4\pi} \int_{-\infty}^\infty  \sqrt{2\pi}\frac{h^{-it}\sigma_{2it}(h)}{\mcV \zeta^*(1-2it)(2\pi h)^{s-\hf}} M_{Y,h,0}(s,t,\delta) \ol{\LA V_{f,\ell}, E(*,1/2+it) \RA} \ dt.\notag
\end{align}
By Theorem 3 in \cite{Watson} we have that
\beq
\langle V_{f,\ell}, u_j\rangle^2 = \frac{\langle f,f\rangle}{8\Lambda(1,\mbox{Ad} f)} \frac{\langle \mu_\ell,\mu_\ell \rangle \langle u_j,u_j \rangle \Lambda(\hf,f \otimes \ol{\mu}_\ell \otimes \ol{u}_j)}{\Lambda(1,\mbox{Ad} \ol{\mu}_\ell)\Lambda(1,\mbox{Ad} \ol{u}_j)}
\eeq
where all the completed $L$-functions are normalized such that $\rho_j(1)=\rho_\ell(1)=1$. From the gamma functions associated with these $L$-functions, Stirling's formula gives us that, when $t_j,t_\ell \in \mt{R}$, 
\begin{align}
&\langle V_{f,\ell}, u_j\rangle \label{youngones} \\
 &\ll_f (1+|t_\ell+t_j|)^{\frac{k-1}{2}}(1+|t_\ell-t_j|)^{\frac{k-1}{2}}\log(1+|t_j|)\log(1+|t_\ell|)e^{-\pif ||t_j|-|t_\ell||} \sqrt{\left|L(\thf,f \otimes \ol{\mu}_\ell \otimes \ol{u}_j)\rt|} \notag
\end{align}
or more simply
\beq
\langle V_{f,\ell}, u_j\rangle \ll_f (1+|t_\ell+t_j|)^{\frac{k-1+\beta_0}{2}}(1+|t_\ell-t_j|)^{\frac{k-1+\beta_0}{2}}\log(1+|t_j|)\log(1+|t_\ell|)e^{-\pif ||t_j|-|t_\ell||}, \label{youngone}
\eeq
where the Generalized Lindel\"{o}f Hypothesis would give that $\beta_0=0$. The convexity bound gives $\beta_0 \leq 1$. 
By a more straightforward computation we get
\bal
\ol{\langle V_{f,\ell}, E(*,s)\rangle}=& \ol{ \langle y^{\kt}\ol{f}\mu_{\ell,k}(z),E(*,s)\rangle}  = \frac{\ol{\rho_{\ell,k}(1)}}{(4\pi)^{\kt-1}\mcV }\frac{\Lambda(s,f\otimes \ol{\mu_{\ell}})}{\zeta^*(2s)},
\end{align}
where $\zeta^*(2s)=\pi^{-s}\Gamma(s)\zeta(2s)$ is the completed zeta function and 
\beq
\Lambda(s,f\otimes \ol{\mu}_{\ell,k}):=(2\pi)^{-2s}\Gamma(s+\tfrac{k-1}{2}+it_\ell)\Gamma(s+\tfrac{k-1}{2}-it_\ell)L(s,f\otimes \ol{\mu_{\ell,k}})=\Lambda(1-s,f\otimes \ol{\mu_{\ell,k}}).
\eeq
From this we get, where again we denote $s=\sigma+ir$, that
\bal
\ol{\langle V_{f,\ell}, E(*,s)\rangle}  \ll_f \frac{(1+|t_\ell+r|)^{\sigma+\kt-1}(1+|t_\ell-r|)^{\sigma+\kt-1} \log(1+|t_\ell|)L(s,f\otimes \ol{\mu_{\ell}})}{(1+|r|)^{\sigma-\hf}\zeta(2s)e^{\pif ||r|-|t_\ell||}} . \label{youngtwoone}
\end{align}
and so when $t \in \mt{R}$
\bal
\ol{\langle V_{f,\ell}, E(*,\thf+it)\rangle}  & \ll_f  (1+|t_\ell+t|)^{\frac{k-1}{2}}(1+|t_\ell-t|)^{\frac{k-1}{2}} \log(1+|t_\ell|) e^{-\pif ||t|-|t_\ell||}L(\thf+it,f\otimes \ol{\mu}_{\ell,k}) \notag \\
& \ll_f (1+|t_\ell+t|)^{\frac{k-1+\beta_0'}{2}}(1+|t_\ell-t|)^{\frac{k-1+\beta_0'}{2}} \log(1+|t_\ell|) e^{-\pif ||t|-|t_\ell||} \label{youngtwo}
\end{align}
where again the Generalized Lindel\"{o}f Hypothesis implies that $\beta_0'=0$, and the convexity bound gives $\beta_0' \leq 1$. Henceforth we let let $\beta_1=\max(\beta_0,\beta_0')$.

Summing over eigenvalues of Maass forms, Weyl's law gives us that
\bal
\sum_{|t_j| \sim T} 1 \ll T^2 \label{weyl}
\end{align}
where $|t_j| \sim T$ denotes that $T/2 \leq |t_j| < 2T$. By diadically dividing the integral $[0,T]$ we can also get that
\bal
\sum_{|t_j| < T} 1 = \sum_{m=0}^\infty \lt( \sum_{|t_j|\sim 2^{-(2m+1)}T}1 \rt) \ll \sum_{m=0}^\infty \lt(\frac{T}{2^{2m+1}}\rt)^2 \ll T^2.  \label{weyl2}
\end{align}
Combining \eqref{youngone},\eqref{youngtwo} and \eqref{weyl2}, we can prove the following proposition:
\begin{proposition} \label{ellind} Let $V_{f,\ell}=y^{\kt}\ol{f}\mu_{\ell,k}$, where $f$ is a weight $k$ holomorphic form on $\Gamma=SL_2(\mt{Z})$ and $\mu_{\ell,k}$ is a weight $k$ Maass cusp form on $\Gamma$ as in \eqref{isometry}. Given an orthonormal basis, $u_j$, of Maass forms for $L^2(\Gamma \backslash \mt{H})$, we have that
\beq
\sum_{|t_j|< T} |\rho_j(-1)\LA V_{f,\ell},u_j \RA|^2 + \int_{-T}^T e^{\pi|t|} |\LA V_{f,\ell}, E(*,1/2+it)\RA |^2 \ dt \ll_f e^{\pi|t_\ell|}T^{2k+2\beta_1}\log(T)^6. \label{jutila}
\eeq
when $T \gg 1$ and the implied constant is independent of $\ell$.
\end{proposition}
\begin{proof}
First consider the sum 
\beq
\sum_{|t_j|\sim T} |\rho_j(-1)\LA V_{f,\ell},u_j \RA|^2,
\eeq 
this satisfies, by \eqref{youngone} and the fact that $\rho_j(-1)\sim e^{\pif |t_j|}\log(1+|t_j|)$, the inequality
\bal
& \sum_{|t_j|\sim T} |\rho_j(-1)\LA V_{f,\ell},u_j \RA|^2  \\  
& \ll_f \sum_{|t_j|\sim T} (1+|t_\ell+t_j|)^{k-1+\beta_1}(1+|t_\ell-t_j|)^{k-1+\beta_1}e^{-\pi (||t_j|-|t_\ell||-|t_j|)} \log(1+|t_\ell|)^2\log(1+|t_j|)^4. \notag
\end{align}
Supposing that $|t_\ell| \leq 4T$, then 
\bal
 \sum_{|t_j|\sim T}& (1+|t_\ell+t_j|)^{k-1+\beta_1}(1+|t_\ell-t_j|)^{k-1+\beta_1}e^{-\pi (||t_j|-|t_\ell||-|t_j|)}\log(1+|t_j|)^4\log(1+|t_\ell|)^2 \notag \\
 &  \leq  e^{\pi |t_\ell|} \sum_{|t_j|\sim T} (1+|t_j|+4T)^{2k-2+2\beta_1}\log(T)^6 \ll e^{\pi |t_\ell|} T^{2k+2\beta_1}\log(T)^6
\end{align}
using Weyl's law. If we suppose that $|t_\ell|>4T$ then
\bal
&\sum_{|t_j|\sim T} (1+|t_\ell+t_j|)^{k-1+\beta_1}(1+|t_\ell-t_j|)^{k-1+\beta_1}e^{-\pi (||t_j|-|t_\ell||-|t_j|)} \log(1+|t_j|)^4\log(|t_\ell|)^2 \\
& \leq (1+2|t_\ell|)^{2k-2+2\beta_1}\log(|t_\ell|)^2 \sum_{|t_j|\sim T} \log(1+|t_j|)^4 \leq (1+2|t_\ell|)^{2k-2+2\beta_1}\log(|t_\ell|)  T^{2} \log(T)^4 \notag \\
&\ll e^{\pi |t_\ell|}T^{2k+2\beta_1}\log(T)^6 \notag
\end{align}
where the implied constant is determined by the weight, $k$. Breaking up the interval $[0,T]$ diadically and taking the limit as $T \to \infty$, we get the desired bound for the sum in \eqref{jutila}. The argument for the integral is nearly identical, only using \eqref{youngtwo} instead of \eqref{youngone}. 
\end{proof}
Recalling our construction of $M_k(s,t,\delta)$ and keeping Proposition {\ref{ybound} in mind, we define $\mc{I}_{f,\ell}(s;h,\delta)$ by 
\beq
 \mc{I}_{f,\ell}(s;h,\delta) = \mc{I}_{f,\ell}^{\mbox{\tiny cusp}}(s;h,\delta)+\mc{I}_{f,\ell}^{\mbox{\tiny con}}(s;h,\delta)-\lim_{Y\rightarrow \infty} \LA  P_{h,Y}^{(2)}(*;s;\delta),1 \RA \label{dless}
\eeq
where 
 \beq
  \mc{I}_{f,\ell}^{\mbox{\tiny cusp}}(s;h,\delta)= \sum_{j=1}^\infty \sqrt{2\pi}\frac{\ol{\rho_j(-h)}}{\mcV (2\pi h)^{s-\hf}}M_{0}(s,t_j,\delta) \ol{\LA V_{f,\ell}, u_j \RA} \label{sum2}
\eeq
and
\beq
 \mc{I}_{f,\ell}^{\mbox{\tiny con}}(s;h,\delta)=\frac{1}{4\pi} \int_{-\infty}^\infty \sqrt{2\pi} \frac{h^{-it}\sigma_{2it}(h)}{\mcV \zeta^*(1-2it)(2\pi h)^{s-\hf}} M_{0}(s,t,\delta) \ol{\LA V_{f,\ell}, E(*,1/2+it) \RA} \ dt \label{int2}
\eeq
when $\re s > \hf$, since we can take $\im t_j=0$ for all $t_j$ in $SL_2(\mt{Z})$. We recall from \eqref{mhere} that 
\bal
\label{extra}\LA  P_{h,Y}^{(2)}(*;s;\delta),1 \RA= \frac{\ol{\rho_{\ell,k}(-1)}}{\mcV (2\pi)^{s+\kt-1}} \sum_{h>m\geq 1}  \frac{a(h-m)\ol{\lambda_\ell(m)}}{m^{s+(k-1)/2}} M_{Y,m,-k}(s+\tkt,t_\ell,\tfrac{h\delta}{m}). 
 \end{align} 
From Proposition \ref{ybound}, we know that for $\re s > \hf-\kt+|\im t_\ell | $, the function $M_{Y,m,-k}(s+\tkt,t_\ell,\tfrac{h\delta}{m})$  above converges absolutely and locally uniformly to $M_{-k}(s+\tkt,t_\ell,\tfrac{h\delta}{m})$ as $Y \rightarrow \infty$. 

By the upper bound \eqref{prop1} of Proposition \ref{props} for bounded $\re s$ and $\delta$, and $s$ at least a distance $\vep>0$ from the poles, $M_k(s,t,\delta)$ decays faster than any power of $|t|$ as $|t| \rightarrow \infty$. This together with (\ref{jutila}) and (\ref{weyl}) we know that the right-hand side of (\ref{sum2}) converges absolutely and locally uniformly, as does the integral in (\ref{int2}). The same argument, together with Proposition \ref{ybound}, shows that choosing $A$ sufficiently large will ensure that the expressions for $\mc{I}_{f,\ell, Y,\delta}^{\mbox{\tiny cusp}}(s;h)$ and $\mc{I}_{f,\ell, Y,\delta}^{\mbox{\tiny con}}(s;h)$ converge absolutely and uniformly for any fixed $Y$. Since the differences, $|\mc{I}_{f,\ell}^{\mbox{\tiny cusp}}(s;h,\delta)-\mc{I}_{f,\ell, Y,\delta}^{\mbox{\tiny cusp}}(s;h)|$ and $ |\mc{I}_{f,\ell}^{\mbox{\tiny con}}(s;h,\delta)-\mc{I}_{f,\ell, Y,\delta}^{\mbox{\tiny con}}(s;h)|$ converge absolutely and vanish as $Y \to \infty$, provided $\re s > \hf$ , allowing the interchange of the limit and the sum. Thus 
\beq 
\mc{I}_{f,\ell}^{\mbox{\tiny cusp}}(s;h,\delta)+\mc{I}_{f,\ell}^{\mbox{\tiny con}}(s;h,\delta)=\lim_{Y \rightarrow \infty} (\mc{I}_{f,\ell, Y,\delta}^{\mbox{\tiny cusp}}(s;h)+\mc{I}_{f,\ell, Y,\delta}^{\mbox{\tiny con}}(s;h)) \label{lim1}
\eeq
when $\re s > \hf  $.Thus we get that for $\re s>\hf$,
 \bal
 \label{lim2}\lim_{Y \rightarrow \infty}& \mc{I}_{f,\ell, Y,\delta}(s;h) =  \mc{I}_{f,\ell}(s;h,\delta).
 \end{align}
 
Now it is not difficult to show that $\mc{I}_{f,\ell}(s;h,\delta)$ has a meromorphic continuation to all $s \in \mt{C}$, indeed we already have a meromorphic continuation of $M_k$ and absolute and locally uniform convergence of \eqref{sum2} and \eqref{int2} so we only have to understand the meromorphic continuation of $\mc{I}_{\ell,\delta}^{\mbox{\tiny{con}}}(s;h)$ as we move $s$ past poles in $z$. Rewriting \eqref{int2}, we we see that when $\re s > \hf$,
 \beq
 \mc{I}_{f,\ell}^{\mbox{\tiny con}}(s;h,\delta)=\frac{(2\pi)^{1-s}}{\mcV h^{s-\hf}}\lt(\frac{1}{2\pi i } \int_{(0)}  \frac{h^{-z}\sigma_{2z}(h)}{2\zeta^*(1-2z)\zeta^*(1+2z)} M_{0}(s,z/i,\delta) \ol{\LA V_{f,\ell}, E^*(*,1/2+z) \RA} \ dz \rt)
\eeq
 where $E^*(s,\hf+z)=\zeta^*(1+2z)E(s,\hf+z)$. From \eqref{prop1} we have that the above integral converges for all $s \in \mt{C}$ except for when $\re s -\hf \in \mt{Z}_{\leq 0}.$ Since we start when $\re s > \hf$, to meromorphically continue  $\mc{I}_{f,\ell}^{\mbox{\tiny con}}(s;h,\delta)$ further left we need to repeatedly shift the line of integration. Take $s$ to be in the region $\hf <\re s <  \hf +C$, where $C$ is a curve of the form $C(x)=\frac{c}{\log(2+|x|)}+ix$ for some constant $c>0$ such that $\zeta^*(2-2s)$ has no zeros in this region. We can then deform the line of integration from $\re z =0$ to  $\wt{C}$, a curve from $-i \infty$ to $i \infty$ between $-\ol{C}$ and $C$  such that both the poles of $M_0(s,z/i,\delta)$ at $z=\pm( \hf-s)$ are passed over to get
  \begin{align}
 \mc{I}_{f,\ell}^{\mbox{\tiny con}}(s;h,\delta)&=  \lt.\frac{1}{2\pi i } \int_{\wt{C}}  \frac{(2\pi)^{1-s}h^{\hf-s-z}\sigma_{2z}(h)}{2\mcV \zeta^*(1-2z)\zeta^*(1+2z)} M_{0}(s,z/i,\delta) \ol{\LA V_{f,\ell}, E^*(*,1/2+z) \RA} \ dz\rt.  \label{adove} \\
 & \lt.  + \frac{(4\pi)^{1-s}\sigma_{1-2s}(h)\Gamma(2s-1)}{\zeta^*(2s)\zeta^*(2-2s)\Gamma(s)} \ol{\LA V_{f,\ell}, E^*(*,1-s) \RA} \lt[1+\mc{O}_{\sigma,b}\lt((1+|r|)^{1-2\sigma-b}\delta\rt)  \rt]  \rt. \notag 
\end{align}
where $s=\sigma+ir$. The error term is computed using \eqref{zeepoles} and has a meromorphic continuation to all $\sigma< \hf+\vep$; it has the same poles as the residual term, though the residues of these poles will vanish as $\delta \to 0$. We see that \eqref{adove} holds for $ \hf -\ol{C} < \sigma< \hf +C$, except when $s=\hf$, when the curve $\wt{C}$ cannot be defined. In this exceptional instance we see it is a simple matter to instead shift the line of integration to $-\ol{C}$ and we will only get half the residual term; since this residue and the integral do not have a pole at $s=\hf$, we are able to safely ignore this exception.

Once we shift $\sigma$ left past $\hf$, we can then shift the line of integration back to $\re z =0$ to get the formula 
 \begin{align}
 \mc{I}_{f,\ell}^{\mbox{\tiny con}}(s;h,\delta)&=  \lt.\frac{1}{2\pi i } \int_{(0)}  \frac{(2\pi)^{1-s}h^{\hf-s-z}\sigma_{2z}(h)}{2\mcV \zeta^*(1-2z)\zeta^*(1+2z)} M_{0}(s,z/i,\delta) \ol{\LA V_{f,\ell}, E^*(*,1/2+z) \RA} \ dz\rt.  \notag \\
 & \lt.  + \frac{(4\pi)^{1-s}\sigma_{1-2s}(h)\Gamma(2s-1)}{\zeta^*(2s)\zeta^*(2-2s)\Gamma(s)} \ol{\LA V_{f,\ell}, E^*(*,1-s) \RA} \lt[1+\mc{O}_{\sigma,b}\lt((1+|r|)^{1-2\sigma-b}\delta\rt)  \rt]  \rt. \notag 
\end{align}
which holds when $-\thf < \re s < \thf$. The above makes use of the fact that $h^{-z}\sigma_{2z}(h)=h^{z}\sigma_{-2z}(h)$ as well as the functional equation $E^*(z,s)=E^*(z,1-s)$. This argument can be reproduced as we pass over every line for the form $\re s=\hf-m$ and so we get that  when $\hf -m-\ol{C} < \re s< \hf-m+C$   
  \begin{align}
& \mc{I}_{f,\ell}^{\mbox{\tiny con}}(s;h,\delta)= \lt.\frac{1}{2\pi i } \int_{\wt{C}}  \frac{(2\pi)^{1-s}h^{\hf-s-z}\sigma_{2z}(h)}{2\zeta^*(1-2z)\zeta^*(1+2z)} M_{0}(s,z/i,\delta) \ol{\LA V_{f,\ell}, E^*(*,1/2+z) \RA} \ dz\rt.  \notag \\
 & \lt.+\sum_{n=0}^m \frac{(4\pi)^{1-s}(-1)^n h^{n}\sigma_{1-2s-2n}(h)\Gamma(2s+n-1)}{2\mcV n!\zeta^*(2s+2n)\zeta^*(2-2s-2n)\Gamma(s+n)} \ol{\LA V_{f,\ell}, E^*(*,1-s-n)) \RA}  \rt. \notag \\
 & \ \ \ \ \ \ \ \times  \lt[\frac{\Gamma(1-s)}{\Gamma(1-s-n)}+\mc{O}_{\sigma,b}\lt((1+|r|)^{1-2\sigma-b}\delta\rt)  \rt]. \label{thetag1}
\end{align}
Again there are exceptions when $s-\hf \in \mt{Z}_{\leq 0}$ where instead we integrate along $-\ol{C}$, get half the new residual term at $n=m$, and observe that the poles at $s=\hf-m$ are due to the $n \neq m$ terms. Thus these exceptions at points are also safely ignorable. 
When $-\hf -m < \re s < \hf -m$, 
  \begin{align}
& \mc{I}_{f,\ell}^{\mbox{\tiny con}}(s;h,\delta)= \lt.\frac{1}{2\pi i } \int_{(0)}  \frac{(2\pi)^{1-s}h^{\hf-s-z}\sigma_{2z}(h)}{2\zeta^*(1-2z)\zeta^*(1+2z)} M_{0}(s,z/i,\delta) \ol{\LA V_{f,\ell}, E^*(*,1/2+z) \RA} \ dz\rt.  \notag \\
 & \lt.+\sum_{n=0}^m \frac{(4\pi)^{1-s}(-1)^n h^{n}\sigma_{1-2s-2n}(h)\Gamma(2s+n-1)}{\mcV n!\zeta^*(2s+2n)\zeta^*(2-2s-2n)\Gamma(s+n)} \ol{\LA V_{f,\ell}, E^*(*,1-s-n)) \RA}  \rt. \notag \\
 & \ \ \ \ \ \ \  \times  \lt[\frac{\Gamma(1-s)}{\Gamma(1-s-n)}+\mc{O}_{\sigma,b}\lt((1+|r|)^{1-2\sigma-b}\delta\rt)  \rt]. \label{thetag2}
\end{align}
In all cases, we observe that the integral component remains absolutely convergent so long as $\delta > 0$, and so we have a meromorphic continuation of $ \mc{I}_{f,\ell}^{\mbox{\tiny con}}(s;h,\delta)$, and thus $\mc{I}_{f,\ell}(s;h,\delta)$, to all $s \in \mt{C}$.  We also observe that there are poles of $ \mc{I}_{f,\ell}^{\mbox{\tiny con}}(s;h,\delta)$ due to the terms 
\beq
\frac{\Gamma(2s+n-1)}{\zeta^*(2s+2n)\Gamma(s+n)}.
\eeq
 More specifically we see that the uncanceled poles due to the $\Gamma(2s+n-1)$ are at $s=\hf -m$ for $m \in \mt{Z}_{>0}$ and that the poles due to $\zeta^*(2s+2n)$ are of the form $s=\frac\varrho2-n$ where $\varrho$ is a nontrivial zero of $\zeta^*(s)$. 

 For notational simplicity, let 
 \begin{align}
 \Omega_{f,\ell}(s;h,\delta) =&\lt.\sum_{n=0}^{\lf \hf -\sigma \rf} \frac{(4\pi)^{1-s}(-1)^n h^{n}\sigma_{1-2s-2n}(h)\Gamma(2s+n-1)}{\mcV n!\zeta^*(2s+2n)\zeta^*(2-2s-2n)\Gamma(s+n)} \ol{\LA V_{f,\ell}, E^*(*,1-s-n)) \RA}  \rt. \notag \\
 & \ \ \ \ \ \ \  \times  \lt[\frac{\Gamma(1-s)}{\Gamma(1-s-n)}+\mc{O}_{\sigma,b}\lt((1+|r|)^{1-2\sigma-b}\delta\rt)  \rt],
 \end{align}
which we observe is piecewise-meromorphic. Using (\ref{dless}),(\ref{sum2}),(\ref{int2}), \eqref{extra}, we have the expansion
\begin{subequations} \label{blarrr}
\bal
\mc{I}_{f,\ell}(s;h,\delta) \ \  \  =& \lt. \sum_{j=1}^\infty \frac{(2\pi)^{1-s}\ol{\rho_j(-h)}} {\mcV h^{s-\hf}}M_0(s,t_j,\delta)\ol{\LA V_{f,\ell}, u_j\RA} \rt. \label{spec0} \\
& +\lt.\frac{1}{2\pi i } \int_{C_\sigma}  \frac{(2\pi)^{1-s}h^{\hf-s-z}\sigma_{2z}(h)}{2\zeta^*(1-2z)\zeta^*(1+2z)} M_{0}(s,z/i,\delta) \ol{\LA V_{f,\ell}, E^*(*,1/2+z) \RA} \ dz\rt.  \label{specy1} \\
&+\Omega_{f,\ell}(s;h,\delta) \label{specy4} \\
&\lt. -\frac{\ol{\rho_{\ell,k}(-1)}}{\mcV (2\pi)^{s+\kt-1}} \sum_{ m =1}^h  \frac{a(h-m)\ol{\lambda_\ell(m)}}{m^{s+(k-1)/2}}  M_{-k}(s+\tkt,t_\ell,\tfrac{h\delta}{m}) \rt. \label{specy3} 
\end{align}
\end{subequations}
for all $s \in \mt{C}$ where
  \beq \label{csigma}
 C_\sigma = \lt\{
 \begin{array}{ll}
 (0) & \mbox{ when } \thf-\sigma \notin \mt{Z}_{\geq 0} \\
\wt{C} & \mbox{ when } \thf -\sigma \in \mt{Z}_{\geq 0}.
 \end{array}
\rt.
\eeq

Now, using the representation of $ \mc{I}_{f,\ell, Y,\delta}(s;h)$ in (\ref{intend}) we see that 
\begin{subequations}
\bal
&\lim_{Y \rightarrow \infty}  \mc{I}_{f,\ell, Y,\delta}(s;h) \notag \\
& = \frac{\ol{\rho_{\ell,k}(1)}}{(2\pi)^{s+\kt-1} } \sum_{m> 0} \frac{a(m+h)\ol{\lambda_\ell(m)}}{m^{s+(k-1)/2}}  \int_{0}^{\infty} e^{-y(1+h\delta/m)}W_{\kt,it_\ell}(2y)y^{s+\kt-1} \frac{dy}{y}.  \label{lub}\\
&= \frac{\ol{\rho_{\ell,k}(1)}\Gamma(s+(k-1)/2+it_\ell)\Gamma(s+(k-1)/2-it_\ell)}{\mcV (4\pi)^{s+\kt-1}\Gamma(s)} D_{f,\ell}^+(s;h,\delta) \label{lub2}
\end{align}
\end{subequations}
where 
\begin{align}
&D_{f,\ell}^+(s;h,\delta)= \sum_{m >0} \frac{a(m+h)\ol{\lambda_\ell(m)}}{m^{s+\frac{k-1}{2}}}\lt(1+\frac{\delta h}{2m}  \rt)^{-(s+\frac{k-1}{2}+it_\ell)}F\lt(s+\tfrac{k-1}{2}+it_\ell,\tfrac{1-k}{2}+it_\ell;s;\tfrac{h\delta}{2m+h\delta}\rt) \notag  \\
&=\sum_{m >0} \frac{a(m+h)\ol{\lambda_\ell(m)}}{m^{s+\frac{k-1}{2}}}\lt(1+\frac{\delta h}{2m}  \rt)^{-(s+\frac{k-1}{2}+it_\ell)} \lt[1+ \mc{O}_{\sigma,k,t_\ell}\lt((1+|r|)^{\sigma+k}\lt(\frac{h\delta}{2m+h\delta}\rt) \rt)\rt] \label{done}
\end{align}
converges absolutely and locally uniformly when $\re s >1$. The error term is obtained by uncomplicated analysis of \eqref{bunk} via Stirling's approximation.  We have now established that the function  $\mc{I}_{f,\ell}(s;h,\delta)$ defined by (\ref{dless}) has a meromorphic continuation comprised of absolutely convergent sums and integrals  for all $s \in \mt{C}$. Since $\mc{I}_{f,\ell}(s;h,\delta)$ is related to $D_{f,\ell}^+(s;h,\delta)$ by 
\bal 
D_{f,\ell}^+(s;h,\delta) \label{mer2}= \frac{\mcV (4\pi)^{s+\kt-1}\Gamma(s)}{\ol{\rho_{\ell,k}(1)}\Gamma(s+(k-1)/2+it_\ell)\Gamma(s+(k-1)/2-it_\ell)}  \mc{I}_{f,\ell}(s;h,\delta)
\end{align}
when $\re s >1$, the meromorphic continuation of $D_{f,\ell}^+(s;h,\delta)$ is related to that of $\mc{I}_{f,\ell}(s;h,\delta)$ via the above formula. We also note that the gamma functions in the denominator of \eqref{mer2} cancel out the poles due to $M_{-\kt}$ in \eqref{spec0}. 

By (\ref{prop1}), (\ref{jutila}), and (\ref{dless}) we know that the series expressions given for $ \mc{I}_{f,\ell}(s;h,\delta)$ and $D_{f,\ell}^+(s;h,\delta)$ have meromorphic continuations with absolutely and locally uniformly convergent sums and integrals when $A>2+|\sigma|+\frac{k}{2},$ but the upper bound can have a factor of $\delta^{-A}$ in it. As we hope to eventually let $\delta \rightarrow 0$, it is useful to identify an interval in for $\re s$  of absolute and uniform convergence where we have an upper bound that is independent of $\delta$.  

Combining (\ref{ll1}) and (\ref{ll2}) we have, for $ \sigma < 1-\kt$ at least $\vep$ away from the poles of $\mc{I}_{f,\ell}(s;h,\delta)$, regardless of the relation between $\delta$ and $(1+|t|)^{-2}$, when $A>2+|\sigma|+\frac{k}{2},$
\beq
\label{ll3} M_k(s,t,\delta) \ll (1+|t|)^{2\sigma-2+k+2\ep}(1+|r|)^{9A+2}e^{-\frac{\pi}{2}|r|}.
\eeq 
Using \eqref{weyl}, Proposition \ref{ellind} and the Cauchy-Schwarz inequality, we have for any $T \gg 1$,
\beq 
\lt| \sum_{|t_j|< T} \ol{\rho_j(-h)\LA V_{f,\ell},u_j \RA }\rt| \ll (T^2)^{\hf}\lt(  \sum_{|t_j|\sim T}\lt|\rho_j(-h)\LA V_{f,\ell},u_j \RA \rt|^2  \rt)^{\hf} \ll_f h^{\theta} T^{k+\beta_1+1}\log(T)^3e^{\frac{\pi}{2}|t_\ell|} \label{jutila2}
\eeq
where the implied constant is independent of $\ell$ and $\theta$ is the best-known constant such that 
$\lambda_\ell(h) \ll h^\theta$. Combining this with (\ref{ll3}) we see that 
\beq
\lt| \sum_{|t_j| \sim T} \ol{\rho_j(-h)}M_0(s,t_j,\delta)\ol{\LA V_{f,\ell},u_j\RA} \rt| \ll_{f,A,b,\vep} h^{\theta}T^{2\sigma +2\ep+k+\beta_1-1}(1+|r|)^{9A+2}e^{-\frac{\pi}{2}(|r|-|t_\ell|)}. \label{cbound}
\eeq 
Breaking up the interval from $1$ to $T$ diadically and letting $T \rightarrow \infty$ it is clear from 
(\ref{cbound}) that the cuspidal part of the spectral expansion of $D_{f,\ell}^+(s;h,\delta)$ converges absolutely and locally uniformly, independent of $\delta$, when $\sigma <\hf-\kt-\tfrac{\beta_1}{2}$. Similar reasoning gives us
\beq
\int_{-T}^T  \frac{h^{-it}\sigma_{2it}(h)}{\zeta^*(1-2it)} M_0(s,t,\delta) \ol{\LA V_{f,\ell}, E(*,\thf+it) \RA} \ dt \ll_{f,A,b,\vep} h^\vep T^{2\sigma +2\vep+k+\beta_1-2} (1+|r|)^{9A+2}e^{-\frac{\pi}{2}(|r|-|t_\ell|)} \label{cobound}
\eeq 
when $s \neq \hf-r$ for $r \geq 0$. When $ \hf -r-\ol{C} \leq \re s \leq \hf-r+C$, we see that using \eqref{youngtwoone} gives us that
\begin{align}
& \int\limits_{\wt{C}_T}  \frac{h^{-z}\sigma_{2z}(h)}{\zeta^*(1-2z)} M_0(s,z/i,\delta) \ol{\LA V_{f,\ell}, E(*,\thf+z) \RA} \ dz \label{cobound2}\\
& \ \ \ \ \ \ \ \ \ \ \ \ \ \   \ll_{f,A,b,\vep} h^\vep T^{2\sigma +2\vep+k+\beta_1-2} (1+|r|)^{10A-3b+\vep}e^{-\frac{\pi}{2}(|r|-|t_\ell|)}, \notag
\end{align} 
where $\wt{C}_T$ is just $\wt{C}$ but stops at $|\im z|=T$. So for $\sigma< \hf-\kt-\tfrac{\beta_1}{2}-\vep$ we see that these integrals are well-defined as $T \rightarrow \infty$ for the domains in $s$ where they occur. So we have that the continuous part of the spectral expansions of  $\mc{I}_{f,\ell}(s;h,\delta)$ and $D_{f,\ell}^+(s;h,\delta)$ converge in the same region as the cuspidal part. 

Thus \eqref{spec0},(\ref{cbound}), (\ref{cobound}), \eqref{cobound2},  \eqref{ll3} for $M_{-k}$, and Stirling's formula provide the upper bounds for $\mc{I}_{f,\ell}(s;h,\delta)$ and $D_{f,\ell}^+(s;h,\delta)$ given in Proposition \ref{prop2} below, which also summarizes the above discussion.

\begin{proposition} \label{prop2} The functions $\mc{I}_{f,\ell}(s;h,\delta)$ and $D_{f,\ell}^+(s;h,\delta)$ defined in (\ref{dless}),(\ref{done}), and (\ref{mer2}) have meromorphic continuations to all $s \in \mt{C}$.   Specifically $\mc{I}_{f,\ell}(s;h,\delta)$ has simple poles at $s=\hf + it_j-r$ for all $j \in \mt{Z}_{\neq 0}$ where $r \in \mt{Z}_{\geq 0}$, as well as at $s=\hf-r$ for $r \in \mt{Z}_{>0}$. It also has poles when $\tfrac{\vho}{2}-r$ for $r \in \mt{Z}_{\geq 0}$, where $\vho$ is a nontrivial zero of the zeta function. From this, the function $D_{f,\ell}^+(s;h,\delta)$ has poles with the same order at the same locations, excepting when $t_j=t_\ell$ for $r\geq \kt$. In addition to this, $D_{f,\ell}^+(s;h,\delta)$ has simple poles at $s=-r$, where $r \in \mt{Z}_{\geq 0}$. 

In particular, letting
\beq
\wt{R}_{f,\ell}(s_0;h,\delta):= \stackrel[s=s_0]{}{\mbox{\emph{Res}} }\mc{I}_{f,\ell}(s;h,\delta)
\eeq
we have 
\begin{subequations} \label{rrbound}
\bal
& \wt{R}_{f,\ell}(1/2+it_j-r;h,\delta) \label{rboundb}\\
&= \frac{(-1)^r (4\pi)^{\hf+r-it_j} h^{r- it_j}\Gamma(\hf-it_j+r)\Gamma(2it_j-r)\ol{\rho_j(-h)\LA V_{f,\ell}, u_j \RA }}{r! \mcV \Gamma(\hf+it_j)\Gamma(\hf-it_j)} \notag \\
& \  \ +\mc{O}_{m,b}\lt( (1+|t_\ell| )^{k-1+\beta_1}(1+|t_j|)^{k-\frac{3}{2}-b+r+\beta_1}\log(1+|t_j|)^4\log(1+|t_\ell|)^2 h^{r+\theta} e^{-\frac{\pi}{2}||t_j|-|t_\ell||}\delta \rt) \notag 
\end{align}
if $t_j \neq \pm t_\ell$ or $r < \kt$, otherwise
\bal
\wt{R}_{f,\ell}(1/2&\pm it_\ell-r;h,\delta) \label{rbound} \\
&=\frac{(-1)^r (4\pi)^{\hf+r\mp it_\ell}h^{r\mp it_\ell} \Gamma(\hf \mp it_\ell+r)\Gamma(\pm 2it_\ell-r)\ol{\rho_\ell(-h)\LA V_{f,\ell}, u_\ell \RA }}{r! \mcV \Gamma(\hf+it_\ell)\Gamma(\hf-it_\ell)}  \notag \\ 
&- \frac{(-1)^{r-\kt}\ol{\rho_{\ell,k}(-1)}(4\pi)^{\hf+r-\kt\mp it_\ell}\Gamma(\thf \mp it_\ell+r)\Gamma(\pm 2it_\ell-r+\tkt)}{\mcV (r-\tkt)!\Gamma(\thf+\tkt+it_\ell)\Gamma(\thf+\tkt-it_\ell)}\sum_{m=1}^h \frac{a(h-m)\ol{\lambda_\ell(m)}}{m^{\kt-r\pm it_\ell}} \notag \\
& +\mc{O}\lt( (1+|t_\ell| )^{r+\kt-\hf-b}\log(1+|t_\ell|)h^{\hf+r+\theta}\delta \rt). \notag 
\end{align}
\end{subequations}
We note that the zeros that occur when $t_\ell=t_j$ and $r \geq \tkt$ are cancelled out in $D_{f,\ell}^+(s;h,\delta)$ by the accompanying $\Gamma(s+\tfrac{k-1}{2}\pm it_\ell)$ factors in the denominator.  


Let $A> 1+\kt+|\sigma|+\beta_1$.  Then for $s$ a distance at least $\vep>0$ from the poles we have the bounds 
\begin{align}
& |\mc{I}_{f,\ell}(s;h,\delta)-\mc{I}_{f,\ell}^{\mbox{\tiny con}}(s;h,\delta)|  \ll_{f,A,\vep} \delta^{-A}(1+|r|)^{9A+2} h^{1-\sigma+\theta} e^{-\frac{\pi}{2}(|r|-|t_\ell|)}, \label{abounds} \\
&|\mc{I}_{f,\ell}(s;h,\delta)-\Omega_{f,\ell}(s;h,\delta)|   \ll_{f,A,\vep} \delta^{-A}(1+|r|)^{9A+2} h^{1-\sigma+\theta}e^{-\frac{\pi}{2}(|r|-|t_\ell|)}. \label{abounds2}
\end{align}
In the region $\re s =\sigma < \hf-\kt-\tfrac{\beta_1}{2}-\vep$, for a distance of at least $\vep>0$ from the poles, the expansion for $\mc{I}_{f,\ell}(s;h,\delta) $ given in \eqref{spec0} converges absolutely and locally uniformly and satisfies the upper bounds
\bal
& \lt| \mc{I}_{f,\ell}(s;h,\delta)-\mc{I}_{f,\ell}^{\mbox{\tiny con}}(s;h,\delta) \rt|  \ll_{f,A,b,\vep}  (1+|r|)^{9A+2}h^{1-\sigma+\theta}(1+|t_\ell|)^{2\sigma-2+2\vep}e^{-\frac{\pi}{2}(|r|-|t_\ell|)}, \label{aless}\\
&|\mc{I}_{f,\ell}(s;h,\delta)-\Omega_{f,\ell}(s;h,\delta)|  \ll_{f,A,b,\vep}  (1+|r|)^{9A+2}h^{1-\sigma+\theta}(1+|t_\ell|)^{2\sigma-2+2\vep}e^{-\frac{\pi}{2}(|r|-|t_\ell|)}. \label{aless2}
\end{align}
\end{proposition}
\begin{proof} 
In Proposition \ref{props}, (\ref{prop1}) gives that for $A>1+\kt+|\sigma|+|\beta_1|$, fixed $\delta>0$ and $s$ at least a distance of $\vep>0$ from the poles, the meromorphic continuation of $M_k(s,t,\delta)$ decays faster than any power of $|t|$. Thus the series expression for $\mc{I}_{f,\ell}(s;h,\delta)$ given in (\ref{dless}) converges absolutely and locally uniformly, giving a meromorphic continuation of $\mc{I}_{f,\ell}(s;h,\delta)$ to all $s \in \mt{C}$, with possible simple poles at the points $s=\hf + it_j-r$ for each $t_j$ and $r \in \mt{Z}_{\geq 0}$. The residues at these points follow from the corresponding residual computations of $M_k(s,t,\delta)$ given in Proposition \ref{props}. The corresponding meromorphic continuation of $D_{f,\ell}^+(s;h,\delta)$ follows from (\ref{mer2}). 

From \eqref{cbound} and \eqref{cobound}, which made use of \eqref{ll3}, we have bounds for \eqref{spec0} and \eqref{specy1}, and we can similarly bound \eqref{specy3} using \eqref{ll3}, Stirling's approximation formula, and known growth properties for Fourier coefficients of holomorphic cusp forms and Maass forms. All of these bounds together give us \eqref{aless} and \eqref{aless2}. By using \eqref{prop1} instead of \eqref{ll3} in these estimates, we get \eqref{abounds} and \eqref{abounds2}. 
\end{proof}

\section{The Limit As $\delta\rightarrow 0$}
The entirety of this section is devoted to proving the following proposition, which concludes this chapter by describing the meromorphic continuation of $D_{f,\ell}^+(s;h)$. 
\begin{proposition} \label{prop5} For  $\sigma = \re s < \hf-\tkt-\tfrac{\beta_1}{2}-\vep$, let
\beq \label{duhh}
\mc{I}_{f,\ell}(s;h):=\lim_{\delta \rightarrow 0} \mc{I}_{f,\ell}(s;h,\delta),
\eeq
which converges uniformly in $s$ and $h$. The function $\mc{I}_{f,\ell}(s;h)$ has a meromorphic continuation to all $s \in \mt{C}$.
When $\sigma > 1$ we have
\beq
D_{f,\ell}^+(s;h):=\lim_{\delta \rightarrow 0} D_{f,\ell}^+(s;h,\delta) = \sum_{m>0} \frac{a(m+h)\ol{\lambda_\ell(m)}}{m^{s+(k-1)/2}}, \label{lhs}
\eeq
and the function $\mc{I}_{f,\ell}(s;h)$ is related to $D_{f,\ell}^+(s;h)$ in this region by
\bal
\mc{I}_{f,\ell}(s;h) \label{mer3}=G_\ell(s)^{-1}D_{f,\ell}^+(s;h). 
\end{align}
where
\beq
G_\ell(s):= \frac{\mcV (4\pi)^{s+\kt-1}\Gamma(s)} {\ol{\rho_{\ell,k}(1)}\Gamma(s+(k-1)/2+it_\ell)\Gamma(s+(k-1)/2-it_\ell)} \label{shro},
\eeq
and thus $D_{f,\ell}^+(s;h)$ also has a continuation to all $s \in \mt{C}$.

The function $\mc{I}_{f,\ell}(s;h)$ has poles when $s=\thf\pm it_j-r$ where $r \in \mt{Z}_{\geq 0}$. Letting 
\beq
\wt{R}_{f,\ell}(s_0;h):=\stackrel[s=s_0]{}{\mbox{\emph{Res}} }\mc{I}_{f,\ell}(s;h) 
\eeq
we see that
\beq
c_{r,j,f,\ell}(h):=\wt{R}_{f,\ell}(\thf+it_j-r;h) = \frac{(-1)^r (4\pi)^{\hf+r-it_j} h^{r-it_j} \Gamma(\hf-it_j+r)\Gamma(2it_j-r)\ol{\rho_j(-h)\LA V_{f,\ell}, u_j \RA }}{r! \mcV \Gamma(\hf+it_j)\Gamma(\hf-it_j)} \label{round2}
 \eeq
 if $t_j =\pm t_\ell$ or $r < \tkt$, otherwise
\bal
&c_{r,\pm \ell,f, \ell}(h)= \frac{ (4\pi)^{\hf+r\mp it_\ell}h^{r\mp it_\ell} \Gamma(\hf \mp it_\ell+r)\Gamma(\pm 2it_\ell-r)\ol{\rho_\ell(-h)\LA V_{f,\ell}, u_\ell \RA }}{(-1)^rr! \mcV \Gamma(\hf+it_\ell)\Gamma(\hf-it_\ell)} \notag \\ 
& \mkern20mu - \frac{(-1)^{r-\kt}\ol{\rho_{\ell,k}(-1)}(4\pi)^{\hf+r-\kt\mp it_\ell}\Gamma(\thf \mp it_\ell+r)\Gamma(\pm 2it_\ell-r+\tkt)}{\mcV (r-\tkt)!\Gamma(\thf+\tkt+it_\ell)\Gamma(\thf+\tkt-it_\ell)}\sum_{m=1}^h \frac{a(h-m)\ol{\lambda_\ell(m)}}{m^{\kt-r\pm it_\ell}}.  \label{rbound3}  
\end{align}
These $c_{r,j,f,\ell}(h)$ satisfy the average upper bound
 \beq
 \sum_{j \neq 0} |c_{r,j,f,\ell}(h)|^2 \ll (1+|t_\ell|)^{k+1+\beta_1+\vep}h^{1+2r+2\theta}. \label{rbound4}
 \eeq
Furthermore $\mc{I}_{f,\ell}(s;h)$ also has poles due to $\Omega_{f,\ell}(s;h)$, which is given by,
    \beq
\Omega_{f,\ell}(s;h) :=\sum_{n=0}^{\lf \hf -\sigma \rf } \frac{(4\pi)^{1-s}(-1)^n h^{n}\sigma_{1-2s-2n}(h)\Gamma(2s+n-1)\Gamma(1-s)}{n!\zeta^*(2s+2n)\zeta^*(2-2s-2n)\Gamma(s+n)\Gamma(1-s-n)} \ol{\LA V_{f,\ell}, E^*(*,s+n)) \RA} \label{omega14}
\eeq
where we take the empty sum to be zero. Thus $\mc{I}_{f,\ell}(s;h)$ also has simple poles at $\hf -r$ for $r \in \mt{Z}_{>0}$ and not necessarily simple poles at $\frac \varrho 2 - r$ for $r \in \mt{Z}_{\geq 0}$ and $\varrho$ is any nontrivial zero of $\zeta(s)$. The poles of $\mc{I}_{f,\ell}(s;h)$ are also poles of $D_{f,\ell}^+(s;h)$,  except when $s=\hf \pm it_\ell -r$ for $r \in \mt{Z}_{\geq \kt}$. All coinciding poles are of the same order.
 The function $D_{f,\ell}^+(s;h)$ also has simple poles when $s \in \mt{Z}_{\leq 0}$. 


For $s$ at least $\vep>0$ away from the polar points, and 
$\sigma<\hf-\tkt-\frac{\beta_1}{2}$, $\mc{I}_{f,\ell}(s;h)$ is given by the following convergent sums and integral:
\bal
\mc{I}_{f,\ell}(s;h) =& \sum_{j=1}^\infty  \frac{(4\pi)^{1-s} \ol{\rho_j(-1)\lambda_j(h)}h^{\hf-s}\Gamma(s-\thf+it_j)\Gamma(s-\thf-it_j)\Gamma(1-s)\ol{\LA V_{f,\ell},u_j \RA }}{\mcV \Gamma(\hf+it_j)\Gamma(\hf-it_j)} \label{spec3} \\
&+ \frac{1}{2\pi i } \int\limits_{C_\sigma}  \frac{(4\pi)^{1-s}\sigma_{2z}(h)h^{\hf-s-z}\Gamma(s-\hf-z)\Gamma(s-\hf+z)\Gamma(1-s)\ol{\LA V_{f,\ell}, E^*(*,\hf+z) \RA }}{2\mcV \zeta^*(1+2z)\zeta^*(1-2z)\Gamma(\hf+z)\Gamma(\hf-z)} \ dz \notag  \\ 
&  - \frac{\ol{\rho_{\ell,k}(-1)}\Gamma(s+\tfrac{k-1}{2}+it_\ell)\Gamma(s+\tfrac{k-1}{2}-it_\ell)\Gamma(1-s)}{\mcV (4\pi)^{s+\kt-1}\Gamma(\hf+\kt+it_\ell)\Gamma(\hf+\kt-it_\ell)} \sum_{m = 1}^h \frac{a(h-m)\ol{\lambda_\ell(m)}}{m^{s+\frac{k-1}{2}}} \notag \\
& + \Omega_{f,\ell}(s;h) \notag 
\end{align} 
where $C_\sigma$ is as in \eqref{csigma}. 

For $s$ at least $\vep>0$ away from the points $\hf +it_j -r$, and in the vertical strip $c < \sigma <2$, for any $c$ with $-\thf -\tkt <c< \hf-\tkt-\frac{\beta_1}{2}$, $\mc{I}_{f,\ell}(s;h)$ is given by
\bal
\mc{I}_{f,\ell}(s;h)&=\label{anconint} \frac{1}{2\pi i} \int_{2-i\infty}^{2+i\infty} \frac{\mc{I}_{f,\ell}(u;h)}{(u-s)} \ du -  \frac{1}{2\pi i} \int_{c-i\infty}^{c+i\infty} \frac{(\mc{I}_{f,\ell}(u;h)-\Omega_{f,\ell}(u;h))}{(u-s)} \ du  \\
&-\sum_{r=0}^\kt \int_{C'_\sigma} \frac{(-1)^r(4\pi)^{\hf+z+r}\sigma_{-2z}(h)h^r\Gamma(-2z-r)\Gamma(\hf+z+r)}{r! \mcV  \zeta^*(1+2z)\zeta^*(1-2z)\Gamma(\hf+z)\Gamma(\hf-z)(\hf-z-r-s)} \ol{\LA V_{f,\ell}, E^*(*,\thf+z)) \RA} \ dz\notag \\
&+\Omega_{f,\ell}(s;h)-\sum_{j \neq 0}\sum_{r=0}^{k/2} \frac{c_{r,j,f,\ell}(h)}{(\hf+it_j-r-s)}, \notag 
\end{align}
where when $\re u= 2$, $\mc{I}_{f,\ell}(u;h)$ is given by (\ref{mer3}) and when $\re u = c$, $\mc{I}_{f,\ell}(u;h)$ is given by (\ref{spec3}). The curve $C'_\sigma$ is essentially the same as $C_\sigma$ but without the ignorable exception when  $s-\hf \in \mt{Z}_{\geq 0}$.

Let $A > 1 + \kt + |\sigma|+ \beta_1$, then for $\sigma<\hf-\kt-\frac{\beta_1}{2}$ and $s$ at least a distance of $\vep>0$ away from the poles at $s=\hf+it_j -r$, we have the bound
\bal
|\mc{I}_{f,\ell}(s;h)-\Omega_{f,\ell}(s;h) | \ll_{f,A,\vep} (1+|s|)^{10A}h^{1-\sigma+\theta}e^{-\frac{\pi}{2}(|\im s | -|t_\ell|)}. \label{aless3}
\end{align}
For $\sigma>1+\ep$, $\mc{I}_{f,\ell}(s;h)$ satisfies the bound
\bal
\mc{I}_{f,\ell}(s;h) \ll_{f,\vep} h^{(k-1)/2}(1+|t_\ell|)^{2\sigma+k-2}(1+|s|)^{\sigma+k-\frac{3}{2}}e^{-\frac{\pi}{2}||\im s|-|t_\ell||}. \label{easy2}
\end{align}
When $c<\sigma<2$, and $s$ at least a distance of $\vep>0$ away from the poles at $s=\hf+it_j -r$, we have the bound
\beq
|\mc{I}_{f,\ell}(s;h) -\Omega_{f,\ell}(s;h)| \ll_{f,c} h^{1-c+\theta}(1+|t_\ell|)^{k}e^{\pif |t_\ell|} \label{later1}
\eeq
 \end{proposition}
\begin{proof}[Proof of Proposition \ref{prop5}]
Referring to (\ref{ll1}), we see that for any fixed $t$ and $\re s < 1/2$, $\lim_{\delta \rightarrow 0} M_0(s,t,\delta)$ and $\lim_{\delta \rightarrow 0} M_{-k}(s,t,\delta)$ exist with
\beq
M_{0}(s,z/i):=\lim_{\delta \rightarrow 0} M_{0}(s,z/i,\delta) = \frac{2^{1-s}\Gamma(s-\frac{1}{2}-z)\Gamma(s-\hf+z)\Gamma(1-s)}{\Gamma(\frac{1}{2}+z)\Gamma(\hf-z)}
\eeq
\beq
M_{-k}(s,z/i):=\lim_{\delta \rightarrow 0} M_{-k}(s,z/i,\delta) = \frac{2^{1-s}\Gamma(s-\frac{1}{2}-z)\Gamma(s-\hf+z)\Gamma(1-s+\kt)}{\Gamma(\frac{1}{2}+\frac{k}{2}+z)\Gamma(\hf+\kt-z)}.
\eeq
From this, \eqref{blarrr}, and \eqref{aless2}, it is clear that when $s<\hf-\kt -\frac{\beta_1}{2}$, the spectral expansion for $\mc{I}_{f,\ell}(s;h,\delta)$ converges absoultely and uniformly, as $\delta \to 0$, to $\mc{I}_{f,\ell}(s;h)$; this gives us \eqref{duhh} and (\ref{spec3}). We also see that $\mc{I}_{f,\ell}(s;h)$ satisfies the upper bound \eqref{aless3}, which we get by an application of Stirling's approximation and using the same reasoning applied in Proposition \ref{prop2}.  Similarly by \eqref{done} and \eqref{mer2}, we get \eqref{lhs} and \eqref{mer3}. We get \eqref{easy2} via Stirling's approximation and the known growth of Fourier coefficients on average.

It remains to be shown that the functions given by \eqref{duhh} and \eqref{mer3} are actually the same function, that they agree and are well-defined in the region $c <\sigma < 2$. The remainder of the proof is devoted to this.

Fix $\delta>0$. Let $  \mc{I}^{\mbox{\tiny else}}_{f,\ell}(s;h,\delta):= \mc{I}_{f,\ell}(s;h,\delta)-\mc{I}^{\mbox{\tiny con}}_{f,\ell}(s;h,\delta)$, where $\mc{I}^{\mbox{\tiny con}}_{f,\ell}(s;h,\delta)$ is as it it is given in \eqref{thetag1} and \eqref{thetag2}. Using the information in Proposition \ref{prop2}, the Cauchy residue theorem allows us to express components of $\mc{I}_{f,\ell}(s;h,\delta)$ with $c<\sigma < 2 $ as follows. For $T \gg 1$, when $|\im s|<T$,  
  \bal 
  \mc{I}^{\mbox{\tiny else}}_{f,\ell}(s;h,\delta)  =& I _{1,\delta}(T)-I _{2,\delta}(T)-I _{3,\delta}(T)+I _{4,\delta}(T)\\
  &  - \sum_{|t_j|<T} \sum_{r=0}^\kt \frac{ \wt{R}_{f,\ell}(1/2+it_j-r;h,\delta)}{(\hf+it_j-r-s)}    \notag
  \end{align} 
  where $\wt{R}_{f,\ell}(1/2+it_j-r;h,\delta)$ is as in \eqref{rrbound} and
\bal
& I _{1,\delta}(T) = \frac{1}{2\pi i} \int_{2-iT}^{2+iT} \frac{\mc{I}^{\mbox{\tiny else}}_{f,\ell}(u;h,\delta)}{(u-s)} \ du\\
  &I _{2,\delta}(T)= \frac{1}{2\pi i} \int_{c+iT}^{2+iT} \frac{\mc{I}^{\mbox{\tiny else}}_{f,\ell}(u;h,\delta)}{(u-s)} \ du \notag \\ 
&I _{3,\delta}(T) =\frac{1}{2\pi i} \int_{c-iT}^{c+iT} \frac{\mc{I}^{\mbox{\tiny else}}_{f,\ell}(u;h,\delta)}{(u-s)} \ du \notag \\
&I _{4,\delta}(T)= \frac{1}{2\pi i} \int_{c-iT}^{2-iT} \frac{\mc{I}^{\mbox{\tiny else}}_{f,\ell}(u;h,\delta)}{(u-s)} \ du. \notag
\end{align}
Here $T$ is such that the contour lines for $I_{2,\delta}(T)$ and $I_{4,\delta}(T)$ never gets closer to the poles of $\mc{I}^{\mbox{\tiny else}}_{f,\ell}(u;h,\delta)$ at $u=\hf+it_j-r$ than  $\alpha T^{-2}$, for some fixed $\alpha>0$, which is permitted by Weyl's Law. 
Thus from \eqref{nearp1} and (\ref{abounds}) we have that
\bal
I_{2,\delta}(T), I_{4,\delta}(T) \ll_{A,c}&  T^{A'}h^{1-c+\theta}\delta^{-A}(1+(1+|t_\ell|)^{2+2\vep-2A})e^{\pi|t_\ell|}e^{-\frac{\pi}{2}T}\\
& +\alpha^{-3}T^{A'}h^{1-c+\theta}\delta^{-A}e^{-\frac{\pi}{2}T}, \notag
\end{align}
where $A>1+\kt+|\re u|$ and $A'$ is a constant dependent on $A$. So we see that $I_{2,\delta}(T)$ and $I_{4,\delta}(T)$ vanish uniformly as $T \to \infty$. We also see that \eqref{abounds} also gives that $I _{1,\delta}(T)$ and  $I _{3,\delta}(T)$ are absolutely and uniformly convergent as $T \to \infty$. Furthermore, \eqref{cbound}, \eqref{rrbound} and Stirling's approximation give us that the infinite sum over residues is also absolutely and uniformly convergent as $T \to \infty$. Thus for all $s$ such that $ c < \sigma <2$, 
  \bal 
  \mc{I}^{\mbox{\tiny else}}_{f,\ell}(s;h,\delta)  =& \frac{1}{2\pi i} \int_{(2)} \frac{\mc{I}^{\mbox{\tiny else}}_{f,\ell}(u;h,\delta)}{(u-s)} \ du  - \frac{1}{2\pi i} \int_{(c)} \frac{\mc{I}^{\mbox{\tiny else}}_{f,\ell}(u;h,\delta)}{(u-s)} \ du \label{else} \\
  &  - \sum_{j \neq 0} \sum_{r=0}^\kt \frac{ \wt{R}_{f,\ell}(1/2+it_j-r;h,\delta)}{(\hf+it_j-r-s)}.    \notag
  \end{align} 
  Now we observe that for $\frac{1}{2}<\sigma<2$, 
  \bal
   \frac{1}{2\pi i} \int\limits_{(2)} & \frac{\mc{I}^{\mbox{\tiny con}}_{f,\ell}(u;h,\delta)}{(u-s)} \ du  =\\
   & \notag \lt(\frac{1}{2\pi i} \rt)^2\int\limits_{(2)}  \int\limits_{(0)}  \frac{(2\pi)^{1-u}h^{\hf-u-z}\sigma_{2z}(h)M_{0}(u,z/i,\delta)}{2\zeta^*(1-2z)\zeta^*(1+2z)(u-s)}  \ol{\LA V_{f,\ell}, E^*(*,1/2+z) \RA} \ dz du.
  \end{align}
By Stirling's approximation, \eqref{prop1}, and \eqref{youngtwo} the above double integral is absolutely and uniformly convergent. So we can interchange the order of integration and then shift the line of integration of $u$ left from $\re u =2$ to $\re u =c$ to get
 \bal
 &  \frac{1}{2\pi i} \int_{(2)} \frac{\mc{I}^{\mbox{\tiny con}}_{f,\ell}(u;h,\delta)}{(u-s)} \ du  =\mc{I}^{\mbox{\tiny con}}_{f,\ell}(s;h,\delta)  \\
 & \ \ \ \ +\lt(\frac{1}{2\pi i} \rt)^2\int\limits_{(c)}  \int\limits_{(0)}  \frac{(2\pi)^{1-u}h^{\hf-u-z}\sigma_{2z}(h)M_{0}(u,z/i,\delta)}{2\zeta^*(1-2z)\zeta^*(1+2z)(u-s)}  \ol{\LA V_{f,\ell}, E^*(*,1/2+z) \RA} \ dz du. \notag \\
 & \ \ \ \ + \sum_{r=0}^{\kt} \frac{1}{2\pi i} \int_{(0)}  \frac{(-1)^r(4\pi)^{\hf+z+r}\sigma_{-2z}(h)h^r\Gamma(-2z-r)\Gamma(\hf+z+r)}{r! \mcV  \zeta^*(1+2z)\zeta^*(1-2z)\Gamma(\hf+z)\Gamma(\hf-z)(\hf-z-r-s)} \ol{\LA V_{f,\ell}, E^*(*,\thf+z)) \RA} \ dz \notag \\
 & \ \ \ \ + \sum_{r=0}^\kt \frac{1}{2\pi i } \int_{(0)} \frac{\Or_{f,r,b}\lt(h^{r+\vep} (1+|\im z|)^{r-b}e^{-\pif |\im z|} \delta \rt)}{\zeta^*(1+2z)\zeta^*(1-2z)(\hf-z-r-s)}  \ol{\LA V_{f,\ell}, E^*(*,\thf+z)) \RA}   \ dz \notag
  \end{align}
where the residual terms were obtained by \eqref{resfin}. From an argument analogous to the one that gave us the meromorphic continuation in \eqref{zeepoles}, we see that the poles of the error term in $z$ are a subset of the poles of the main term, and their residues have the accompanying $\delta$ term. Rewriting the above equality we get the identity
\begin{align}
&\mc{I}_{f,\ell}^{\mbox{\tiny con}}(s;h,\delta) =  \frac{1}{2\pi i} \int_{(2)} \frac{\mc{I}^{\mbox{\tiny con}}_{f,\ell}(u;h,\delta)}{(u-s)} \ du  -  \frac{1}{2\pi i} \int_{(c)} \frac{(\mc{I}^{\mbox{\tiny con}}_{f,\ell}(u;h,\delta)-\Omega_{f,\ell}(u;h,\delta))}{(u-s)} \ du  + \mc{O}_{s,f,\ell}(\delta) \notag \\
 & \ \ \ \ - \sum_{r=0}^{\kt} \frac{1}{2\pi i} \int_{(0)}  \frac{(-1)^r(4\pi)^{\hf+z+r}\sigma_{-2z}(h)h^r\Gamma(-2z-r)\Gamma(\hf+z+r)}{r! \mcV  \zeta^*(1+2z)\zeta^*(1-2z)\Gamma(\hf+z)\Gamma(\hf-z)(\hf-z-r-s)} \ol{\LA V_{f,\ell}, E^*(*,\thf+z)) \RA} \ dz
\end{align}
when $\hf < \sigma < 2$. Following an argument that is nearly identical to the one that produced \eqref{thetag1} and \eqref{thetag2}, we get the meromorphic continuation to all $s$ such that $c < \sigma < 2$ by meromorphically continuing the sum of integrals, thus
\begin{align}
&\mc{I}_{f,\ell}^{\mbox{\tiny con}}(s;h,\delta) =  \frac{1}{2\pi i} \int_{(2)} \frac{\mc{I}^{\mbox{\tiny con}}_{f,\ell}(u;h,\delta)}{(u-s)} \ du  -  \frac{1}{2\pi i} \int_{(c)} \frac{(\mc{I}^{\mbox{\tiny con}}_{f,\ell}(u;h,\delta)-\Omega_{f,\ell}(u;h,\delta))}{(u-s)} \ du  \label{coni} \\
 & \ \ \ \ - \sum_{r=0}^{\kt} \frac{1}{2\pi i} \int_{C'_\sigma}  \frac{(-1)^r(4\pi)^{\hf+z+r}\sigma_{-2z}(h)h^r\Gamma(-2z-r)\Gamma(\hf+z+r)}{r! \mcV  \zeta^*(1+2z)\zeta^*(1-2z)\Gamma(\hf+z)\Gamma(\hf-z)(\hf-z-r-s)} \ol{\LA V_{f,\ell}, E^*(*,\thf+z)) \RA} \ dz  \notag \\
 & \ \ \ \ + \Omega_{f,\ell}(s;h) + \mc{O}_{s,f,r,b}(\delta)  \notag
\end{align}
where the poles of the error term are a subset of the poles of  $\Omega_{f,\ell}(s;h)$. Combining \eqref{coni} and \eqref{else} we get that for $c < \sigma < 2$, 
\begin{subequations}\label{aconint}
\bal
&\mc{I}_{f,\ell}(s;h,\delta) \\
&= \frac{1}{2\pi i} \int_{2-i\infty}^{2+i\infty} \frac{\mc{I}_{f,\ell}(u;h,\delta)}{(u-s)} \ du -  \frac{1}{2\pi i} \int_{c-i\infty}^{c+i\infty} \frac{(\mc{I}_{f,\ell}(u;h,\delta)-\Omega_{f,\ell}(s;h,\delta))}{(u-s)} \ du \label{aconinta}  \\
&-\sum_{r=0}^\kt \int_{C'_\sigma} \frac{(-1)^r(4\pi)^{\hf+z+r}\sigma_{-2z}(h)h^r\Gamma(-2z-r)\Gamma(\hf+z+r)\ol{\LA V_{f,\ell}, E^*(*,\thf+z)) \RA}}{r! \mcV  \zeta^*(1+2z)\zeta^*(1-2z)\Gamma(\hf+z)\Gamma(\hf-z)(\hf-z-r-s)}  \ dz \label{aconintb}  \\
&+\Omega_{f,\ell}(s;h)+ \mc{O}_{s,f,r,b}(\delta) -\sum_{j \neq 0}\sum_{r=0}^{k/2} \frac{ \wt{R}_{f,\ell}(1/2+it_j-r;h,\delta)}{(\hf+it_j-r-s)}.  \label{aconintc}
\end{align}
\end{subequations}
we will use this alternate form of $\mc{I}_{f,\ell}(s;h,\delta)$ to help take the limit as $\delta \to 0$. 

Combining \eqref{mer2} with \eqref{done} and applying Stirling's approximation gives us the first integral in \eqref{aconinta} converges absolutely and uniformly as $\delta \to 0$. We similarly get convergence of the second integral in \eqref{aconinta} via \eqref{aless2}. From \eqref{rrbound} we can get that  $\wt{R}_{f,\ell}(1/2+it_j-r;h,\delta) \to c_{r,j,f,\ell}(h)$ as $\delta \to 0$ and  \eqref{cbound}, \eqref{rrbound} and Stirling's approximation give us \eqref{rbound4} and that the sum over $j \neq 0$ still converges in the limit. All of this together shows that \eqref{aconint} converges to \eqref{anconint} as $\delta \to 0$, giving a meromorphic continuation of $\mc{I}_{f,\ell}(s;h)$ into the region $c < \sigma <2$, which agrees with the other two convergent formulas for $\mc{I}_{f,\ell}(s;h)$ in their respective regions. It is also clear from \eqref{anconint} and \eqref{spec3} that the poles of $\mc{I}_{f,\ell}(s;h)$ are as they are described in the statement of the proposition. We get the bound \eqref{later1} by inputting \eqref{aless3} and \eqref{easy2} into the first few integrals of\eqref{anconint}, applying \eqref{rbound4} to the sum over the cuspidal residues, and the sum over integrals is bounded using the method demonstrated in Proposition \ref{ellind}.

The poles of $D_{f,\ell}^+(s;h)$ can be classified via \eqref{mer3} and our knowledge about the poles of $\mc{I}_{f,\ell}(s;h)$. 
  \end{proof}

 \chapter{The Double Sum \lowercase{$\uppercase{Z}^+_{f,\ell}(s,w)$}} 
 \section{The Region of Spectral Convergence}
 Now that we have the meromorphic continuation of $D^+_{f,\ell}(s;h)$, we want to make use of it in the spectral expansion of $T^+$ as in \eqref{refy2}. As described in the introduction, we accomplish this by examining and meromorphically continuing
 \beq
 Z_{f,\ell}^+(s,w):=\sum_{m,h>0} \frac{a(m+h)\ol{\lambda_\ell(m)}}{m^{s+\frac{k-1}{2}}h^{w+\frac{k-1}{2}}} = \sum_{h=0}^\infty \frac{D_{f,\ell}^+(s,h)}{h^{w+\frac{k-1}{2}}} \label{anothersum}, \ \ \mbox{when} \ \ \re(s,w) > 1
 \eeq
 which, as was described in the introduction, will feed directly into \eqref{refy2}. 
 Letting $\sigma:= \re s$ and $\omega : = \re w$, we see that the above series easily converges absolutely when $\sigma,\omega>1$.  Throughout this chapter, it will be helpful to refer to the complex variable $w_2 =:s+w+\kt-1$ with $\omega_2:=\re w_2$. 
 
 To begin to meromorphically continue this function to the rest of the complex plane, we will need the following proposition.
 \begin{proposition}\label{sumswithin}
The function $Z_{f,\ell}^+(s,w)$ has a meromorphic continuation to the region given by 
\begin{align}
& \lt\{\begin{array}{ll}
 \omega_2 > \sigma+1+\tkt+\frac{\beta_1}{2}+\theta & \mbox{ if \ }  \sigma\geq \thf-\tkt-\frac{\beta_1}{2}\\
\omega_2> \tfrac{3}{2}+\theta & \mbox{ if \ } \sigma<\thf-\tkt-\frac{\beta_1}{2}
 \end{array}\rt\} \\
& \ \ \ \ \ \ \ \ \ \ \ \ \ \quad \quad \quad \quad \ \notag \bigcup
\lt\{ \begin{array}{ll}
\sigma <\thf-\tkt-\frac{\beta_1}{2} & \mbox{ if \ } \omega_2 >1+\vep \\
\sigma < \frac{\omega_2-k-\beta_1-\vep}{2} & \mbox{ if } \  -\vep<\omega_2<1+\vep \\
\sigma<\omega_2-\kt-\frac{\beta_1}{2}& \mbox{ if } \ \omega_2 \leq -\vep
\end{array}\rt\}
\end{align}
or, a fortiori, in the more simply defined sub-region given by $\omega_2>\sigma+1+\kt+\frac{\beta_1}{2}+\theta$. 

When $\sigma<\hf-\kt-\frac{\beta_1}{2}$ in this region, $Z_{f,\ell}^+(s,w)$ is given by 
\bal
&Z_{f,\ell}^+(s,w)= \frac{(4\pi)^{\kt}\Gamma(s)}{\ol{\rho_{\ell,k}(1)}}\sum_{j} \frac{\ol{\rho_j(-1)}\Gamma(s-\thf+it_j)\Gamma(s-\thf+it_j)\Gamma(1-s)\ol{\LA V_{f,\ell},u_j \RA }}{\Gamma(s+\tfrac{k-1}{2}+it_\ell)\Gamma(s+\tfrac{k-1}{2}-it_\ell)\Gamma(\hf+it_j)\Gamma(\hf-it_j)}L(w_2,\ol{u_j}) \notag \\
&+\frac{1}{2\pi i } \int_{C_{s,w}}  \frac{(4\pi)^{\kt}\Gamma(s)\Gamma(s-\hf-z)\Gamma(s-\hf+z)\Gamma(1-s)\ol{\LA V_{f,\ell}, E^*(*,\hf+z) \RA }\zeta(w_2,z)}{2\ol{\rho_{\ell,k}(1)}\zeta^*(1-2z)\zeta^*(1+2z)\Gamma(s+\tfrac{k-1}{2}+it_\ell)\Gamma(s+\tfrac{k-1}{2}-it_\ell)\Gamma(\hf+z)\Gamma(\hf-z)}   \ dz \notag  \\ 
& + \wt{\Omega}_{f,\ell}(s,w) +\eta_{f,\ell}(s,w) \notag \\
& - \lt. \frac{\epsilon_{\ell,k} \Gamma(s)\Gamma(1-s)}{2\pi i} \int_{(-\vep)}\frac{\Gamma(-z)\Gamma\lt(w+\tfrac{k-1}{2}+z\rt)L(w+z,f)L\lt(s-z+\tfrac{k-1}{2},\ol{\mu_\ell}\rt)}{\Gamma\lt(w+\frac{k-1}{2}\rt)\Gamma(\hf-\kt+it_\ell)\Gamma(\hf+\kt-it_\ell)}\ dz \rt.  \label{spec100} 
\end{align}
where $\epsilon_{\ell,k}=\pm 1$ as in \eqref{thatthing} and we choose $\vep>0$ to be small such that $\omega + \kt -\thf -\vep \notin \mt{Z}_{\leq 0}$, $C_{s,w}$ will be described following \eqref{nicetry}, and
\beq
\wt{\Omega}_{f,\ell}(s,w):=G_\ell(s) \sum_{h=1}^\infty \frac{\Omega_{f,\ell}(s;h)}{h^{w+\frac{k-1}{2}}}.
\eeq
where $\Omega_{f,\ell}(s;)$ is as in \eqref{omega14} and
\bal
&\eta_{f,\ell}(s,w):= \\
& \ \ \ \frac{\pi^{w_2-\hf}(4\pi)^\kt \Gamma(s) \Gamma(1-s)\Gamma(s+\hf-w_2)\Gamma(s-\tfrac{3}{2}+w_2)\ol{\LA V_{f,\ell}, E^*(*,w_2-\thf) \RA }\psi_{(-\infty,1]}(\omega_2)}{\ol{\rho_{\ell,k}(1)}\zeta^*(3-2w_2)\Gamma(s+\frac{k-1}{2}+it_\ell)\Gamma(s+\frac{k-1}{2}-it_\ell)\Gamma(\tfrac{3}{2}-w_2)\Gamma(w_2-\hf)^2} \notag 
\end{align}
where $\psi_{(-\infty,1]}$ is the characteristic function on $(-\infty,1]$. 
This continuation has polar lines in $s$ wherever $D_{f,\ell}^+(s;h)$ has poles. In particular, when $s=\hf+it_j-r$ for $r \in \mt{Z}_{\geq 0}$, unless $t_j = \pm t_\ell$ and $r\geq \kt$, we have the residues
\beq
\resie{s=\hf + it_j -r} Z_{f,\ell}^+(s,w)=G_\ell(\thf+ it_j-r)c_{r,j,f,\ell}(1)L(w_2,\ol{u_j}), \label{polarlines}
\eeq
where $c_{r,j,f,\ell}(1)$ is as in \eqref{round2}. Furthermore, we have polar lines of the form $s+w_2-\frac{3}{2} \in \mt{Z}_{\leq 0}$ with residues
\begin{align}
&\resie{w_2=\frac{3}{2}-s-r} Z_{f,\ell}^+(s,w)= \label{otherstupidpole}\frac{(-1)^r \pi^{1-s-r}(4\pi)^\kt \Gamma(s) \Gamma(1-s)\Gamma(2s+r-1)\ol{\LA V_{f,\ell}, E^*(*,1-s-r) \RA }}{r! \ol{\rho_{\ell,k}(1)}\zeta^*(2s+2r)\Gamma(s+\frac{k-1}{2}+it_\ell)\Gamma(s+\frac{k-1}{2}-it_\ell)\Gamma(s+r)\Gamma(1-s-r)^2}.
\end{align}

Let $A_2 > 1+|\sigma|+|\omega|+k$ and let $\pmb{P}(A_2)$ be an unspecified piecewise affine function in $A_2$.  In the region of convergence where   $\sigma < \hf-\tkt-\frac{\beta_1}{2}$ we have the bound
\beq
Z_{f,\ell}^+(s,w) \ll_{A_2} \lt[(1+|s|)(1+|w|)(1+|t_\ell|)\rt]^{\pmb{P}(A_2)}e^{\pi |t_\ell|}
 \label{tell1}
\eeq
and for $\sigma > 1+\vep$ we have the bound
\beq
Z_{f,\ell}^+(s,w) \ll_\vep O(1). \label{tell3}
\eeq 
 \end{proposition}
 \begin{proof} We begin by inserting the continuation of $D_{f,\ell}^+(s,h)$  from Proposition \ref{prop5} into \eqref{anothersum}. We observe that for fixed $s$, not at any of the poles of $G_\ell(s)\Omega_{f,\ell}(s;h)$, the function $\wt{\Omega}_{f,\ell}(s,w)$ is absolutely convergent for sufficiently large $\omega$ and has a meromorphic continuation to all $w \in \mt{C}$ as we can interchange the order of summation for $w_2$ sufficiently big, as
 \beq
 \sum_{h=1}^\infty \frac{h^{n}\sigma_{1-2s-2n}}{h^{w+\frac{k-1}{2}-n}} = \zeta(w_2,\thf-s-n)
 \eeq
 where
 \beq
  \zeta(s,w) := \zeta(s+w)\zeta(s-w).
 \eeq
Letting
\beq
D^*_{f,\ell}(s;h) := D_{f,\ell}^+(s;h)-\Omega_{f,\ell}(s;h)
\eeq
we now direct our attention to the sum
\beq
\sum_{h=0}^\infty \frac{D^*_\ell(s,h)}{h^{w+\frac{k-1}{2}}} = Z^+_{f,\ell}(s,w) - \wt{\Omega}_{f,\ell}(s,w).
\eeq
Using the bounds given in \eqref{aless3}, \eqref{easy2}, and \eqref{later1}, we see the above sum is absolutely convergent in $h$ when 
 \beq
 \lt\{\begin{array}{ll}
 \omega > 2+\frac{\beta_1}{2}+\theta & \mbox{ if \ }  \sigma\geq \thf-\tkt-\frac{\beta_1}{2} \\
\sigma+\omega+\tfrac{k-1}{2}-\theta-1>1 & \mbox{ if \ }  \sigma<\thf-\tkt-\frac{\beta_1}{2}
 \end{array}\rt\}
 \eeq
 and $\sigma$ is some $\vep>0$ distance away from the poles of $D_{f,\ell}^+(s;h)$ not due to $\wt{\Omega}_{f,\ell}(s,w)$. We see that this gives us convergence in the equivalent region of
  \beq
 \lt\{\begin{array}{ll}
 \omega_2 > \sigma+1+\tkt+\frac{\beta_1}{2}+\theta & \mbox{ if \ }  \sigma\geq \thf-\tkt -\frac{\beta_1}{2}\\
\omega_2> \tfrac{3}{2}+\theta & \mbox{ if \ }  \sigma<\thf-\tkt-\frac{\beta_1}{2}
 \end{array}\rt\}.
 \eeq
Given the meromorphic continuation of $\wt{\Omega}_{f,\ell}(s,w)$ for fixed $s$, we also see that this absolute convergence gives a meromorphic continuation of $Z^+_{f,\ell}(s,w)$ to all $s,w$ in the above region with polar lines in $s$ due to those of $D_{f,\ell}^+(s;h)$. It is easy to show that \eqref{polarlines} holds from the representations of $D_{f,\ell}^+(s;h)$ given in Proposition \ref{prop5}, and furthermore that all other poles that are independent of $w_2$ in this region are due to the $\Gamma(s)$ in $G_\ell(s)$, the  and the poles of $\Omega_{f,\ell}(s;h)$. 

When $\sigma < \hf-\kt-\frac{\beta_1}{2}$ we can insert the spectral expansion for $G_\ell(s)\mc{I}_{f,\ell}(s,h)$ from Proposition \ref{prop5} into \eqref{anothersum} to get
\bal
&Z_{f,\ell}^+(s,w)= \frac{(4\pi)^{\kt}\Gamma(s)}{\ol{\rho_{\ell,k}(1)}}\sum_{j} \frac{\ol{\rho_j(-1)}\Gamma(s-\thf+it_j)\Gamma(s-\thf+it_j)\Gamma(1-s)\ol{\LA V_{f,\ell},u_j \RA }}{\Gamma(s+\tfrac{k-1}{2}+it_\ell)\Gamma(s+\tfrac{k-1}{2}-it_\ell)\Gamma(\hf+it_j)\Gamma(\hf-it_j)}\sum_{n=1}^\infty \frac{\ol{\lambda_j(h)}}{h^{w_2}}\notag \\
&+\frac{1}{2\pi i } \int_{C_\sigma}  \frac{(4\pi)^{\kt}\Gamma(s)\Gamma(s-\hf-z)\Gamma(s-\hf+z)\Gamma(1-s)\ol{\LA V_{f,\ell}, E^*(*,\hf+z) \RA }\lt(\sum_{h=1}^\infty \frac{\sigma_{2z}(h)}{h^{w_2+z}}\rt) }{2\ol{\rho_{\ell,k}(1)}\zeta^*(1-2z)\zeta^*(1+2z)\Gamma(s+\tfrac{k-1}{2}+it_\ell)\Gamma(s+\tfrac{k-1}{2}-it_\ell)\Gamma(\hf+z)\Gamma(\hf-z)} dz \notag  \\ 
&  -\epsilon_{\ell,k} \frac{\Gamma(s)\Gamma(1-s)}{\Gamma(\hf-\kt+it_\ell)\Gamma(\hf+\kt-it_\ell)} \sum_{h=1}^\infty \sum_{h > m \geq 1} \frac{a(h-m)\ol{\lambda_\ell(m)}}{m^{s+\frac{k-1}{2}}h^{w+\frac{k-1}{2}}}\notag \\
& + \wt{\Omega}_{f,\ell}(s,w)  \label{spec10}
\end{align}
 where, again, interchanging the order of summation is allowed for $\omega_2$ sufficiently big. We will show that \eqref{spec100} will follow from \eqref{spec10} when the appropriate substitutions are made.
Since 
\beq
\sum_{h=1}^\infty  \frac{\ol{\lambda_j(h)}}{h^{w_2}} = L(w_2,\ol{u_j})
\eeq
when $\omega_2>1$, we can substitute it into the first line of \eqref{spec10}. By the functional equation of $L(w_2,\ol{u_j})$ and Stirling's approximation we get the convexity bound
\beq
L(w_2,\ol{u_j}) \ll_\vep \lt\{ 
\begin{array}{ll}
1 & \mbox{if } \  \omega_2 > 1+\vep \\
(1+|w_2|)^{1-\omega_2+\vep}(1+|t_j|)^{1-\omega_2+\vep} & \mbox{if } \ -\vep <\omega_2 < 1+\vep \\
(1+|w_2|)^{1-2\omega_2}(1+|t_j|)^{1-2\omega_2} & \mbox{if } \  \omega_2 < -\vep
\end{array}
\rt. \label{usefulcompa}
\eeq
for $\vep>0$. Again making use of Stirling's formula, we see that the sum over eigenvalues in \eqref{spec10} converges when
\beq
\lt\{ \begin{array}{ll}
\sigma < \thf-\tkt-\frac{\beta_1}{2} & \mbox{ if \ } \omega_2 >1+\vep \\
\sigma < \frac{\omega_2-k-\beta_1-\vep}{2} & \mbox{ if } \  -\vep<\omega_2<1+\vep \\
\sigma<\omega_2-\kt-\frac{\beta_1}{2}& \mbox{ if } \ \omega_2 \leq -\vep
\end{array}\rt\} \label{regionconvexa2}
\eeq
within $\vep>0$ of the poles at $s-\hf-it_j \in \mt{Z}_{\leq 0}$. 

Similarly, since
\beq
\sum_{n=1}^\infty \frac{\sigma_{2it}(n)}{n^{w_2+it}}=\zeta(w_2,it)
\eeq
we can substitute the above equality into the second line of \eqref{spec10}. From the functional equation of $\zeta(s)$ and bounds due to convexity we have that
\beq
\zeta(w_2\pm it) \ll_\vep \lt\{ 
\begin{array}{ll}
1 & \mbox{if } \  \omega_2 > 1+\vep \\
(1+|w_2|)^{\frac{1}{2}-\frac{\omega_2}{2}+\vep}(1+|t|)^{\frac{1}{2}-\frac{\omega_2}{2}+\vep} & \mbox{if } \ -\vep <\omega_2 < 1+\vep \\
(1+|w_2|)^{\frac{1}{2}-\omega_2} (1+|t|)^{\frac{1}{2}-\omega_2} & \mbox{if } \  \omega_2 < -\vep. 
\end{array} \rt. \label{neverdone2}
\eeq
Making use of Stirling's formula, we see that the integral (and thus the meromorphic continuation of the integral) on the second line of \eqref{spec10} converges absolutely when
\beq
\lt\{ \begin{array}{ll}
\sigma<1-\tkt & \mbox{ if \ } \omega_2 >1+\vep \\
\sigma<  \frac{\omega_2+1-k-\vep}{2} & \mbox{ if } \  -\vep<\omega_2<1+\vep \\
\sigma<\omega_2+\thf-\kt & \mbox{ if } \ \omega_2<-\vep
\end{array}\rt\} \label{regionconvexa3}
\eeq
which encompasses the region in \eqref{regionconvexa2}. We see that, other than those poles due to $\Gamma(s)$, this integral component contributes no other poles in $s$ as they are included in $\wt{\Omega}_{f,\ell}(s,w)$.  We do, however, have to understand how the function given by the integral changes as we move past the line $\omega_2 =1$. By taking $\sigma$ to be sufficiently negative and $\omega_2 >1$,we consider the relevant integral
\begin{subequations} \label{superint42}
\bal
& Z_{f,\ell,\mbox{\tiny con}}^+(s,w):= \label{nicetry}  \\
& \frac{1}{2\pi i } \int_{C_\sigma}  \frac{(4\pi)^{\kt}\Gamma(s)\Gamma(s-\hf-z)\Gamma(s-\hf+z)\Gamma(1-s)\ol{\LA V_{f,\ell}, E^*(*,\hf+z) \RA }\zeta(w_2,z)}{2\ol{\rho_{\ell,k}(1)}\zeta^*(1-2z)\zeta^*(1+2z)\Gamma(s+\tfrac{k-1}{2}+it_\ell)\Gamma(s+\tfrac{k-1}{2}-it_\ell)\Gamma(\hf+z)\Gamma(\hf-z)}   \ dz. \notag
\end{align}
Let $C_{s,w}$ denote a vertically-alligned contour similar to $C_\sigma$ in that it occupies the narrow region around $\re z =0$ where $\zeta^*(1\pm 2z)$ is proven to be zero-free, which we'll call $B_0$. Usually $C_{s,w}$ will be just the curve $C_\sigma$ but whenever $\omega_2$ gets sufficiently close to $1$, let $C_{s,w}$ be the deformation of $C_\sigma$ that stays in the region $B_0$ but passes over both poles at $z=\pm(1-w_2)$ without passing over the poles at $z=\pm(\hf-s-r)$. Thus we have that, for fixed $s$, the integral $ Z_{f,\ell,\mbox{\tiny con}}^+(s,w)$ continues to   
\bal
& Z_{f,\ell,\mbox{\tiny con}}^+(s,w)= \\
& \frac{1}{2\pi i } \int_{C_{s,w}}  \frac{(4\pi)^{\kt}\Gamma(s)\Gamma(s-\hf-z)\Gamma(s-\hf+z)\Gamma(1-s)\ol{\LA V_{f,\ell}, E^*(*,\hf+z) \RA }\zeta(w_2,z)}{2\ol{\rho_{\ell,k}(1)}\zeta^*(1-2z)\zeta^*(1+2z)\Gamma(s+\tfrac{k-1}{2}+it_\ell)\Gamma(s+\tfrac{k-1}{2}-it_\ell)\Gamma(\hf+z)\Gamma(\hf-z)}   \ dz \notag \\
  & \ \ \ \ + \frac{(4\pi)^\kt \Gamma(s) \Gamma(1-s)\Gamma(s+\hf-w_2)\Gamma(s-\tfrac{3}{2}+w_2)\ol{\LA V_{f,\ell}, E^*(*,w_2-\thf) \RA }\zeta(2w_2-1)}{\ol{\rho_{\ell,k}(1)}\zeta^*(2w_2-1)\zeta^*(3-2w_2)\Gamma(s+\frac{k-1}{2}+it_\ell)\Gamma(s+\frac{k-1}{2}-it_\ell)\Gamma(\tfrac{3}{2}-w_2)\Gamma(w_2-\hf)}. \notag
\end{align}
\end{subequations}
when $\omega_2 \leq 1$. In the case where $w_2$ becomes close to $1$, we might also be concerned about the behavior of the function as $w_2 \to 1$ as $C_{s,w}$ is not defined. However,  we see that we can choose to shift over just the pole at $z = 1-w_2$, picking up half the above residue, and then taking the limit as $w_2 \to 1$, which we see is not a polar line. 
So the integral remains well-defined for $\omega_2 \leq 1$, and in this region and contributes no additional poles in $s$ or $w_2$. The residue appears to contribute polar lines due to $\Gamma(s+\thf-w_2)$ and $\Gamma(s-\frac{3}{2}+w_2)$. To better understand these potential polar lines we need to consider $\wt{\Omega}_{f,\ell}(s,w)$ as, by construction,  $Z_{f,\ell,\mbox{\tiny con}}^+(s,w) + \wt{\Omega}_{f,\ell}(s,w)$ is a meromorphic function in $s$ and $w$ in the region specified by \eqref{regionconvexa2}. We observe that in this region 
 \bal
& \wt{\Omega}_{f,\ell}(s,w) := \label{omega141}\\
&\sum_{n=0}^{\lf \hf -\sigma \rf } \frac{(4\pi)^{\kt}\Gamma(s)(-1)^n \zeta(w_2,\hf-s-n)\Gamma(2s+n-1)\Gamma(1-s)\ol{\LA V_{f,\ell}, E^*(*,s+n)) \RA}}{n!\ol{\rho_{\ell,k}(1)} \Gamma(s+\frac{k-1}{2}+it_\ell)\Gamma(s+\frac{k-1}{2}-it_\ell) \zeta^*(2s+2n)\zeta^*(2-2s-2n)\Gamma(s+n)\Gamma(1-s-n)},   \notag
\end{align}
which similarly appears to have polar lines due to $ \zeta(w_2,\hf-s-n)$, indeed these coincide with the lines above which appear at $s+ \hf -w_2 \in \mt{Z}_{\leq 0}$ and $s-\frac{3}{2}+w_2 \in \mt{Z}_{\leq 0}$. However, we note that $\zeta(w_2+\hf-s-n)$ and $\Gamma(s+\hf-w_2)$ contribute exact opposite residues which cancel out; this is to be expected since lines  $s+ \hf -w_2 \in \mt{Z}_{\leq 0}$ extend into the region where $\sigma$ and $\omega_2$ are very large and we see this does not occur in the representation of $Z^+_{f,\ell}(s,w)$ given in \eqref{anothersum}. Alternately, we see that $\zeta(w_2+s+n-\hf)$ and $\Gamma(s-\frac{3}{2}+w_2)$ contribute identical residues but in disjoint regions, that is $\wt{\Omega}_{f,\ell}(s,w)$ contributes the polar line when $\omega_2 >1$ and $Z_{f,\ell,\mbox{\tiny con}}^+(s,w)$ contributes the polar line when $\omega_2 < 1$. When $\omega_2 =1$, both parts contribute half of the residual term. This gives us \eqref{otherstupidpole}.

For the purposes of proving \eqref{tell1}, we need to bound $\wt{\Omega}_{f,\ell}(s,w)$ and the residual term in \eqref{superint42}. A pretty straightforward application of Stirling's Approximation and \eqref{youngtwoone} gives us
\beq
\wt{\Omega}_{f,\ell}(s,w) \ll_{A_2} [(1+|s|)(1+|w|)(1+|t_\ell|)]^{\pmb{P}(A_2)} \label{timetoeat}
\eeq
and 
\bal
& \frac{(4\pi)^\kt \Gamma(s) \Gamma(1-s)\Gamma(s+\hf-w_2)\Gamma(s-\tfrac{3}{2}+w_2)\ol{\LA V_{f,\ell}, E^*(*,w_2-\thf) \RA }\zeta(2w_2-1)}{\ol{\rho_{\ell,k}(1)}\zeta^*(2w_2-1)\zeta^*(3-2w_2)\Gamma(s+\frac{k-1}{2}+it_\ell)\Gamma(s+\frac{k-1}{2}-it_\ell)\Gamma(\tfrac{3}{2}-w_2)\Gamma(w_2-\hf)} \\
& \ll_{A_2} [(1+|s|)(1+|w|)(1+|t_\ell|)]^{\pmb{P}(A_2)} \notag \\
& \ \ \ \ \ \ \ \ \ \ \ \ \ \ \ \ \times e^{-\pif (||\im w_2|-|t_\ell||-||\im s|-|t_\ell||-3|\im w_2|+|\im s| + |\im s+\im w_2| +|\im s - \im w_2| )},  \notag
\end{align}
where $A_2$ and $\pmb{P}(A_2)$ are as in the statement of the proposition. Since 
\beq
||\im w_2|-|t_\ell||-|\im w_2| \geq -| t_\ell|,
\eeq
\beq
|\im s + \im w_2|+|\im s - \im w_2| - 2|\im w_2 | \geq 0
\eeq
and
\beq
|\im s| - ||\im s|-|t_\ell || \geq -|t_\ell|
\eeq
we get
\bal
& \frac{(4\pi)^\kt \Gamma(s) \Gamma(1-s)\Gamma(s+\hf-w_2)\Gamma(s-\tfrac{3}{2}+w_2)\ol{\LA V_{f,\ell}, E^*(*,w_2-\thf) \RA }\zeta(2w_2-1)}{\ol{\rho_{\ell,k}(1)}\zeta^*(2w_2-1)\zeta^*(3-2w_2)\Gamma(s+\frac{k-1}{2}+it_\ell)\Gamma(s+\frac{k-1}{2}-it_\ell)\Gamma(\tfrac{3}{2}-w_2)\Gamma(w_2-\hf)} \notag \\
& \ \ \ \ \ll_{A_2} [(1+|s|)(1+|w|)(1+|t_\ell|)]^{\pmb{P}(A_2)}   e^{\pi |t_\ell|}.  \label{wellthatsweird}
\end{align}

Now we look at the last line of \eqref{spec10}, from which we have the shifted sum
\beq
 \sum_{h \geq 1} \sum_{h>m\geq 1} \frac{a(h-m)\ol{\lambda_\ell(m)}}{m^{s+\frac{k-1}{2}}h^{w+\frac{k-1}{2}}}. 
\eeq
The above series is clearly absolutely convergent for any $s$ provided $\omega$ is sufficiently positive. Using the substitution $h-m \to h$ we can rewrite the above sum as
\beq
 \sum_{h,m=1}^\infty \frac{a(h)\ol{\lambda_\ell(m)}}{m^{s+\frac{k-1}{2}}(m+h)^{w+\frac{k-1}{2}}} = \sum_{h,m=1}^\infty \frac{a(h)\ol{\lambda_\ell(m)}}{m^{s+\frac{k-1}{2}}h^{w+\frac{k-1}{2}}}\lt(1+\frac{m}{h}\rt)^{-(w+\frac{k-1}{2})}. \label{huh1a}
\eeq 
which is still a convergent series. From 6.422(3) in \cite{GR} we have that
\beq
\frac{1}{2 \pi i \Gamma(\beta)}\int_{(-\gamma)} \Gamma(-z)\Gamma(\beta+z)t^z \ dz = (1+t)^{-\beta} \label{magic2}
\eeq 
when $0>-\gamma>-\re(\beta_1)$ and $| \arg t|<\pi$. Using \eqref{magic2}, letting $\sigma$ be sufficiently big, and the fact that $|\arg \frac{h}{m}|=0$, we can rewrite \eqref{huh1a} as
\beq
\frac{1}{2\pi i \Gamma\lt(w+\frac{k-1}{2}\rt)}\sum_{h,m=1}^\infty \frac{a(h)\ol{\lambda_\ell(m)}}{m^{s+\frac{k-1}{2}}h^{w+\frac{k-1}{2}}} \int_{(-\gamma)}\Gamma(-z)\Gamma\lt(w+\tfrac{k-1}{2}+z\rt)\lt(\frac{m}{h}\rt)^z \ dz. \label{huh2a}
\eeq
Since
\beq
 \int_{(-\gamma)}\Gamma(-z)\Gamma\lt(w+\tfrac{k-1}{2}+z\rt)\lt(\frac{m}{h}\rt)^z \ dz \ll_{w,k} \lt(\frac{h}{m}\rt)^\gamma
\eeq
we can interchange the order of the sums and integrals in \eqref{huh2a}, provided $\omega$ and $\gamma$ are both sufficiently large so that $s+\gamma+\frac{k-1}{2} >1$ and $w+\frac{k-1}{2}-\gamma>1$. We get
\beq
\frac{1}{2\pi i \Gamma\lt(w+\frac{k-1}{2}\rt)}  \int_{(-\gamma)}\Gamma(-z)\Gamma\lt(w+\tfrac{k-1}{2}+z\rt)\sum_{h,m=1}^\infty \frac{a(h)\ol{\lambda_\ell(m)}}{m^{s-z+\frac{k-1}{2}}h^{w+z+\frac{k-1}{2}}}  \ dz \label{huh3a}
\eeq
and observe that the nested sum
\beq
\sum_{h,m=1}^\infty \frac{a(h)\ol{\lambda_\ell(m)}}{m^{s-z+\frac{k-1}{2}}h^{w+z+\frac{k-1}{2}}} =\sum_{m=1}^\infty \frac{\ol{\lambda_\ell(m)}}{m^{s-z+\frac{k-1}{2}}}\sum_{h=1}^\infty \frac{a(h)}{h^{w+z+\frac{k-1}{2}}}=L(w+z,f)L\lt(s-z+\tfrac{k-1}{2},\ol{\mu_\ell}\rt)
\eeq
is the product of $L$-functions associated to $f$ and $\ol{\mu_\ell}$. This product has an analytic continuation to all $s$ and $w$ for any fixed $z$. So we can shift the line of integration from $\re z = -\gamma$ to $\re z = -\vep$ to get
\bal
 \sum_{h \geq 1} &\sum_{h>m\geq 1}  \frac{a(h-m)\ol{\lambda_\ell(m)}}{m^{s+\frac{k-1}{2}}h^{w+\frac{k-1}{2}}}  \notag
\\ &= \frac{1}{2\pi i \Gamma\lt(w+\frac{k-1}{2}\rt)}  \int_{(-\vep)}\Gamma(-z)\Gamma\lt(w+\tfrac{k-1}{2}+z\rt)L(w+z,f)L\lt(s-z+\tfrac{k-1}{2},\ol{\mu_\ell}\rt) \ dz, \label{supermagica}
\end{align}
since the $L$-functions are well-defined everywhere with at most polynomial growth in $\im z$ and the other two gamma factors in the integrand give more than sufficient decay for convergence. The integrand has no poles in $s$ and $w$ since $\Gamma(s+\frac{k-1}{2})L(s,f)$ and $L\lt(s,\ol{\mu_\ell}\rt)$ are analytic. Thus we see that the third line of \eqref{spec10} is meromorphic for all $(s,w) \in \mt{C}^2$ with it's only poles resulting from $\Gamma(s)$. 

We see that we now must also be concerned with the growth of \eqref{supermagica} in terms of $t_\ell$. To this end we use Stirling's formula to get that
\bal
& \int_{(-\vep)}\Gamma(-z)\Gamma\lt(w+\tfrac{k-1}{2}+z\rt)L(w+z,f)L\lt(s-z+\tfrac{k-1}{2},\ol{\mu_\ell}\rt)\ dz \notag \\ 
& \ll_{\sigma,\omega,k, \vep} \int_{-\infty}^\infty \frac{(1+|t|)^{\vep-\hf}(1+|\im w+t|)^{\omega+\frac{k-1}{2}-\vep-\hf}}{e^{\frac{\pi}{2}(|\im w+t|+|t|)}}|L(w-\vep+it,f)L(s+\vep-it+\tfrac{k-1}{2},\ol{\mu_\ell})| \ dt.
\end{align}
As the $L$-functions contribute, at most, polynomial growth in $t$ and $\im w$, and $|\im w+t|+|t|=\max(|\im w|, |\im w +2t|)$, we see that we get the largest contribution from the integral when $\im w$ and $t$ have opposite sign and $|t| \leq |\im w|$. Thus, by continued application of Stirling's approximation,
\bal
& \frac{1}{2\pi i \Gamma\lt(w+\frac{k-1}{2}\rt)}  \int_{(-\vep)}\Gamma(-z)\Gamma\lt(w+\tfrac{k-1}{2}+z\rt)L(w-z,f)L\lt(s+z+\tfrac{k-1}{2},\ol{\mu_\ell}\rt) \ dz \notag \\
& \ \ \ \ \ \ \  \ll_{A_2} (1+|\im s|)^{\pmb{P}({A_2})}(1+|\im w|)^{\pmb{P}({A_2})}(1+|t_\ell|)^{\pmb{P}({A_2})} \label{hahahaha}
\end{align}
where $A_2$ and $\pmb{P}(A_2)$  are as described in the proposition. 

The bound \eqref{tell3} follows pretty trivially from \eqref{anothersum} and the fact that the Ramanujan conjecture is known on average.  To get the bound in \eqref{tell1}, we merely reproduce the argument that gave us \eqref{aless3}, but also accounting for $G_\ell(s)$, \eqref{usefulcompa}, \eqref{neverdone2}, \eqref{hahahaha}, \eqref{timetoeat} and \eqref{wellthatsweird}.
\end{proof}
\section{The Expanded Region of Dirichlet Series Convergence}
Now we will show that we can come up with an analytic continuation of $Z_{f,\ell}^+(s,w)$ into another overlapping region by enlarging the region where $Z_{f,\ell}^+(s,w)$ converges as a multiple Dirichlet series. 
\begin{proposition} \label{sumswithout} Fix $A \gg 1$. Let $w_2 = s+w+\frac{k}{2}-1$. The function $Z_{f,\ell}^+(s,w)$ has an analytic continuation to the region where $\omega_2:=\re w_2 >3/2$ and  $\sigma:=\re s >1$. Furthermore, in the sub-region where $1 \geq \omega:= \re w$ and $\omega_2 > \frac{3}{2}+\vep$ we can let  $K \in \mt{R}$ such that $1+\vep> \re w +K > 1$ and $Z_{f,\ell}^+(s,w)$ can be described as the analytic continuation of the series
\bal
Z_{f,\ell}^+(s,w)=&\sum_{m\geq h} \frac{a(m+h)\ol{\lambda_\ell(m)}}{m^{s+\frac{k-1}{2}}h^{w+\frac{k-1}{2}}}
+\sum_{0 \leq j \leq K} \binom{w+\tfrac{k-1}{2}+j-1}{j} \sum_{m< h} \frac{A(m+h)\ol{\lambda_\ell(m)}}{m^{s+\frac{k-1}{2}-j}(m+h)^{w+j}}  \notag \\
& +\sum_{j>K} \binom{w+\tfrac{k-1}{2}+j-1}{j} \sum_{m< h} \frac{A(m+h)\ol{\lambda_\ell(m)}}{m^{s+\frac{k-1}{2}-j}(m+h)^{w+j}}.
\end{align}
Everywhere in the larger region we have the bound
\beq \label{restatement}
Z_{f,\ell}^+(s,w) \ll_{A_3} (1+|w|)^{\pmb{P}(A_3)}
\eeq
where $A_3 > 1+ |\omega|$ and $\pmb{P}(A_3)$ is an unspecified affine function in $A_3$ with no notable growth in the $s$ and $t_\ell$ aspects.

\end{proposition} 
\begin{proof}
Let
\beq
Z_{f,\ell}^+(s,w)=S_1(s,w)+S_2(s,w)
\eeq
where 
\beq
S_1(s,w):=\sum_{m\geq h} \frac{a(m+h)\ol{\lambda_\ell(m)}}{m^{s+\frac{k-1}{2}}h^{w+\frac{k-1}{2}}}
\eeq
and
\beq
S_2(s,w):=\sum_{ m<h} \frac{a(m+h)\ol{\lambda_\ell(m)}}{m^{s+\frac{k-1}{2}}h^{w+\frac{k-1}{2}}}.
\eeq
Both are absolutely convergent in the region where $\sigma, \omega>1$. Furthermore, since
\beq
S_1(s,w)=\sum_{m=1}^\infty \frac{\ol{\lambda_\ell(m)}}{m^s} \sum_{h=1}^m \frac{A(m+h)(1+\tfrac{h}{m})^{\frac{k-1}{2}}}{h^{w+\frac{k-1}{2}}}
\eeq
and
\beq
 \sum_{h=1}^m \frac{A(m+h)(1+\tfrac{h}{m})^{\frac{k-1}{2}}}{h^{w+\frac{k-1}{2}}} \ll_k \lt\{
 \begin{array}{ll}
\ln(m) & \mbox{ if } \ \re w +\tfrac{k-1}{2} \geq 1 \\
 m^{1-(w+\frac{k-1}{2})} & \mbox{ if } \re w +\tfrac{k-1}{2} < 1,
  \end{array}
  \rt.  \label{trivialbd}
\eeq
in this region, we see that $S_1(s,w)$ is absolutely convergent if $\sigma >1$ and $\omega_2 > \frac{3}{2}$ and has no notable growth in the $t_\ell$, $s$ or $w$  aspects.

Now we can rewrite $S_2(s,w)$ as 
\beq
S_2(s,w) = \sum_{m< h} \frac{A(m+h)\ol{\lambda_\ell(m)}}{m^{s+\frac{k-1}{2}}(m+h)^{w}}\lt(\frac{m+h}{h}\rt)^{w+\frac{k-1}{2}}.
\eeq
Note that
\beq
\lt(\frac{m+h}{h}\rt)^{w+\frac{k-1}{2}} = \lt(1-\frac{m}{m+h}\rt)^{-(w+\frac{k-1}{2})} = \sum_{j\geq0} \binom{w+\tfrac{k-1}{2}+j-1}{j} \lt( \frac{m}{m+h}\rt)^{j}
\eeq
is absolutely convergent. For some fixed $K \in \mt{R}$ we see that
\beq
S_2(s,w)=S_3^K(s,w)+S_4^K(s,w)
\eeq
where
\beq
S_3^K(s,w):=\sum_{0 \leq j \leq K} \binom{w+\tfrac{k-1}{2}+j-1}{j} \sum_{m< h} \frac{A(m+h)\ol{\lambda_\ell(m)}}{m^{s+\frac{k-1}{2}-j}(m+h)^{w+j}}
\eeq
and
\beq
S_4^K(s,w):=\sum_{j>K} \binom{w+\tfrac{k-1}{2}+j-1}{j} \sum_{m< h} \frac{A(m+h)\ol{\lambda_\ell(m)}}{m^{s+\frac{k-1}{2}-j}(m+h)^{w+j}},
\eeq
where everything remains absolutely convergent if $\sigma, \omega >1$.
Indeed, we see that
\bal
S_4^K(s,w) &= \sum_{j>K} \binom{w+\tfrac{k-1}{2}+j-1}{j} \sum_{m< h} \frac{A(m+h)\ol{\lambda_\ell(m)}}{m^{s+\frac{k-1}{2}-K}(m+h)^{w+K}}\lt(\frac{m}{m+h}\rt)^{j-K} \notag \\
& \ll \lt( \sum_{j>K} \binom{w+\tfrac{k-1}{2}+j-1}{j}\lt(\frac{1}{2}\rt)^{j-K} \rt)\lt(\sum_{m< h} \frac{A(m+h)\ol{\lambda_\ell(m)}}{m^{s+\frac{k-1}{2}-K}(m+h)^{w+K}}\rt),
\end{align}
where the sum over $j$ is absolutely convergent and
\beq
\sum_{m< h} \frac{A(m+h)\ol{\lambda_\ell(m)}}{m^{s+\frac{k-1}{2}-K}(m+h)^{w+K}}=\sum_{h=1}^\infty \frac{1}{h^{w+K}} \sum_{m=1}^{h-1} \frac{A(m+h)\ol{\lambda_\ell(m)}}{m^{s+\frac{k-1}{2}-K}(1+\frac{m}{h})^{w+K}}.
\eeq
Using a result comparable to that given in $\eqref{trivialbd}$ we see that $S_4^K(s,w)$ is absolutely convergent when $\omega_2 > \frac{3}{2}$ and $\omega +K >1$, and has no growth in the $t_\ell$, $s$ or $w$  aspects. Furthermore $S^K_3(s,w)$ converges for $\sigma, \omega$ sufficiently sufficiently positive. 

Now let
\beq
S_3^K(s,w)=\sum_{0 \leq j \leq K} \binom{w+\tfrac{k-1}{2}+j-1}{j} S_3(s,w,j)
\eeq
where 
\beq
S_3(s,w,j):=  \sum_{m< h} \frac{A(m+h)\ol{\lambda_\ell(m)}}{m^{s+\frac{k-1}{2}-j}(m+h)^{w+j}}.
\eeq
For any fixed $j$ with $0\leq j \leq K$, we can write $S_3(s,w,j)$ as 
\beq
S_3(s,w,j)=S_5(s,w,j)-S_6(s,w,j)-S_7(s,w,j)-S_8(s,w,j) \label{esses}
\eeq
when each of the Dirichlet series are absolutely convergent, with 
\beq
S_5(s,w,j):=\sum_{m=1}^\infty \sum_{h \in \mt{Z}} \frac{A(m+h)\ol{\lambda_\ell(m)}}{m^{s+\frac{k-1}{2}-j}(m+h)^{w+j}},
\eeq
\beq
S_6(s,w,j):=\sum_{m=1}^\infty \sum_{h =1}^m \frac{A(m+h)\ol{\lambda_\ell(m)}}{m^{s+\frac{k-1}{2}-j}(m+h)^{w+j}},
\eeq
\beq
S_7(s,w,j):=\sum_{m=1}^\infty \sum_{h =1}^\infty \frac{A(m-h)\ol{\lambda_\ell(m)}}{m^{s+\frac{k-1}{2}-j}(m-h)^{w+j}},
\eeq
and
\beq
S_8(s,w,j):=\sum_{m=1}^\infty \frac{A(m)\ol{\lambda_\ell(m)}}{m^{s+w+\frac{k-1}{2}}}.
\eeq

We see that 
\beq
S_6(s,w,j)=\sum_{m=1}^\infty \frac{\ol{\lambda_\ell(m)}}{m^{s+w+\frac{k-1}{2}}}\sum_{h=1}^m A(m+h)\lt(1+\tfrac{h}{m}\rt)^{-(w+j)} \ll \sum_{m=1}^\infty \frac{1}{m^{s+w+\frac{k-1}{2}-1}}
\eeq
which is absolutely convergent when $\re w_2 > \frac{3}{2}$ and has no growth in the $t_\ell$, $s$ or $w$  aspects. Since
\beq
S_7(s,w,j)=\sum_{m=1}^\infty \sum_{h=1}^{m-1} \frac{A(h)\ol{\lambda_\ell(m)}}{m^{s+\frac{k-1}{2}-j}h^{w+j}}=\sum_{m=1}^\infty \frac{\ol{\lambda_\ell(m)}}{m^{s+\frac{k-1}{2}-j}}\sum_{h=1}^{m-1} \frac{A(h)}{h^{w+j}},
\eeq
and using the bound like the one from \eqref{trivialbd} we see that $S_7(s,w,j)$ converges absolutely when $\re s+\frac{k-1}{2}-j>1$ and $\re w_2 > \frac{3}{2}$ and also has no growth in the $t_\ell,s,$ or $w$ aspects.

Now we note that
\beq
S_5(s,w,j)=L(w+j,f)L(s+\tfrac{k-1}{2}-j,\ol{\mu_\ell}) \label{overhere}
\eeq
which is absolutely convergent when $\re s+\frac{k-1}{2}-j>1$ and $\re w +j >1$, but has an analytic continuation to all $(s,w) \in \mt{C}^2$.  We do get growth in the $w$ aspect when $\omega+j\leq 1$, which we'll note later. There is only growth in the $s$ and $t_\ell$ aspects when $\re s+ \frac{k-1}{2}-j \leq 1$, and we'll see that we're ultimately unconcerned with these cases. Finally, we see that 
\beq
S_8(s,w,j)=\frac{L(s+w+\tfrac{k-1}{2},f\otimes \ol{\mu_\ell})}{\zeta(2s+2w+k-1)}
\eeq
which is absolutely convergent for $\re w_2 > \frac{1}{2}$, where it has no growth in the $t_\ell$, $s$ or $w$  aspects, but also has a meromorphic continuation to all $(s,w)\in \mt{C}^2$ with potential poles at $w_2=\frac{\varrho}{2}-\hf$ where $\varrho$ is a nontrivial zero of the zeta function.

Putting all of this together we get that $S_3(s,w,j)$ is absolutely convergent when $\re s +\frac{k-1}{2}-j>1$ and $\re w +j >1$ and has an analytic continuation into the region where $1 \geq \re w+j$ and $\re w_2 > \frac{3}{2}.$ This gives us that $S_3^K(s,w)$ is absolutely convergent when $\re s +\frac{k-1}{2}-K>1$ and $\re w>1 $ and has an analytic continuation into the region $1 \geq \re w$, $\re s +\frac{k-1}{2}-K>1$ and $\re w_2 >\frac{3}{2}$.

Using the region of convergence for $S_4^K(s,w)$ we have that $S_2(s,w)$ has an analytic continuation into the region bound by $\re s + \frac{k-1}{2}-K>1$ and $\re w + K >1$ for every $K \in \mt{R}$. This is equivalent to the region given by $\re w_2 > \frac{3}{2}$, giving the area of convergence described by the proposition. The only notable growth comes from \eqref{overhere} and gives us the bound as noted in the proposition. 
\end{proof}

\section{Employing Bochner's Convexity Theorem}
Now we use all of the results in the previous sections to prove the following proposition, which gives a meromorphic continuation of $Z^+_{f,\ell}(s,w)$ to all $(s,w) \in \mt{C}^2$. 
\begin{proposition} \label{zealotry} 
The function $Z_{f,\ell}^+(s,w)$ has a meromorphic continuation to all of $(s,w) \in \mt{C}^2$. All of its polar lines are as described in Proposition \ref{sumswithin}. 
In the region where $\omega_2>\frac{3}{2}$ and $\sigma>1$, we have
\beq
Z_{f,\ell}^+(s,w) \ll_{A_3} (1+|\omega|)^{\pmb{P}(A_3)} \label{regi1}
\eeq
where $A_3>1+|\omega|$ and $\pmb{P}(A_3)$ is an unspecified piecewise linear polynomial of $A_3$. Everywhere else when $\sigma>\hf-A$ and $\omega_2 > \frac{3}{2}+\frac{k}{2}+\frac{\beta_1}{2}+\theta-A$ we have the bound 
\beq
Z_{f,\ell}^+(s,w) \ll_{A,A_2} [(1+|\im s|)(1+|t_\ell|)(1+|\im w | )]^{\pmb{P}(A_2)}e^{\pi |t_\ell|\epsilon_A} \label{regi3}
\eeq
at least $\vep>0$ away from its poles, where $A_2>1+|\sigma|+|\omega|+\kt$, and $\epsilon_A \to 0$ as $A \to \infty$. 

\end{proposition}
\begin{proof}
Fix $A\gg 1$ and let $\sigma>\hf-A$ and $\omega_2 > \frac{3}{2}+\frac{k}{2}+\frac{\beta_1}{2}+\theta-A$ , then let $\wt{Z}_{f,\ell}^+(s,w)$ be defined as such:
\beq
\wt{Z}_{f,\ell}^+(s,w):=B_A(s,w)(Z_{f,\ell}^+(s,w)-P_\ell(s,w)) \label{ztilde}
\eeq
where
\beq
P_\ell(s,w):=\sum_{j \neq 0 } \sum_{0\leq r \leq A} \frac{G_\ell(\hf+it_j-r)c_{r,j,f,\ell}(1)L(w_2,\ol{u_j})}{(s-(\hf-r+it_j))}e^{(s-\hf-it_j+r)^2}
\eeq
and 
\bal \label{bs}
&B_A(s,w)= \\
& \lt(\frac{1}{(s-\hf)}\prod_{\ell=0}^{\lf A+1 \rfloor} \frac{\zeta(2s+2\ell)(s+\ell-\hf)^2(s+\ell)}{(s+\ell+1)(s+\ell+2)\cdots(s+\lf A-\hf \rf)}\rt)\lt(\prod_{\ell=0}^{\lf 2A-\hf-\kt-\frac{\beta_1}{2}-\theta \rf} (s+\omega_2-\tfrac{3}{2}+\ell) \rt) \notag
\end{align}
are well-defined meromorphic functions for all $(s,w)\in\mt{C}^2$. We also note that by \eqref{youngone}, \eqref{round2}, \eqref{usefulcompa}, and Stirling's formula we have that 
\beq
P_\ell(s,w) \ll_{A_2} [(1+|t_\ell|)(1+|s|)(1+|w|)]^{\pmb{P}(A_2)} \label{Pbound}
\eeq
and 
\beq
B_A(s,w) \ll_{A_2} [(1+|s|)(1+|w|)]^{\pmb{P}(A_2)} \label{Bbound}
\eeq
where $A_2> 1+|\sigma|+|\omega|+\kt$. 
By construction we observe that
\beq
\resi{s=\hf+it_j-r} Z_{f,\ell}^+(s,w)=\resi{s=\hf+it_j-r} P_\ell(s,w).
\eeq
Furthermore $B_A(s,w)$ is analytic and for $\sigma > \hf -A$ only has zeros that only coincide with poles of $Z^+_{f,\ell}(s)$. This gives us that $\wt{Z}_{f,\ell}^+(s,w)$ is an analytic function in $s$ and $w$ when $\sigma > \hf-A$ where $Z_{f,\ell}^+(s,w)$ is defined. By Propositions \ref{sumswithin} and \ref{sumswithout}, we have that $\wt{Z}_{f,\ell}^+(s,w)$ is defined for the region bound by $\omega_2>\sigma+1+\tkt+\frac{\beta_1}{2}+\theta$ and $\sigma > \hf -A$ as well as of the region bound by $\omega_2>\frac{3}{2}$ and $\sigma>1$. These regions overlap and have a convex hull spanning $\sigma>\hf-A$ and $\omega_2 > \frac{3}{2}+\frac{k}{2}+\frac{\beta_1}{2}+\theta-A$. By Bochner's Theorem of Analytic Continuation \cite{Bochner}, $\wt{Z}_{f,\ell}^+(s,w)$ has an analytic continuation into this region, and so $Z_{f,\ell}^+(s,w)$ can be extended meromorphically into this region as well. Taking $A$ to be arbitrarily large, we have that $Z_{f,\ell}^+(s,w)$ has a meromorphic continuation to all $(s,w) \in \mt{C}^2$. 

The bound \eqref{regi1} is just a restatement of \eqref{restatement}. Combining this and \eqref{tell1} with \eqref{Pbound} and \eqref{Bbound}, we can use the Phragm\'{e}n-Lindel\"{o}f theorem between some arbitrarily distant point with real part $(\sigma_0,\omega_0)$ (where we have at most polynomial growth in the $|\im s|$ and $|\im w|$ aspects) with distance away constrained only by the size of $A$, and some point in the region where \eqref{regi1} holds. This gives us the bound
\beq
Z_{f,\ell}^+(s,w) \ll_{A,A_2} [(1+|\im s|)(1+|t_\ell|)(1+|\im w | )]^{\pmb{P}(A_2)}e^{\pi |t_\ell|\alpha_A(\sigma,\omega)} \label{regi13}
\eeq
where $\alpha_A: \mt{R}^2 \to [0,1]$ is a function such that $0<\alpha_A<1$ when $(s,w)$ is away from the boundaries of the above region. Indeed, for any $\vep >0$ we can take $A$ to be large enough that $\alpha_A(s,w)<\vep$ at any $(s,w)\in\mt{C}^2$. Without loss of generality, this gives us \eqref{regi3}. 
\end{proof}

 \chapter{The Eisenstein Series Variant}
 \section{The Weighted Eisenstein Series}
 In this chapter we produce a meromorphic continuation of $Z^+_{f,u}(s,w)$, which is wholly analogous to $Z^+_{f,\ell}(s,w)$ but with divisor functions, $\sigma_{2u}(h)$, used in place of the Fourier coefficients of a Maass form, $\lambda_\ell(h)$. Our methods are nearly identical to those used in continuing $Z^+_{f,\ell}(s,w)$, and so most proofs in this chapter will be abbreviated and make note of changes when they arise.
   
 Consider the partially holomorphic, even weight $k>0$ Eisenstein series
 \beq
 E^{(k)}(z,s)=y^s+\sum_{(c,d)=1, c>0} \frac{y^s |cz+d|^k}{|cz+d|^{2s}(cz+d)^k}.
 \eeq
 where $k$ remains a positive, even integer. We can compute the Fourier coefficients when $|n|\neq 0$ to be 
 \beq
 \int_{0}^1 E^{(k)}(z,s)e^{-2\pi  i n x} \ dx = \frac{1}{\zeta(2s)}\int_0^1 \sum_{c>0, d\in\mt{Z}}  \frac{y^s |cz+d|^k}{|cz+d|^{2s}(cz+d)^k} e^{-2\pi i n x} \ dx
 \eeq
 letting $d=mc+r$ for $1 \leq r \leq c$ we have 
 \bal
  \int_{0}^1 E^{(k)}(z,s)e^{-2\pi  i n x} \ dx &= \frac{1}{\zeta(2s)}\sum_{c=1}^\infty c^{-2s} \sum_{r=1}^c \sum_{m \in \mt{Z}} \int_0^1 \frac{y^se^{-2\pi i n x}|z+m+\tfrac{r}{c}|^k}{|z+m+\tfrac{r}{c}|^{2s}(z+m+\tfrac{r}{c})^k} \ dx  \\
  & =\frac{1}{\zeta(2s)} \sum_{c=1}^\infty c^{-2s} \sum_{r=1}^c \sum_{m \in \mt{Z}}\int_{m+\frac{r}{c}}^{1+m+\frac{r}{c}} \frac{y^s e^{-2\pi i n (x-\frac{r}{c})} |z|^k}{|z|^{2s}z^k} \ dx \notag \\
  &= \frac{1}{\zeta(2s)}\sum_{c=1}^\infty c^{-2s} \sum_{r=1}^c e^{2\pi i n r/c} \int_{-\infty}^\infty \frac{y^se^{-2\pi i n x}(x+iy)^{\kt}(x-iy)^{\kt}}{(x^2+y^2)^{s}(x+iy)^k} \ dx \notag \\
  &= \frac{y^{1-s}\sigma_{1-2s}(|n|)}{\zeta(2s)}\int_{-\infty}^\infty e^{-2\pi i n xy}(x+i)^{-s-\kt}(x-i)^{-s+\kt}  \ dx \notag \\
  & =  \frac{(-1)^\kt y^{1-s}\sigma_{1-2s}(|n|)}{\zeta(2s)}\int_{-\infty}^\infty e^{-2\pi i n x y} (1+ix)^{-(s-\kt)}(1-ix)^{-(s+\kt)} \ dx. \notag
 \end{align}
 From 3.384(9) in \cite{GR} we can evaluate the integral to get 
 \beq
  \int_{0}^1 E^{(k)}(z,s)e^{-2\pi  i n x} \ dx =\frac{(-1)^\kt \pi^s \sigma_{1-2s}(|n|)|n|^{s-1}}{\zeta(2s)\Gamma(s+\mbox{sgn}(n) \tkt)}W_{\mbox{\tiny sgn}(n) \kt,\hf-s}(4\pi |n| y).
 \eeq
Similarly in the case when $n=0$ we get 
\bal 
  \int_{0}^1 E^{(k)}(z,s)\ dx &= y^s+\frac{y^{1-s}\zeta(2s-1)}{\zeta(2s)}\int_{-\infty}^\infty (x-i)^{-(s-\kt)}(x+i)^{-(s+\kt)} \ dx  \\
  & = y^s \!+\! \frac{(-1)^{\kt}4^{1-s}\pi \zeta(2s-1)\Gamma(2s-1) }{\zeta(2s)\Gamma(s+\tkt)\Gamma(s-\tkt)}y^{1-s} = y^s\! +\! \frac{\pi^{s-1}\zeta(2-2s)\Gamma(1-s+\kt)}{\pi^{-s}\zeta(2s)\Gamma(s+\tkt)}y^{1-s} \notag
\end{align}
From this, it is not difficult to show that $E^{*(k)}(z,s):=\pi^{-s}\zeta(2s)\Gamma(s+\tkt) E^{(k)}(z,s)$ has an analytic continuation to all $s \in \mt{C}$ and $E^{*(k)}(z,s)= E^{*(k)}(z,1-s)$. Further we have
\bal
E^{*(k)}(z,\tfrac{1}{2}+u) =&\pi^{-\hf-u}\zeta(1+2u)\Gamma(\thf+u+\tkt) y^{\hf+u} +   \pi^{-\hf+u}\zeta(1-2u)\Gamma(\thf-u+\tkt) y^{\hf-u} \notag\\
&+ \sum_{n > 0} (-1)^\kt \sigma_{-2u}(|n|)|n|^{-\hf+u} W_{\mbox{\tiny sgn}(n) \kt, u}(4\pi |n| y)e^{2\pi i n x} \notag \\ 
& +\sum_{n < 0} \frac{(-1)^\kt \sigma_{-2u}(|n|)|n|^{-\hf+u}\Gamma(\hf+u+\kt)}{\Gamma(\hf+u-\kt)} W_{\mbox{\tiny sgn}(n) \kt, u}(4\pi |n| y)e^{2\pi i n x}. \label{wke}
\end{align}
From this, our aim is to construct a meromorphic continuation of 
\beq
D^+_{f,u}(s;h):=\sum_{h=1}^\infty \frac{a(m+h)\sigma_{2u}(m)}{m^{s+\frac{k-1}{2}+u}} \label{deeplus}
\eeq
much as we did with $D_{f,\ell}^+(s;h)$. We will discuss this more in the next section.

\section{Continuing $D^+_{f,u}(s;h)$.}

The procedure for continuing $D^+_{f,u}(s;h)$ is largely unchanged from the argument which produced the continuation of $D_{f,\ell}^+(s;h)$. By replacing $\mu_{\ell,k}(z)$, as given in \eqref{maassk}, with $E^{*(k)}(z,\hf-\ol{u})$, as given in \eqref{wke}, we prompt the changes
\beq
\lambda_\ell(|h|) \to \sigma_{2\ol{u}}(|h|)|h|^{-\ol{u}}, \ \ \ \rho_{\ell,k}(\pm 1) \to \frac{(-1)^{\kt}\Gamma(\hf +\kt+\ol{u})}{\Gamma(\hf\pm \tkt+\ol{u})}, \ \ \ it_\ell \to -\ol{u} \label{thesubs}
\eeq
and the notational change $V_{f,\ell} \to V_{f,u} = y^{\kt} \ol{f(z)}E^{*(k)}(z,\hf-\ol{u})$ back in \eqref{tag1}. Accounting for these replacements the analog of Proposition \ref{prop5} follows with very few additional arguments.

\begin{proposition} \label{eisy}
Fix $\vep>0$, let $\sigma = \re s$. The series $D^+_{f,u}(s,h)$ given in \eqref{deeplus} converges absolutely when $\re s >1+|\re u|$ and has a meromorphic continuation to all $(s,u)\in\mt{C}^2$. It has  simple poles when $s=1/2+it_j-r$ for $r\geq 0$, where $t_j$ correspond to the eigenvalues of Maass forms, with residues of the form \beq
\resie{s=\hf+it_j-r}=G_u(\thf +it_j-r)d_{r,j,f,u}(h),
\eeq
where
\bal
&d_{r,j,f,u}(h):=  \frac{(-1)^r (4\pi)^{\hf+r-it_j} h^{r-it_j} \Gamma(\hf-it_j+r)\Gamma(2it_j-r)\ol{\rho_j(-h)\LA V_{f,u}, u_j \RA }}{r! \mcV \Gamma(\hf+it_j)\Gamma(\hf-it_j)} \label{dpoles}
\end{align}
and 
\beq
G_u(s):=\frac{(-1)^\kt (4\pi)^{s+\kt-1}\Gamma(s)}{\Gamma(s+\tfrac{k-1}{2}+u)\Gamma(s+\tfrac{k-1}{2}-u)}. \label{dagamma}
\eeq
These poles do not occur if $t_j = \pm u$ and $r \geq \kt$. The $d_{r,j,f,u}(h)$ satisfy the average upper bound
 \beq
 \sum_{j \neq 0} |d_{r,j,f,u}(h)|^2 \ll (1+|u|)^{k+1+\beta_2(u)+\vep}h^{1+2r+2\theta}. \label{rbound9}
 \eeq
 where
 \beq
\beta_2(u) = \lt\{
\begin{array}{ll}
2|\re u| & \mbox{ if } |\re u|>\hf \\
1 & \mbox{ otherwise}.
\end{array}
\rt. 
\eeq
There are also poles due to $\Omega_{f,\ell}(s;h)$, which is given by,
\beq
\Omega_{f,u}(s;h) :=\sum_{n=0}^{\lf \hf -\sigma \rf } \frac{(4\pi)^{1-s}(-1)^n h^{n}\sigma_{1-2s-2n}(h)\Gamma(2s+n-1)\Gamma(1-s)}{n!\zeta^*(2s+2n)\zeta^*(2-2s-2n)\Gamma(s+n)\Gamma(1-s-n)} \ol{\LA V_{f,u}, E^*(*,s+n)) \RA}. \label{omega24}
\eeq
Thus, $D_{u}^+(s;h)$ also has simple poles at $\hf -r$ for $r \in \mt{Z}_{>0}$ and not necessarily simple poles at $\frac \varrho 2 - r$ for $r \in \mt{Z}_{\geq 0}$ and $\varrho$ is any nontrivial zero of $\zeta(s)$. There are also poles when $s \in \mt{Z}_{\leq 0}$ due to the $\Gamma(s)$ in $G_u(s)$. 

When $s$ is at least $\vep>0$ away from the polar points and $ \re s < \hf-\kt-\tfrac{\beta_2(u)}{2}$, $D^+_{f,u}(s,h)$ is given by
\bal
&D^+_{f,u}(s,h)=G_u(s)\lt[ \sum_{j=1}^\infty  \frac{ \ol{\rho_j(-1)\lambda_j(h)}h^{\hf-s}\Gamma(s-\thf+it_j)\Gamma(s-\thf-it_j)\Gamma(1-s)\ol{\LA V_{f,u},u_j \RA }}{(4\pi)^{s-1} \Gamma(\hf+it_j)\Gamma(\hf-it_j)} \label{spec13}\rt. \\
&\lt. + \frac{1}{2\pi i } \int_{C_\sigma}  \frac{\sigma_{2z}(h)h^{\hf-s-z}\Gamma(s-\hf-z)\Gamma(s-\hf+z)\Gamma(1-s)\ol{\LA V_{f,u}, E^*(*,\hf+z) \RA }}{2(4\pi)^{s-1} \zeta^*(1+2z)\zeta^*(1-2z)\Gamma(\hf+z)\Gamma(\hf-z)} \ dz   +\Omega_{f,u}(s;h) \rt] \notag  \\ 
&  - \frac{(-1)^\kt \Gamma(s)\Gamma(1-s)}{\mcV\Gamma(\hf+\kt-u)\Gamma(\hf-\kt+u)} \sum_{m = 1}^{h-1} \frac{a(h-m)\sigma_{2u}(m)}{m^{s+\frac{k-1}{2}+u}} \notag 
\end{align} 
where $C_\sigma$ is as in \eqref{csigma}. 

For $c< \re s < 2+|\re u|$ for any $c$ such that $-\thf -\tkt -\frac{\beta_2(u)}{2} <c< \hf-\tkt-\frac{\beta_2(u)}{2}$, $D^+_{f,u}(s,n)$ is given by
\bal
&D^+_{f,u}(s;h)=\label{eisy2} G_u(s)\lt[ \frac{1}{2\pi i} \int\limits_{(2+|\re u|)} \frac{D^+_{f,u}(v;h)}{G_u(v)(v-s)} \ dv -  \frac{1}{2\pi i} \int\limits_{(c)} \frac{(\frac{D^+_{f,u}(v;h)}{G_u(v)}-\Omega_{f,u}(v;h))}{(v-s)} \ dv \rt.  \\
&\lt. -\sum_{r=0}^{\lf \kt+|\re u| \rf} \int_{C'_\sigma} \frac{(-1)^r(4\pi)^{\hf+z+r}\sigma_{-2z}(h)h^r\Gamma(-2z-r)\Gamma(\hf+z+r)}{r! \mcV  \zeta^*(1+2z)\zeta^*(1-2z)\Gamma(\hf+z)\Gamma(\hf-z)(\hf-z-r-s)} \ol{\LA V_{f,u}, E^*(*,\thf+z)) \RA} \ dz\notag \rt.\\
&+\Omega_{f,u}(s;h)-\sum_{j \neq 0}\sum_{r=0}^{\lf \kt +|\re u|  \rf} \frac{d_{r,j,f,u}(h)}{(\hf+it_j-r-s)}\notag\\
&+\sum_{r=0}^{\lf \hf-\kt-c+\re u \rf}\frac{(4\pi)^{\hf-u-r}\Gamma(\hf+\kt-u-r)\Gamma(2u-r)}{\Gamma(\hf-\kt+u)\Gamma(\hf+\kt-u)}\sum_{m=1}^{h-1} \frac{a(h-m)m^r\sigma_{2u}(m)}{(s+\frac{k-1}{2}+u+r)}\\
&\notag +\sum_{r=0}^{\lf \hf-\kt-c-\re u \rf}\frac{(4\pi)^{\hf+u-r}\Gamma(\hf+\kt+u-r)\Gamma(-2u-r)}{\Gamma(\hf-\kt+u)\Gamma(\hf+\kt-u)}\sum_{m=1}^{h-1} \frac{a(h-m)m^r\sigma_{-2u}(m)}{{(s+\frac{k-1}{2}-u+r)}} \notag \\
&\lt. + \frac{a(h)}{\sqrt{\pi}} \lt(\frac{\pi^{-u}\zeta(1+2u)\Gamma(\tfrac{k+1}{2}+u)}{(s+\tfrac{k-1}{2}+u)} +\frac{\pi^{u}\zeta(1-2u)\Gamma(\tfrac{k+1}{2}-u)}{(s+\tfrac{k-1}{2}-u)} \rt) \rt] \notag 
\end{align}
where when $\re v= 2+|\re u|$, $D^+_{f,u}(v;h)$ is given by (\ref{deeplus}) and when $\re v = c$, $D_u(v;h)$ is given by (\ref{spec13}). The curve $C'_\sigma$ is essentially the same as $C_\sigma$ but without the exception when  $s-\hf \in \mt{Z}_{\geq 0}$.

Let $A_4 > 1 + \kt + |\sigma|+ \beta_2(u)$, then for $\sigma<\hf-\kt-\frac{\beta_2}{2}$ and $s$ at least a distance of $\vep>0$ away from the poles at $s=\hf+it_j -r$, we have the bound
\bal
|D_{u}^+&(s;h)-G_u(s)\Omega_{f,u}(s;h) |\\
& \ll_{f,A_4,\vep} [(1+|s|)(1+|\im u|)]^{\pmb{P}(A_4)}h^{1-\sigma+|\re u|}e^{-\frac{\pi}{2}((|\im s | -|t_\ell|)-||\im s|-|t_\ell||)}. \label{aless4}
\end{align}
For $\sigma>1+|\re u|+\ep$, $D^+_{f,u}(s;h)$ satisfies the bound
\bal
D^+_{f,u}(s;h) \ll_{f,\re u, \vep} h^{(k-1)/2} \label{easy3}
\end{align}
When $c<\sigma<2$, and $s$ at least a distance of $\vep>0$ away from the poles at $s=\hf+it_j -r$, we have the bound
\beq
|D^+_{f,u}(s;h) -G_u(s)\Omega_{f,u}(s;h)| \ll_{f,c} h^{1-c+\frac{\beta_2(u)}{2}}[(1+|s|)(1+|t_\ell|)]^{\pmb{P}(A_4)}e^{\pif ( |t_\ell|+||\im s|-|\im t_\ell||) }.\label{later2}
\eeq

\end{proposition}
\begin{proof}
The proof of Proposition \ref{eisy} follows from some modifications of Proposition \ref{prop5} and the propositions that feed into that.  
 What few other changes there are, like the change in bounds in the regions of absolute convergence, owe to things like the growth of $\LA V_{f,u},u_j \RA$ and the constant term in the Fourier expansion of $E^{*(k)}$.

We reconsider the unfolding integral $\LA P_{h,Y}(*;s;\delta),V_{f,u} \RA$ for this case, which is completely analogous to \eqref{tag1} but with an additional constant term that must be subtracted off along with $ \LA  P_{h,Y}^{(2)}(*;s;\delta),1 \RA$, which is
\beq
\frac{a(h)}{\sqrt{\pi}} \int_{Y^{-1}}^Y y^{s+\frac{k-1}{2}}\lt(\pi^{-u}\zeta(1+2u)\Gamma(\tfrac{k+1}{2}+u) y^{u} +   \pi^{u}\zeta(1-2u)\Gamma(\tfrac{k+1}{2}-u) y^{-u}\rt)e^{-2\pi yh \delta} \frac{dy}{y}. 
\eeq
This expression is, for sufficiently positive $\sigma$, uniformly convergent to 
\bal
&\kappa_u(s;\delta)\label{kon} \\
&:=\frac{a(h)}{\sqrt{\pi}} \lt(\frac{\pi^{-u}\zeta(1+2u)\Gamma(\tfrac{k+1}{2}+u) \Gamma(s+\tfrac{k-1}{2}+u)}{(2\pi h \delta)^{s+\frac{k-1}{2}+u}}+\frac{\pi^{u}\zeta(1-2u)\Gamma(\tfrac{k+1}{2}-u)\Gamma(s+\tfrac{k-1}{2}-u)}{(2\pi h \delta)^{s+\frac{k-1}{2}-u}}\rt) \notag
\end{align}
as $Y \to \infty$. We see this has a meromorphic continuation to all $s \in \mt{C}$ for $\delta >0$ and vanishes as $\delta \to 0$ when $\re s < \hf -\kt-|\re u|$, thus it does not contribute to \eqref{spec13}. However we also see that
\bal
 &\frac{1}{2\pi i} \int_{(2+|\re u|)} \frac{\kappa_u(v;\delta)}{(v-s)} \ dv -  \frac{1}{2\pi i} \int_{(c)} \frac{\kappa_u(v;\delta)}{(v-s)} \ dv \\
& \ \ \ \ \ \ = \frac{a(h)}{\sqrt{\pi}} \lt(\frac{\pi^{-u}\zeta(1+2u)\Gamma(\tfrac{k+1}{2}+u)}{(\hf-\kt-u-s)} +\frac{\pi^{u}\zeta(1-2u)\Gamma(\tfrac{k+1}{2}-u)}{(\hf-\kt+u-s)} \rt) + \mc{O}_{c,u}(\delta) \notag
\end{align}
and thus we have the contribution to \eqref{eisy2}. We see that the extra term does not contribute additional poles, as the potential poles at $s = \hf -\kt \pm u$ are cancelled out by the corresponding zeros in $G_u(s)$, and the poles in $u$ due to the  $\Gamma(\tfrac{k+1}{2}\pm u)$ are cancelled by the trivial zeros of the zeta functions. 
 
For $\sigma$ sufficiently positive we have
\bal
 &\LA u_j, y^{\kt}\ol{f}E^{*(k)}(*,\ol{s})\RA \notag 
 = \iint\limits_{\Gamma \setminus \mt{H}}  u_j(z) y^\kt f(z) \ol{E^{*(k)}(z,\ol{s})} \ \frac{dxdy}{y^2} \notag 
\\ & =\pi^{-s}\Gamma(s+\tkt)\zeta(2s)\int_0^1 \int_0^\infty  y^{s+\kt-1} f(z) u_j(z) \frac{dxdy}{y} \notag \\
& =\frac{\Gamma(s+\tkt)\zeta(2s)}{\pi^s} \int_0^1  \int_0^\infty y^{s+\kt-1}\lt(\sum_{n=1}^\infty a(n)e^{2\pi i n z} \rt)\lt(\sum_{|m|\neq 0} \rho_j(m) 2\sqrt{y}K_{it_j}(2\pi |m| y)e^{2\pi i m x}\rt) \frac{dxdy}{y} \notag 
\\ & = 2\pi^{-s}\Gamma(s+\tkt)\zeta(2s) \rho_j(-1) \sum_{n=1}^\infty \frac{a(n)\lambda_j(n)}{(2\pi n)^{s+\frac{k-1}{2}}} \int_0^\infty y^{s+\frac{k-1}{2}}e^{-y} K_{it_j}(y) \frac{dy}{y} \notag 
\\
& = \frac{2\sqrt{\pi}\rho_j(-1)\Gamma(s+\tfrac{k-1}{2}+it_j)\Gamma(s+\tfrac{k-1}{2}-it_j)}{\pi^s (4\pi)^{s+\frac{k-1}{2}}} L(s,f \otimes u_j) 
\end{align}
which is analytic for all $s \in \mt{C}$, giving that 
\beq
\ol{\LA V_{f,u},u_j \RA} =  \frac{2\pi^u \rho_j(-1)\Gamma(\tkt+it_j-u)\Gamma(\tkt-it_j-u)}{ (4\pi)^{\kt-u}} L(\thf-u,f \otimes u_j)
\eeq
and similar reasoning gives
\beq
\ol{\LA V_{f,u},E(*,\thf+z) \RA}   = \frac{2\pi^{u}\Gamma(\tkt+z-u)\Gamma(\tkt-z-u)}{ (4\pi)^{\kt-u}\zeta^*(1+2z)} L\lt(\thf-u-z,f\rt)L\lt(\thf-u+z,f\rt) .
\eeq
Using Stirling's formula we have
\bal
&\ol{\LA V_{f,u},u_j \RA} \label{newbd1}\\
&\ll (1+|t_j+\im u|)^{\frac{k-1+\beta_2(u)}{2}}(1+|t_j-\im u|)^{\frac{k-1+\beta_2(u)}{2}}\log(1+|t_j|)e^{-\frac{\pi}{2}(|t_j+\im u|+|t_j-\im u|-|t_j|)} \notag
\end{align}
where $\beta_2(u)\geq 0$ corresponds to the convexity bound for $L(\hf-u,f \otimes u_j)$.
Similarly, assuming $\re z =0$ and letting $t = \im z$ we have
\bal
&\ol{\LA V_{f,u},E(*,\thf+z) \RA} \label{newbd2}\\
&\ll (1+|t+\im u|)^{\frac{k-1+\beta_2(u)}{2}}(1+|t-\im u|)^{\frac{k-1+\beta_2(u)}{2}}\log(1+|t|)e^{-\frac{\pi}{2}(|t+\im u|+|t-\im u|-|t|)}. \notag
\end{align}
Plugging these bounds into an analogue of Proposition \ref{ellind}, the convergence of the necessary sums and integrals follows as do the bounds. This completes the proof.\end{proof}

\section{Continuing $Z^+_{f,u}(s,w)$}
Following from Proposition \ref{eisy}, the analogues of Propositions \ref{sumswithin}, \ref{sumswithout}, and \ref{zealotry} follow even more easily, producing an analytic continuation of 
\beq
Z_{f,u}^+(s,w):=\sum_{m,h>0} \frac{a(m+h)\sigma_{2u}(m)}{m^{s+\frac{k-1}{2}+u}h^{w+\frac{k-1}{2}}}=\sum_{h>0} \frac{D^+_{f,u}(s,h)}{h^{w+\frac{k-1}{2}}}, \ \ \mbox{when} \  \ \re(s,w)>1.
\eeq
to to an enlarged region in $\mt{C}^3$. We provide the following analogous propositions for easy reference.
\begin{proposition}\label{sumswithin2}
The function $Z_{f,u}^+(s,w)$ has a meromorphic continuation to all $(s,w,u) \in \mt{C}^3$ in the region given by 
\begin{align}
& \lt\{\begin{array}{ll}
 \omega_2 > \sigma+1+\tkt+\beta_2(u) & \mbox{ if \ }  \sigma\geq \thf-\tkt-\frac{\beta_2(u)}{2}\\
\omega_2> \tfrac{3}{2}+\frac{\beta_2(u)}{2} & \mbox{ if \ } \sigma<\thf-\tkt-\frac{\beta_2(u)}{2}
 \end{array}\rt\} \\
& \ \ \ \ \ \ \ \ \ \ \ \ \ \quad \quad \quad \quad \ \notag \bigcup
\lt\{ \begin{array}{ll}
\sigma <\thf-\tkt-\frac{\beta_2(u)}{2} & \mbox{ if \ } \omega_2 >1+\vep \\
\sigma < \frac{\omega_2-k-\beta_2(u)-\vep}{2} & \mbox{ if } \  -\vep<\omega_2<1+\vep \\
\sigma<\omega_2-\kt-\frac{\beta_2(u)}{2}& \mbox{ if } \ \omega_2 \leq -\vep
\end{array}\rt\}
\end{align}
or, a fortiori, in the more simply defined sub-region given by $\omega_2>\sigma+1+\kt+\beta_2(u)$. 

When  in this region and  $\sigma<\hf-\kt-\frac{\beta_2(u)}{2}$, $Z_{f,u}^+(s,w)$ is given by 
\bal
&Z_{f,u}^+(s,w)= (4\pi)^{\kt}\Gamma(s)\sum_{j} \frac{(-1)^\kt \ol{\rho_j(-1)}\Gamma(s-\thf+it_j)\Gamma(s-\thf+it_j)\Gamma(1-s)\ol{\LA V_{f,u},u_j \RA }}{\Gamma(s+\tfrac{k-1}{2}+u)\Gamma(s+\tfrac{k-1}{2}-u)\Gamma(\hf+it_j)\Gamma(\hf-it_j)}L(w_2,\ol{u_j}) \notag \\
&+\frac{1}{2\pi i } \int_{C_{s,w}}  \frac{(4\pi)^{\kt}\Gamma(s)\Gamma(s-\hf-z)\Gamma(s-\hf+z)\Gamma(1-s)\ol{\LA V_{f,u}, E^*(*,\hf+z) \RA }\zeta(w_2+z)\zeta(w_2-z)}{(-1)^\kt 2\zeta^*(1-2z)\zeta^*(1+2z)\Gamma(s+\tfrac{k-1}{2}+u)\Gamma(s+\tfrac{k-1}{2}-u)\Gamma(\hf+z)\Gamma(\hf-z)}   \ dz \notag  \\ 
& + \wt{\Omega}_{f,u}(s,w) +\eta_{f,u}(s,w) \notag \\
& -\Gamma(s)\Gamma(1-s) \lt[\frac{1}{2\pi i } \int_{(-\vep)}\frac{\Gamma(-z)\Gamma\lt(w+\tfrac{k-1}{2}+z\rt)L(w+z,f)\zeta\lt(s-z+\tfrac{k-1}{2},u\rt)}{\Gamma\lt(w+\frac{k-1}{2}\rt)\Gamma(\hf-\kt+u)\Gamma(\hf+\kt-u)}\ dz \rt. \notag \\
&+ \frac{\Gamma(\frac{3}{2}-\kt-s-u)\Gamma(w_2+\kt-1+u)L(w_2-\hf+u,f)\zeta(1-2u)}{\Gamma(w+\frac{k-1}{2})\Gamma(\hf-\kt+u)\Gamma(\hf+\kt-u)} \notag\\
&\lt.+ \frac{\Gamma(\frac{3}{2}-\kt-s+u)\Gamma(w_2+\kt-1-u)L(w_2-\hf-u,f)\zeta(1+2u)}{\Gamma(w+\frac{k-1}{2})\Gamma(\hf-\kt+u)\Gamma(\hf+\kt-u)}\rt]\label{spec101}  
\end{align}
where we choose $\vep$ such that $\omega + \kt -\thf -\vep \notin \mt{Z}_{\leq 0}$, $C_{s,w}$ is as described after \eqref{nicetry}, 
\beq
\wt{\Omega}_{f,u}(s,w):=G_u(s) \sum_{h=1}^\infty \frac{\Omega_{f,u}(s;h)}{h^{w+\frac{k-1}{2}}}.
\eeq
where $\Omega_{f,u}(s;h)$ is as in \eqref{omega24} and
\bal
&\eta_{f,u}(s,w):= \\
& \ \ \ \frac{\pi^{w_2-\hf}(4\pi)^\kt \Gamma(s) \Gamma(1-s)\Gamma(s+\hf-w_2)\Gamma(s-\tfrac{3}{2}+w_2)\ol{\LA V_{f,u}, E^*(*,w_2-\thf) \RA }\psi_{(-\infty,1]}(\omega_2)}{(-1)^\kt\zeta^*(3-2w_2)\Gamma(s+\frac{k-1}{2}+u)\Gamma(s+\frac{k-1}{2}-u)\Gamma(\tfrac{3}{2}-w_2)\Gamma(w_2-\hf)^2}, \notag 
\end{align} 
where  $\psi_{(-\infty,1]}$ is the characteristic function on $(-\infty,1]$. This continuation has polar lines in $s$ wherever $D^+_{f,u}(s;h)$ has poles. In particular, when $s=\hf+it_j-r$ for $r \in \mt{Z}_{\geq 0}$, unless $t_j = \pm u$ and $r\geq \kt$, we have the residues
\beq
\resie{s=\hf + it_j -r} Z_{f,u}^+(s,w)=G_u(\thf+ it_j-r)d_{r,j,f,u}(1)L(w_2,\ol{u_j}), \label{polarlines2}
\eeq
where $d_{r,j,f,u}(1)$ is as in \eqref{dpoles}. Furthermore, we have polar lines of the form $s+w_2-\frac{3}{2} \in \mt{Z}_{\leq 0}$ with residues
\begin{align}
&\resie{w_2=\frac{3}{2}-s-r} Z_{f,u}^+(s,w)= \label{otherstupidpole2}\frac{(-1)^{r+\kt} \pi^{1-s-r}(4\pi)^\kt \Gamma(s) \Gamma(1-s)\Gamma(2s+r-1)\ol{\LA V_{f,u}, E^*(*,1-s-r) \RA }}{r! \zeta^*(2s+2r)\Gamma(s+\frac{k-1}{2}+u)\Gamma(s+\frac{k-1}{2}-u)\Gamma(s+r)\Gamma(1-s-r)^2}. 
\end{align}
 \notag

Let $A_5 > 1+|\sigma|+|\omega|+k+\beta_2(u)$ and let $\pmb{P}(A_5)$ be an unspecified piecewise linear polynomial in $A_5$. For  $\sigma < \hf-\tkt-\frac{\beta_1}{2}$ and  in the region of convergence we have
\beq
Z_{f,u}^+(s,w) \ll_{A_5} \lt[(1+|s|)(1+|w|)(1+|u|)\rt]^{\pmb{P}(A_5)}e^{\pi |\im u|} . 
 \label{tell1}
\eeq\end{proposition}
 \begin{proof}
 The proof of this is mostly identical to the proof of Proposition \ref{sumswithin}, merely substitute $D_{f,u}^+(s,w)$ for $D_{f,\ell}^+(s,w)$ and thus use Proposition \ref{eisy} instead of Proposition \ref{prop5}. The most substantial change is due to the alteration of the analog of \eqref{supermagica}, wherein $L(s-z+\frac{k-1}{2},\ol{\mu_\ell})$, which has no poles, is replaced by $\zeta(s-z+\frac{k-1}{2},u)$, which has poles at $z = s+\kt-\frac{3}{2} \pm u$. Thus shifting the line of integration to $\re z =-\vep$ contributes more terms as can be seen in the last two lines of \eqref{spec101}. These do not contribute any additional poles, however, since again $\Gamma(s+\frac{k-1}{2})L(s,f)$ is everywhere analytic and the possible poles at $u=0$ due to $\zeta(1-2u)$ and $\zeta(1+2u)$ have opposite residues.
 \end{proof}
 \begin{proposition}
 \label{sumswithout2} Let $w_2 = s+w+\frac{k}{2}-1$. The function $Z_{f,u}^+(s,w)$ has an analytic continuation to the region $\re w_2 >\frac{3}{2}+|\re u|$ and $\re s >1+|\re u|$. Furthermore, in the sub-region where $1 \geq \re w$ and $\re w_2 > \frac{3}{2}+|\re u | + \vep$ then we can let  $K \in \mt{R}$ be such that $1+\vep> \re w +K > 1$ and $Z_{f,u}^+(s,w)$ can be described as the analytic continuation of the series
\bal
Z_{f,u}^+(s,w)=&\sum_{m\geq h} \frac{a(m+h)\sigma_{2u}(m)}{m^{s+\frac{k-1}{2}+u}h^{w+\frac{k-1}{2}}}
+\sum_{0 \leq j \leq K} \binom{w+\tfrac{k-1}{2}+j-1}{j} \sum_{m< h} \frac{A(m+h)\sigma_{2u}(m)}{m^{s+\frac{k-1}{2}+u-j}(m+h)^{w+j}} \notag \\
& +\sum_{j>K} \binom{w+\tfrac{k-1}{2}+j-1}{j} \sum_{m< h} \frac{A(m+h)\sigma_{2u}(m)}{m^{s+\frac{k-1}{2}+u-j}(m+h)^{w+j}}.
\end{align}
In this region $Z_{f,u}^+(s,w)$ satisfies the bound
\beq \label{restatement}
Z_{f,\ell}^+(s,w) \ll_{A_3} (1+|w|)^{\pmb{P}(A_4)}
\eeq
where $A_4 > 1+ |\omega|+|\re u|$.
 \end{proposition}
 \begin{proof}
 The proof of this is identical to the proof of Proposition \ref{sumswithout}, only substituting  $\sigma_{2\ol{u}}(m)m^{-\ol{u}}$ for $\lambda_\ell(m)$. Since $\sigma_{2\ol{u}}(m)m^{-\ol{u}} \ll m^{|\re u|}$, the region of convergence is adjusted accordingly. The only other place where problems might occur is in adapting \eqref{overhere}, where $L(s+\frac{k-1}{2}-j,\ol{\mu_\ell})$ becomes $\zeta(s+\frac{k-1}{2}-j,u)$, which has poles. However, we observe that in the original proof we always required $\re s+\frac{k-1}{2}-j>1$ and in the modified proof we simply require $\re s+\frac{k-1}{2}-j>1+|\re u|.$
 \end{proof}
 \begin{proposition}
  \label{zealotry2} 
The function $Z_{f,u}^+(s,w)$ has a meromorphic continuation to all of $(s,w,u) \in \mt{C}^3$. All of its polar lines are as described in Proposition \ref{eisy}. 
In the region where $\omega_2>\frac{3}{2}+|\re u| $ and $\sigma>1+|\re u|$, we have
\beq
Z_{f,u}^+(s,w) \ll_{A_5} (1+|\omega|)^{\pmb{P}(A_4)} \label{regi12}
\eeq
where $A_5>1+|\omega|+|\re u|$ and $\pmb{P}(A_5)$ is an unspecified linear polynomial of $A_5$. Fixing $A \gg 1+|\re u|$, we have that everywhere else when $\sigma>\hf-A$ and $\omega_2 > \frac{3}{2}+\frac{k}{2}+\beta_2(u)-A$ we have the bound 
\beq
Z_{f,\ell}^+(s,w) \ll_{A,A_4} [(1+|\im s|)(1+|\im u|)(1+||\im w | )]^{\pmb{P}(A_4)}e^{\pi |\im u|\ep_A} \label{regi32}
\eeq
at least $\vep>0$ away from its poles, where at least $\vep>0$ away from its poles, where $A_4>1+|\sigma|+|\omega|+\kt+|\re u|$ and $\epsilon_A \to 0$ as $A \to \infty$. \end{proposition}
 \begin{proof}
 The proof of this is near identical to the proof of Proposition \ref{zealotry}, making use of Propositions \ref{sumswithin2} and \ref{sumswithout2} instead of Propositions \ref{sumswithin} and \ref{sumswithout}, respectively. One notable change is that we are continuing into three variables instead of two, but for any compact range of $\re u$ we see that we can get meromorphic continuation of $Z_{f,u}^+(s,w)$ to all $(s,w)\in \mt{C}^2$. Thus we do indeed have meromorphic continuation of $Z_{f,u}^+(s,w)$ to all $(s,w,u) \in \mt{C}^3$. 
 \end{proof}
  \chapter{Continuing The Triple Sums}

\section{The Families of $D^-$ and $Z^-$ Functions}
In this chapter, we will produce a meromorphic continuation of the triple sums
\beq
T_{f_1,f_2,f_3}^{\pm }(s_1,s_2,s_3):=\sum_{m,h,n=1}^\infty  \frac{a_1(m-h)\ol{a_2(m)}a_3(m\pm n)}{m^{s_1+\frac{3}{2}k-1}h^{s_2}n^{s_3\pm\frac{k-1}{2}}} \label{otp}
\eeq
as were given in the introduction, to all $(s_1,s_2,s_3) \in \mt{C}^3$. To do this, we must first better understand the spectral expansions of 
\beq
D^-_{f_1,f_2}(s;h): = \sum_{m=h+1}^\infty \frac{a_1(m-h)\ol{a_2(m)}}{m^{s+k-1}} \label{deeminus}
\eeq
 \beq
D^-_{f,\ell}(s;h):= \sum_{m=1}^\infty \frac{a(m-h)\ol{\lambda_\ell(m)}}{m^{s+\frac{k-1}{2}} } \label{shiftedminus}
 \eeq
 \beq
 D^-_{f,u}(s;h):=\sum_{m=1}^\infty \frac{a(m-h)\sigma_{2u}(m)}{m^{s+\frac{k-1}{2}+u} }
 \eeq
for fixed $h$.  The method for obtaining these constructions are just mild variants on the methods of Selberg \cite{Selberg2} and Good \cite{Gd1,Gd2} and follows with significantly less difficulty than was found in continuing  $D^+_{f,\ell}$ and $D^+_{f,u}$.  From these expansions we will proceed to give continuations of 
\beq
Z^{-}_{f,\ell}(s,w): = \sum_{m,h=1}^\infty  \frac{a(m-h)\ol{\lambda_\ell(m)}}{m^{s+\frac{k-1}{2}}h^{w+\frac{k-1}{2}} } = \sum_{h=1}^\infty \frac{D^-_{f,\ell}(s;h)}{h^{w+\frac{k-1}{2}}}
\eeq
and 
\beq
Z^{-}_{f,u}(s,w): = \sum_{m,h=1}^\infty   \frac{a(m-h)\sigma_{2u}(m)}{m^{s+\frac{k-1}{2}+u} h^{w\frac{k-1}{2}}}= \sum_{h=1}^\infty \frac{D^-_{f,u}(s;h)}{h^{w+\frac{k-1}{2}}},
\eeq
by methods analogous to, but significantly more straightforward than those used to continue $Z^+_{f,\ell}(s,w)$ and $Z^+_{f,u}(s,w)$. 
In the next section we will use the above constructions to build $T^\pm$ as given in \eqref{otp}. 

We begin by continuing $D^-_{f_1,f_2}(s;h)$. Consider the more traditional real analytic Poincar\'{e} series 
 \beq
 P_h(z,s)=\sum_{\gamma \in \Gamma_\infty \backslash \Gamma} (\im \gamma z)^s e^{2\pi i h \gamma z}
 \eeq
 which is absolutely and uniformly convergent for $\re s >1$ and is square integrable in $z$ on the fundamental domain. Since $V_{f_1,f_2}:=y^k\ol{f_1(z)}f_2(z) \in L^2(\Gamma \backslash \mt{H})$, we have that 
 \begin{align}
\langle P_h(*,s),V_{f_1,f_2}\rangle &= \iint\limits_{\Gamma \backslash \mt{H}} \lt( \sum_{\gamma \in \Gamma_\infty \backslash \Gamma} (\im \gamma z)^s e^{2\pi i h \gamma z}  \rt) y^kf_1(z)\ol{f_2(z)} \ \frac{dxdy}{y^2} \\
& = \int_0^\infty \int_0^1 y^{s+k-1} e^{2\pi i h z} f_1(z) \ol{f_2(z)} \frac{dxdy}{y}  \notag \\
& =\sum_{m=h+1}^\infty a_1(m-h)\ol{a_2(m)} \int_0^\infty y^{s+k-1} e^{-4\pi m y} \frac{dy}{y} \notag \\ 
& =  \frac{\Gamma(s+k-1)}{(4\pi)^{s+k-1}} \sum_{m=h+1}^\infty \frac{a_1(m-h)\ol{a_2(m)}}{m^{s+k-1}}. \notag 
\end{align}
Since $P_h$ is square integrable for sufficiently large $\sigma$, we can take its spectral expansion with respect to an orthonormal basis of Maass forms, $\mu_\ell$, as in \eqref{maass1} to get
\begin{align} \label{speceasy}
\langle P_h(*,s),V_{f_1,f_2}\rangle=& \sum_{\ell>0} \langle P_h(*,s),\mu_\ell \rangle \ol{\langle V_{f_1,f_2},\mu_\ell  \rangle} \\ 
&+\frac{1}{4\pi}  \int_{-\infty}^\infty \langle P_h(*,s),E(*,\thf+it)\rangle \ol{\langle V_{f_1,f_2}, E(*,\thf+it)\rangle} \ dt \notag,
\end{align}
where $E(z,s)$ is just the usual weight-zero Eisenstein series. We then compute 
\bal
 \langle P_h(*,s),\mu_\ell \rangle & = \sum_{|m|\neq 0}2\ol{\rho_\ell(m)} \int_0^\infty \int_0^1 y^{s-\hf} e^{2\pi i h z} K_{it_\ell}(2\pi |m| y) e^{-2\pi i mx} \frac{dxdy}{y} \label{discretealpha}\\
 &=\frac{2\ol{\rho_\ell(h)}}{\mcV (2\pi h)^{s-\hf}} \int_0^\infty  y^{s-\hf}e^{-y}K_{it_\ell}(y) \frac{dy}{y} \notag \\
 & = \frac{(4\pi)^{1-s}\ol{\rho_\ell(h)}\sqrt{\pi}}{\mcV h^{s-\hf}}\frac{\Gamma(s-\thf+it_\ell)\Gamma(s-\thf-it_\ell)}{\Gamma(s)} \notag
\end{align}
and similarly, when $t \in \mt{R}$,
\bal
 \langle P_h(*,s),E(*,\thf+it)\rangle & = \frac{2}{\zeta^*(1-2it)} \sum_{|m|\neq 0}\frac{\sigma_{2it}(m)}{|m|^{it}} \int_0^\infty \int_0^1 y^{s-\hf} e^{2\pi i h z} K_{it}(2\pi |m| y) e^{-2\pi i mx} \frac{dxdy}{y} \notag \\
 &=\frac{(4\pi)^{1-s} }{\zeta^*(1-2it)} \frac{\sigma_{2it}(h)}{h^{s-\hf+it}}\frac{\Gamma(s-\thf+it)\Gamma(s-\thf-it)}{\Gamma(s)} .\label{contalpha}
 \end{align}
 Putting this into \eqref{speceasy} we get
 \bal
 D^-_{f_1,f_2}(s;h)=\sum_{m=h+1}^\infty & \frac{a(m-h)\ol{a(m)}}{m^{s+k-1}} =\frac{(4\pi)^{k}h^{\hf-s}}{\Gamma(s+k-1)\Gamma(s)} \label{dminusspec}\\ 
  \times \lt[  \sum_{\ell>0}\rt. & \lt. \ol{\lambda_\ell(h)\rho_\ell(1)}\Gamma(s-\thf+it_\ell)\Gamma(s-\thf-it_\ell) \ol{\langle V_{f_1,f_2},\mu_\ell  \rangle} \rt. \notag \\
 &+ \lt.\frac{1}{2\pi i} \int_{(0)} \frac{\sigma_{2z}(h)h^{-z}\Gamma(s-\thf+z)\Gamma(s-\thf-z)}{2\zeta^*(1-2z)\zeta^*(1+2z)} \ol{\langle V_{f_1,f_2}, E^*(*,\thf+z)\rangle} \ dz \rt], \notag
 \end{align}
when $\sigma=\re s > \hf$, which by Stirling's formula and Theorem 3 in \cite{Watson} is absolutely convergent.  Indeed, we see the sum over eigenvalues is absolutely and uniformly convergent for all $s \in \mt{C}^2$ away from the poles at $s-\hf\pm it_\ell \in \mt{Z}_{\leq 0}$, and so we have a meromorphic continuation of the sum to all $\mt{C}$. The same can be said about the integral, although the meromorphic continuation is complicated by the poles of the integrand when $s-\hf\pm z \in \mt{Z}_{\leq 0}$. Let 
\beq
\mc{I}^{\mbox{\tiny con}}(s;h):=\frac{1}{2\pi i} \int_{(0)} \frac{\sigma_{2z}(h)h^{-z}\Gamma(s-\thf+z)\Gamma(s-\thf-z)}{2\zeta^*(1-2z)\zeta^*(1+2z)} \ol{\langle V_{f_1,f_2}, E^*(*,\thf+z)\rangle} \ dz
\eeq
for $\sigma>\hf$. Using the same argument as was used for continuing the continuous part of the spectrum of $D_{f,\ell}^+(s;h)$ we get that
\bal
\mc{I}^{\mbox{\tiny con}}(s;h):=&\frac{1}{2\pi i} \int_{C_\sigma}  \frac{\sigma_{2z}(h)h^{-z}\Gamma(s-\thf+z)\Gamma(s-\thf-z)}{2\zeta^*(1-2z)\zeta^*(1+2z)} \ol{\langle V_{f_1,f_2}, E^*(*,\thf+z)\rangle} \ dz \\
&+ \sum_{r=0}^{\lf \hf -\sigma \rf} \frac{(-1)^r\sigma_{1-2s-2r}(h)h^{s+r-\hf}\Gamma(2s+r-1)}{2r! \zeta^*(2s+2r)\zeta^*(2-2s-2r)}\ol{\langle V_{f_1,f_2}, E^*(*,s+r)\rangle}. \notag
\end{align}
where $C_\sigma$ is as in \eqref{csigma}. In the following proposition we summarize the facts above as well as state analogous facts about $D^-_{f,\ell}(s;h)$ and $D^-_{f,u}(s;h).$
\begin{proposition}\label{soobvious}
For fixed $h$, the functions $D^-_{f_1,f_2}(s;h)$,  $D^-_{f,\ell}(s;h)$, and $D^-_{f,u}(s;h)$ have meromorphic continuations to all $s \in \mt{C}$ with spectral expansions 
 \bal
 D^-_{f_1,f_2}(s;h)=&  \frac{(4\pi)^{k}h^{\hf-s}}{\Gamma(s+k-1)\Gamma(s)} \label{dminusspec2} \lt[  \sum_{\ell>0}\rt.  \lt. \ol{\lambda_\ell(h)\rho_\ell(1)}\Gamma(s-\thf+it_\ell)\Gamma(s-\thf-it_\ell) \ol{\langle V_{f_1,f_2},\mu_\ell  \rangle} \rt.  \\
 &+ \lt.\frac{1}{2\pi i} \int_{C_\sigma} \frac{\sigma_{2z}(h)h^{-z}\Gamma(s-\thf+z)\Gamma(s-\thf-z)}{2\zeta^*(1-2z)\zeta^*(1+2z)} \ol{\langle V_{f_1,f_2}, E^*(*,\thf+z)\rangle} \ dz \rt.\notag  \\
 &\lt. + \sum_{r=0}^{\lf \hf -\sigma \rf} \frac{(-1)^r\sigma_{1-2s-2r}(h)h^{s+r-\hf}\Gamma(2s+r-1)}{2r! \zeta^*(2s+2r)\zeta^*(2-2s-2r)}\ol{\langle V_{f_1,f_2}, E^*(*,s+r)\rangle} \rt],  \notag
 \end{align}
  \bal
 D^-_{f,\ell}(s;h)=&\frac{(4\pi)^{\kt}h^{\hf-s}}{\ol{\rho_{\ell,k}(1)}\Gamma(s+\tfrac{k-1}{2}+it_\ell)\Gamma(s+\tfrac{k-1}{2}-it_\ell)} \label{dellminusspec2} \\
 &\times  \lt[  \sum_{j>0}\rt.  \lt. \ol{\lambda_j(h)\rho_j(1)}\Gamma(s-\thf+it_j)\Gamma(s-\thf-it_j) \ol{\langle V_{f,\ell},u_j \rangle} \rt.  \notag \\
 &+ \lt.\frac{1}{2\pi i} \int_{C_\sigma} \frac{\sigma_{2z}(h)h^{-z}\Gamma(s-\thf+z)\Gamma(s-\thf-z)}{2\zeta^*(1-2z)\zeta^*(1+2z)} \ol{\langle V_{f,\ell}, E^*(*,\thf+z)\rangle} \ dz \rt.\notag  \\
 &\lt. + \sum_{r=0}^{\lf \hf -\sigma \rf} \frac{(-1)^r\sigma_{1-2s-2r}(h)h^{s+r-\hf}\Gamma(2s+r-1)}{2r! \zeta^*(2s+2r)\zeta^*(2-2s-2r)}\ol{\langle V_{f,\ell}, E^*(*,s+r)\rangle} \rt], \notag
 \end{align}
 and
   \bal
 D^-_{f,u}(s;h)=&\frac{(-1)^\kt (4\pi)^{\kt}h^{\hf-s}}{\Gamma(s+\tfrac{k-1}{2}+u)\Gamma(s+\tfrac{k-1}{2}-u)} \label{duminusspec2} \\
 &\times  \lt[  \sum_{j>0}\rt.  \lt. \ol{\lambda_j(h)\rho_j(1)}\Gamma(s-\thf+it_j)\Gamma(s-\thf-it_j) \ol{\langle V_{f,u},u_j  \rangle} \rt.  \notag \\
 &+ \lt.\frac{1}{2\pi i} \int_{C_\sigma} \frac{\sigma_{2z}(h)h^{-z}\Gamma(s-\thf+z)\Gamma(s-\thf-z)}{2\zeta^*(1-2z)\zeta^*(1+2z)} \ol{\langle V_{f,u}, E^*(*,\thf+z)\rangle} \ dz \rt.\notag  \\
 &\lt. + \sum_{r=0}^{\lf \hf -\sigma \rf} \frac{(-1)^r\sigma_{1-2s-2r}(h)h^{s+r-\hf}\Gamma(2s+r-1)}{2r! \zeta^*(2s+2r)\zeta^*(2-2s-2r)}\ol{\langle V_{f,u}, E^*(*,s+r)\rangle} \rt], \notag
 \end{align}
 where $V_{f,\ell}:=y^{\kt}\ol{f(z)}\mu_{\ell,k}(z)$ and $V_{f,u}:=y^{\kt}\ol{f(z)}E^{*(k)}(z,\hf-\ol{u})$.\end{proposition}
%

\begin{proof}
The construction of the spectral expansion of $D^-_{f_1,f_2}$ was demonstrated above, the same argument applies for $D^-_{f,\ell}$ and $D^-_{f,u}$ by replacing $V_{f_1,f_2}$ with $V_{f,\ell}$ and $V_{f,u}$ respectively. 
\end{proof}
As a corollary of the above proposition, we get a meromorphic continuation and expansion of the $Z^-$ functions and information about their poles.
\begin{corollary} \label{coraline}
The function $Z_{f,\ell}^-(s,w)$ has a meromorphic continuation to all $(s,w) \in \mt{C}^2$. Letting $w_2 = s+w+\kt-1$ and $\omega_2 = \re w_2$ we have that 
  \bal
 Z&^-_{f,\ell}(s,w)=\frac{(4\pi)^{\kt}}{\ol{\rho_{\ell,k}(1)}\Gamma(s+\tfrac{k-1}{2}+it_\ell)\Gamma(s+\tfrac{k-1}{2}-it_\ell)} \label{zeeteeell2} \\
 \times & \lt[  \sum_{j>0}\rt.  \lt. \ol{\rho_j(1)}\Gamma(s-\thf+it_j)\Gamma(s-\thf-it_j) \ol{\langle V_{f,\ell},u_j \rangle}L(w_2,\ol{u_j}) \rt.  \notag \\
 &+ \lt.\frac{1}{2\pi i} \int_{C_{s,w}} \frac{\Gamma(s-\thf+z)\Gamma(s-\thf-z)}{2\zeta^*(1-2z)\zeta^*(1+2z)} \ol{\langle V_{f,\ell}, E^*(*,\thf+z)\rangle} \zeta(w_2,z) \ dz \rt.\notag  \\
 &\lt. + \sum_{r=0}^{\lf \hf -\sigma \rf} \frac{(-1)^r\Gamma(2s+r-1)}{2r! \zeta^*(2s+2r)\zeta^*(2-2s-2r)}\ol{\langle V_{f,\ell}, E^*(*,s+r)\rangle}\zeta(w_2,\thf-s-r) \rt. \notag\\ 
 & \notag \lt. +\frac{ \pi^{w_2-\hf}\Gamma(s+\hf-w_2)\Gamma(s-\tfrac{3}{2}+w_2)\ol{\LA V_{f,\ell}, E^*(*,w_2-\thf) \RA }\psi_{(-\infty,1]}(\omega_2)}{\Gamma(w_2-\hf)\zeta^*(3-2w_2)} \rt]
 \end{align}
 for all $(s,w) \in \mt{C}^2$, where $C_{s,w}$ is as described after \eqref{nicetry} and $\psi_{(-\infty,1]}$ is the characteristic function on $(-\infty,1]$. Similarly, $Z_{f,u}^-(s,w)$ has a meromorphic continuation to all $(s,w,u) \in \mt{C}^3$ and 
   \bal
 Z&^-_{f,u}(s,w)=\frac{(-4\pi)^{\kt}}{\Gamma(s+\tfrac{k-1}{2}+u)\Gamma(s+\tfrac{k-1}{2}-u)} \label{zeeteeell3} \\
 \times & \lt[  \sum_{j>0}\rt.  \lt. \ol{\rho_j(1)}\Gamma(s-\thf+it_j)\Gamma(s-\thf-it_j) \ol{\langle V_{f,u},u_j \rangle}L(w_2,\ol{u_j}) \rt.  \notag \\
 &+ \lt.\frac{1}{2\pi i} \int_{C_{s,w}} \frac{\Gamma(s-\thf+z)\Gamma(s-\thf-z)}{2\zeta^*(1-2z)\zeta^*(1+2z)} \ol{\langle V_{f,u}, E^*(*,\thf+z)\rangle} \zeta(w_2,z) \ dz \rt.\notag  \\
 &\lt. + \sum_{r=0}^{\lf \hf -\sigma \rf} \frac{(-1)^r\Gamma(2s+r-1)}{r! \zeta^*(2s+2r)\zeta^*(2-2s-2r)}\ol{\langle V_{f,u}, E^*(*,s+r)\rangle}\zeta(w_2,\thf-s-r) \rt. \notag\\ 
 & \notag \lt. +\frac{ \pi^{w_2-\hf}\Gamma(s+\hf-w_2)\Gamma(s-\tfrac{3}{2}+w_2)\ol{\LA V_{f,u}, E^*(*,w_2-\thf) \RA }\psi_{(-\infty,1]}(\omega_2)}{\Gamma(w_2-\hf)\zeta^*(3-2w_2)} \rt].
 \end{align}

 From these expansions it is clear that $Z^-_{f,\ell}$ and $Z^-_{f,u}$ have simple poles when $s-\hf-it_j \in \mt{Z}_{\leq 0}$ for all $t_j$ and $s \in \hf \mt{Z}_{<0}$, and not necessarily simple poles when $s-\frac \varrho 2 \in \mt{Z}_{\leq 0}$ for any $\vho$ a zero of $\zeta^*(s)$. There are also simple poles when $s-\frac{3}{2}+w_2 \in \mt{Z}_{\leq 0}$ which have the residues
 \beq
 \resie{w_2=\frac{3}{2}-s-r}Z^-_{f,\ell}(s,w) =  \frac{ (-1)^{r} \pi^{1-s-r}\Gamma(2s+r-1)\ol{\LA V_{f,\ell}, E^*(*,1-s-r) \RA }}{r!\ol{\rho_{\ell,k}(1)}\Gamma(s+\frac{k-1}{2}+it_\ell)\Gamma(s+\frac{k-1}{2}-it_\ell)\Gamma(1-s-r)\zeta^*(2-2s-2r)}
 \eeq 
 and
  \beq
 \resie{w_2=\frac{3}{2}-s-r}Z^-_{f,u}(s,w) =  \frac{ (-1)^{r+\kt} \pi^{1-s-r}\Gamma(2s+r-1)\ol{\LA V_{f,u}, E^*(*,1-s-r) \RA }}{r!\Gamma(s+\frac{k-1}{2}+u)\Gamma(s+\frac{k-1}{2}-u)\Gamma(1-s-r)\zeta^*(2-2s-2r)}.
 \eeq
 For large fixed $A \gg 1$, we have the bounds
 \beq
 Z_{f,\ell}^-(s,w) \ll_{A,A_6} [(1+|\im s|)(1+|\im w|)(1+|t_\ell|)]^{\pmb{P}(A_6)} e^{\pif||t_\ell|-|\im s||\ep_A} \label{uglier}
 \eeq
 \beq
 Z_{f,u}^-(s,w) \ll_{A,A_7} [(1+|\im s|)(1+|\im w|)(1+|\im u|)]^{\pmb{P}(A_7)} e^{\pif||\im u|-|\im s||\ep_A} \label{ugliest}
 \eeq
at least $\vep>0$ away from its poles  when $\sigma,\omega > -A$, $A_6> 1+|\sigma|+|\omega|+\kt $ and $A_7 > 1+|\sigma|+|\omega|+\kt+|\re u|$ and $\ep_A \to 0$ as $A \to \infty$. 
\end{corollary}
\begin{proof}
This follows pretty easily from Proposition \ref{soobvious} and the meromorphic continuations of $L(w_2,\ol{u_j})$ and $\zeta(w_2,z)$. Absolute and uniform convergence of everything is fairly straightforward due to the exponential decay provided by the $\Gamma(s-\hf\pm it_j)$ and $\Gamma(s-\hf \pm z)$ functions, and the residual terms fall out as they did for $Z_{f,\ell}^+(s,w)$ in Proposition \ref{sumswithin} but with fewer complications due to convergence. 

To prove that bounds hold, we note that in the region where $\sigma,\omega>1$, both $ Z_{f,\ell}^-(s,w)$ and $ Z_{f,u}^-(s,w)$ are absolutely convergent as Dirichlet series and are bounded in this region. Furthermore, using Watson's formula for the triple inner-products and Stirling's approximation we see that 
\bal
 & Z^-_{f,\ell}(s,w) \ll_{A_6} [(1+|\im s|)(1+|\im w|)(1+|t_\ell|)]^{\pmb{P}(A_6)} e^{\pif||\im s|-|t_\ell||}
\end{align}
and indeed we can apply the same argument to the rest of the expansion in \eqref{zeeteeell2} and get the same bound. From there we can use the same methods as in Proposition \ref{zealotry} to construct $\wt{Z}^-_{f,\ell}(s,w)$, which is essentially $Z^-_{f,\ell}(s,w)$ with the poles removed and equivalent exponential growth. We can then apply the Phragm\'{e}n-Lindel\"{o}f theorem in the $w$ variable, as we did in Proposition \ref{zealotry}, to get the bound \eqref{uglier}. We can follow nearly identical reasoning to get \eqref{ugliest}. 
\end{proof}

\section{The Continuation of the Triple Sums} \label{joanna}
Let $\re s_i = \sigma_i$. We are now ready, at last, to give a meromorphic continuation of the triple sums
\bal
T_{f_1,f_2,f_3}^\pm(s_1,s_2,s_3)&=\sum_{m,h,n=1}^\infty  \frac{a_1(m-h)\ol{a_2(m)}a_3(h\pm n)}{m^{s_1+\frac{3k}{2}-1}h^{s_2}n^{s_3+\frac{k-1}{2}}} \label{overthee}\\
& =\sum_{m,h,n=1}^\infty  \frac{A_1(m-h)\ol{A_2(m)}A_3(h\pm n)(1-\frac{h}{m})^{\frac{k-1}{2}}(1\pm \frac{n}{h})^{\frac{k-1}{2}}}{m^{s_1+\kt}h^{s_2+\hf-\kt}n^{s_3+\frac{k-1}{2}}} \notag 
\end{align}
which are absolutely convergent as Dirichlet series for $\sigma_1>1-\kt$, $\sigma_2>\frac{3}{2}-\sigma_1$, and $\sigma_3>1$. We observe that for sufficiently large $\sigma_i$,
\beq
 T_{f_1,f_2,f_3}^\pm (s_1,s_2,s_3) = \sum_{h,n=1}^\infty \frac{D^-(s_1+\frac{k}{2};h)a(h\pm n)}{h^{s_2}n^{s_3+\frac{k-1}{2}}}.
\eeq
 Substituting the spectral expansion \eqref{dminusspec2} in for $D^-_{f_1,f_2}(s;h)$ and then interchanging the order of summation, we are now able to state the following proposition.
 \begin{proposition}
 The functions $T_{f_1,f_2,f_3}^\pm (s_1,s_2,s_3)$ as in \eqref{overthee} have meromorphic continuations to all $(s_1,s_2,s_3) \in \mt{C}^3$ and have the expansions
 \bal
&  T_{f_1,f_2,f_3}^\pm(s_1,s_2,s_3)=\label{tfspec1} \\
& \notag \lt. \sum_{\ell>0}\rt.  \lt. \frac{(4\pi)^{k}\ol{\rho_\ell(1)}\Gamma(s_1+\frac{k-1}{2}+it_\ell)\Gamma(s_1+\frac{k-1}{2}-it_\ell) \ol{\langle V_{f_1,f_2},\mu_\ell  \rangle}} {\Gamma(s_1+\frac{3k}{2}-1)\Gamma(s_1+\kt)} \rt. Z_{f_3,\ell}^\pm (s_1+s_2,s_3)  \\
 &+ \lt.\frac{1}{2\pi i} \int_{C_{\sigma_1+\kt}} \frac{(4\pi)^\kt\Gamma(s_1+\frac{k-1}{2}+z)\Gamma(s_1+\frac{k-1}{2}-z)\ol{\langle V_{f_1,f_2}, E^*(*,\thf+z)\rangle}}{2\Gamma(s_1+\frac{3k}{2}-1)\Gamma(s_1+\kt)\zeta^*(1-2z)\zeta^*(1+2z)}  Z_{f_3,z}^\pm (s_1+s_2,s_3)\ dz \rt.\notag  \\
 &\lt. + \sum_{r=0}^{\lf \hf -\sigma_1-\kt \rf} \frac{(-1)^r (4 \pi)^\kt\Gamma(2s_1+k+r-1)\ol{\langle V_{f_1,f_2}, E^*(*,s_1+\kt+r)\rangle}Z^\pm_{f_3,\hf-s_1-\kt-r}(s_1+s_2,s_3) }{2r! \Gamma(s_1+\frac{3k}{2}-1)\Gamma(s_1+\kt) \zeta^*(2s_1+k+2r)\zeta^*(2-2s_1-2r-k)} \rt. .  \notag
 \end{align}
where $C_{\sigma_1+\kt}$ corresponds to $C_\sigma$ in \eqref{csigma}. 

These functions have polar lines when 
\beq
\begin{array}{ll}
s_1 +\kt-\hf -it_j \in \mt{Z}_{\leq 0} & \\
s_1+s_2-\hf-it_j \in \mt{Z}_{\leq 0} & \\
s_1+\kt \in \hf\mt{Z}_{<0}^{\mbox{\tiny odd}}& \\
s_1+s_2 \in \hf \mt{Z}_{< 0}& \\
s_1 +\kt - \frac{\vho}{2} \in \mt{Z}_{\leq 0} & \\
s_1+s_2 - \frac{\vho}{2} \in \mt{Z}_{\leq 0} &  \\
2s_1+2s_2+s_3+\kt-\frac{5}{2} \in \mt{Z}_{\leq 0} &  \\
s_1+s_2 = 0 & \mbox{ (for just $T^+$)}
\end{array} 
\eeq
where $t_j$ is any eigenvalue of a Maass form (positive or negative) and $\vho$ is any zero of $\zeta^*(s)$. All polar lines are simple except possibly for the ones due to $\vho$, which have the same degree as the order of vanishing of $\zeta^*(\vho)$.

Letting $A \gg 1$ be fixed, and letting $A_8>1+|\sigma_1|+|\sigma_2|+|\sigma_3|+\kt$, we see that when the points $(s_1,s_2,s_3) \in \mt{C}^3$ are at least $\vep>0$ away from the poles and $\sigma_i>-A$,  then $T_{f_1,f_2,f_3}^\pm$ satisfies the upper bound
\bal
&T_{f_1,f_2,f_3}^\pm (s_1,s_2,s_3)\ll_{A,A_8,f_1,f_2,f_3,\vep}  [(1+|\im s_1|)(1+|\im s_2|)(1+|\im s_3|)]^{\pmb{P}(A_8)}e^{\epsilon_A(|\im s_1|+|\im s_2|)}. \label{lastly}
\end{align}
where $\pmb{P}(A_8)$ is a fixed but unspecified piecewise linear polynomial of $A_8$ and $\ep_A \to 0$ as $A \to \infty$.
\end{proposition}
\begin{proof}
From the bounds \eqref{regi1} and \eqref{regi3} for $Z_{f,\ell}^+$, \eqref{regi12} and \eqref{regi32} for $Z_{f,u}^+$, \eqref{uglier} for $Z_{f,\ell}^-$, and \eqref{ugliest} for $Z_{f,u}^-$ we see that we have a meromorphic continuations of $T_f^\pm(s_1,s_2,s_3)$ to all $(s_1,s_2,s_3)\in\mt{C}^3$. Knowing what we do about the poles of $Z_{f,\ell}^\pm$ and $Z_{f,u}^\pm$, 
 the relevant poles are easily found from the expression \eqref{tfspec1}. It's worth noting that the possible poles due to $\ol{\langle V_{f_1,f_2}, E^*(*,s+\kt+r)\rangle}$ when $f_1=f_2$ and $s+\kt+r=0$ are always cancelled by denominator terms. 

The bound \eqref{lastly} follows from the expression \eqref{tfspec1}, the bounds on $Z^\pm$ noted above, Watson's formula for the inner products, and Stirling's approximation. 
\end{proof}

\section{Asymptotic Growth of Shifted Sums}

Consider the inverse Mellin transform
\bal
&\lt( \frac{1}{2\pi i }\rt)^3 \! \! \! \iiint\limits_{(2)(\kt+\hf)(2)} T^\pm_{f_1,f_2,f_3}(s_1,s_2,s_3) \Gamma(s_1+\tkt-1)\Gamma(s_2-\tkt)\Gamma(s_3+\tfrac{k-1}{2})X^{s_1+s_2+s_3+\kt-\frac{3}{2}} \ ds_1 ds_2 ds_3\notag \\
& \ \ = \sum_{m,h,n=1}^\infty \frac{a_1(m-h)\ol{a_2(m)}a_3(h\pm n)}{m^k h^\kt } e^{-(\frac{m+h+n}{X})} \notag \\
 & \ \ =  \sum_{m,h,n=1}^\infty \frac{A_1(m-h)\ol{A_2(m)}A_3(h\pm n)(1-\frac{h}{m})^{\frac{k}{2}}(1\pm \frac{n}{h})^{\frac{k}{2}}}{\sqrt{(m-h)(m)(h\pm n)} } e^{-(\frac{m+h+n}{X})}.
\end{align}
Using Proposition \ref{joanna}, we can shift the lines of integration to better understand asymptotic growth of these sums. We start by shifting $s_1$ past the first pole of the integrand due to $\Gamma(s+\kt-1)$
\bal
&\lt( \frac{1}{2\pi i }\rt)^3 \! \! \! \iiint\limits_{(2)(\kt+\hf)(2)} T^\pm_{f_1,f_2,f_3}(s_1,s_2,s_3) \Gamma(s_1+\tkt-1)\Gamma(s_2-\tkt)\Gamma(s_3+\tfrac{k-1}{2})X^{s_1+s_2+s_3+\kt-\frac{3}{2}} \ ds_1 ds_2 ds_3\notag \\
&=\lt( \frac{1}{2\pi i }\rt)^3 \mkern-36mu  \iiint\limits_{(2)(\kt+\hf)(1-\kt-\vep)} \mkern-36mu T^\pm_{f_1,f_2,f_3}(s_1,s_2,s_3) \Gamma(s_1+\tkt-1)\Gamma(s_2-\tkt)\Gamma(s_3+\tfrac{k-1}{2})X^{s_1+s_2+s_3+\kt-\frac{3}{2}} \ ds_1 ds_2 ds_3 \notag \\
& \ \ \  + \lt( \frac{1}{2\pi i }\rt)^2   \iint\limits_{(2)(\kt+\hf)}  T^\pm_{f_1,f_2,f_3}(1-\tkt,s_2,s_3) \Gamma(s_2-\tkt)\Gamma(s_3+\tfrac{k-1}{2})X^{s_2+s_3-\hf} \ ds_2 ds_3. \label{yitze}
\end{align}
Taking the remaining triple integral, we can shift $s_3$ past the pole at $s_3 = \hf-\kt$ and further, just to the right of the pole at $s_3=\frac{5}{2}-\kt-2s_1-2s_2$ which occurs when $\sigma_3= -\hf-\kt+2\vep$, so we get
\bal
&\lt( \frac{1}{2\pi i }\rt)^3 \mkern-36mu  \iiint\limits_{(2)(\kt+\hf)(1-\kt-\vep)} \mkern-18mu T^\pm_{f_1,f_2,f_3}(s_1,s_2,s_3) \Gamma(s_1+\tkt-1)\Gamma(s_2-\tkt)\Gamma(s_3+\tfrac{k-1}{2})X^{s_1+s_2+s_3+\kt-\frac{3}{2}} ds_1 ds_2 ds_3 \notag \\
&= \lt( \frac{1}{2\pi i }\rt)^3 \mkern-72mu  \iiint\limits_{(-\hf-\kt+3\vep)(\kt+\hf)(1-\kt-\vep)} \mkern-60mu T^\pm_{f_1,f_2,f_3}(s_1,s_2,s_3) \Gamma(s_1+\tkt-1)\Gamma(s_2-\tkt)\Gamma(s_3+\tfrac{k-1}{2})X^{s_1+s_2+s_3+\kt-\frac{3}{2}} ds_1 ds_2 ds_3 \notag \\
& \ \ \ +  \lt( \frac{1}{2\pi i }\rt)^2 \mkern-36mu  \iint\limits_{(\kt+\hf)(1-\kt-\vep)} \mkern-18mu T^\pm_{f_1,f_2,f_3}(s_1,s_2,\thf-\tkt) \Gamma(s_1+\tkt-1)\Gamma(s_2-\tkt)X^{s_1+s_2-1} ds_1 ds_2. 
\end{align}
We see that the triple integral above is $\mc{O}(X^{-\hf+2\vep})$.  Considering the double integral above, we note that for $T^\pm_{f_1,f_2,f_3}(s_1,s_2,\thf-\tkt) $ the pole at $s_3=\frac{5}{2}-\kt-2s_1-2s_2$ turns into the pole at $s_1+s_2=1$. So shifting $s_2$ past this pole at $\sigma_2 = \kt+\vep$ and past the pole at $s_2=\kt$ and then continuing to shift just to the right of the poles at $s_1+s_2=\hf+it_j$, we get 
\bal
&\lt( \frac{1}{2\pi i }\rt)^2 \mkern-36mu  \iint\limits_{(\kt+\hf)(1-\kt-\vep)} \mkern-18mu T^\pm_{f_1,f_2,f_3}(s_1,s_2,\thf-\tkt) \Gamma(s_1+\tkt-1)\Gamma(s_2-\tkt)X^{s_1+s_2-1} ds_1 ds_2 \\
& \ \ \ \ = \frac{1}{2\pi i} \int\limits_{(1-\kt-\vep)} R^\pm_{f_1,f_2,f_3}(s_1) \Gamma(s_1+\tkt-1)\Gamma(1-s_1-\tkt) ds_1 \notag \\
& \ \ \ \ \ \ + \frac{1}{2\pi i} \int\limits_{(1-\kt-\vep)} T^\pm_{f_1,f_2,f_3}(s_1,\tkt,\thf-\tkt) \Gamma(s_1+\tkt-1) X^{s_1+\kt-1} ds_1 \notag \\
&\ \ \ \ \ \ \ + \lt( \frac{1}{2\pi i }\rt)^2 \mkern-36mu  \iint\limits_{(\kt-\hf+2\vep)(1-\kt-\vep)} \mkern-18mu T^\pm_{f_1,f_2,f_3}(s_1,s_2,\thf-\tkt) \Gamma(s_1+\tkt-1)\Gamma(s_2-\tkt)X^{s_1+s_2-1} ds_1 ds_2 \notag
\end{align}
where $R^\pm_{f_1,f_2,f_3}(s_1)$ is the residue of $T^\pm_{f_1,f_2,f_3}(s_1,s_2,\thf-\tkt)$ at $s_2=1-s_1$. However, if we reconsider all forms of $Z^\pm(s,w)$ given in Propositions \ref{sumswithin}, \ref{sumswithin2}, and  \ref{soobvious} then we'll note that when we fix $w=\hf-\kt$ then $w_2=s+w+\kt-\hf=s-\hf$ and so all of the poles at $w_2=\frac{3}{2}-s$ become poles when $s=1$, at which point the residues vanish. Following from this, it is easy to see that $R^\pm_{f_1,f_2,f_3}(s_1)$ is identically zero. We also see the first single integral above is a constant and the double integral is $\mc{O}(X^{-\hf+\vep})$. The other single integral above is also $\mc{O}(X^{-\hf+\vep})$ by shifting $s_1$ to $\sigma = \hf-\kt+\vep$. 

Returning to the double integral in \eqref{yitze}, we see that shifting $s_3$ past the pole at $s_3=\hf-\kt$, just to the right of the pole due to $s_3=\frac{5}{2}-\kt-2s_1-2s_2$ at $\sigma_3=-\hf-\kt$ we get 
\bal
 &\lt( \frac{1}{2\pi i }\rt)^2   \iint\limits_{(2)(\kt+\hf)}  T^\pm_{f_1,f_2,f_3}(1-\tkt,s_2,s_3) \Gamma(s_2-\tkt)\Gamma(s_3+\tfrac{k-1}{2})X^{s_2+s_3-\hf} \ ds_2 ds_3 \\
 &= \lt( \frac{1}{2\pi i }\rt)^2  \mkern-36mu  \iint\limits_{(-\hf-\kt+\vep)(\kt+\hf)}  \mkern-36mu T^\pm_{f_1,f_2,f_3}(1-\tkt,s_2,s_3) \Gamma(s_2-\tkt)\Gamma(s_3+\tfrac{k-1}{2})X^{s_2+s_3-\hf} \ ds_2 ds_3 \notag\\
 & \ \ \ \  +  \frac{1}{2\pi i} \int_{(\kt+\hf)} T^\pm_{f_1,f_2,f_3}(1-\tkt,s_2,\thf-\tkt)\Gamma(s_2-\tkt) X^{s_2-\kt} ds_2. \notag
\end{align}
The remaining double integral above is $\mc{O}(X^{-\hf+\vep})$. The last remaining single integral requires a bit more attention since it appears to have a double pole at $s_2 = \kt$, one due to the $\Gamma(s_2-\kt)$ and one due to  $T^\pm_{f_1,f_2,f_3}(1-\tkt,s_2,\thf-\tkt)$ from the $s_3=\frac{5}{2}-\kt-2s_1-2s_2$ singularity. However, returning to Proposition \ref{sumswithin}, Proposition \ref{sumswithin2}, and Corollary \ref{coraline}, we actually observe that the residues of $Z^\pm(1-\kt,s_2,\thf-\tkt)$ at $s_2=\kt$ are all zero. This means that $T^\pm_{f_1,f_2,f_3}(1-\tkt,\kt,\thf-\tkt)$ is actually a special value. Thus
\bal
& \frac{1}{2\pi i} \int\limits_{(\kt+\hf)} T^\pm_{f_1,f_2,f_3}(1-\tkt,s_2,\thf-\tkt)\Gamma(s_2-\tkt) X^{s_2-\kt} ds_2 \notag \\
& \ \ \ \ =T^{\pm}_{f_1,f_2,f_3}(1-\tkt,\tkt,\thf-\tkt)+\frac{1}{2\pi i} \int\limits_{(\kt-\hf+\vep)} T^\pm_{f_1,f_2,f_3}(1-\tkt,s_2,\thf-\tkt)\Gamma(s_2-\tkt) X^{s_2-\kt} ds_2. 
\end{align}
where the remaining integral is $\mc{O}(X^{-\hf+\vep})$. Thus we have proven the main theorem:\begin{theorem} \label{main}
When $A_i(r)$ are the respective normalized coefficients of even, positive weight $k$ holomorphic cusp forms $f_i$, we have that for $X \gg 1$,
\bal
\sum_{m,h,n=1}^\infty & \frac{A_1(m-h)\ol{A_2(m)}A_3(h\pm n)(1-\frac{h}{m})^{\frac{k}{2}}(1\pm \frac{n}{h})^{\frac{k}{2}}}{\sqrt{(m-h)(m)(h\pm n)} } e^{-(\frac{m+h+n}{X})} \label{booyahs} \\
&\mkern200mu =T^{\pm}_{f_1,f_2,f_3}(1-\tkt,\tkt,\thf-\tkt)+\mc{O}_{f_1,f_2,f_3}(X^{-\hf+\vep}) \notag
\end{align}
where $\vep>0$ is arbitrarily small.
\end{theorem}


\bibliography{bibfile}	
\bibliographystyle{alpha}	

\end{document}